\documentclass[11pt]{article}
\usepackage{amssymb,amsmath,epic,eepic}
\title{
\vspace{.5ex}
\begin{center}
\LARGE {\bf First hitting time and place, monopoles and multipoles for
pseudo-processes driven by the equation
$\frac{\partial}{\partial t}=\pm\frac{\partial^N}{\partial x^N}$}
\end{center}\vspace{.21ex}}

\author{Aim\'e LACHAL
\footnote{
\mbox{Postal adress: \textsc{Institut National des
Sciences Appliqu\'ees de Lyon}}
\mbox{B\^atiment L\'eonard de Vinci, 20 avenue Albert Einstein}
\hspace*{10\textwidth} \mbox{} \mbox{69621 Villeurbanne Cedex, \textsc{France}}
\hspace*{10\textwidth} \mbox{} E-mail: {\tt aime.lachal@insa-lyon.fr}
\hspace*{10\textwidth} \mbox{} Web page:
{\tt http://maths.univ-lyon1.fr/$\mbox{}^{\sim}$lachal}
}
\\
\begin{small}
\textsl{Institut National des Sciences Appliqu\'ees de Lyon}
\end{small}
}

\date{}

\newlength{\centrage}
\divide\centrage by 2
\oddsidemargin \centrage

\newtheorem{theo}{Theorem}[section]
\newtheorem{pr}[theo]{Proposition}
\newtheorem{lm}[theo]{Lemma}
\newtheorem{co}[theo]{Corollary}
\newtheorem{de}[theo]{Definition}

\newtheorem{remk}[theo]{Remark}
\newtheorem{ex}[theo]{Example}
\newenvironment{rem}{\begin{remk}\normalfont}{\ \rule{0.5em}{0.5em}\end{remk}}
\newenvironment{exam}{\begin{ex}\normalfont}{\ \rule{0.5em}{0.5em}\end{ex}}

\newcommand{\bpr}[1]{\begin{pr}#1\end{pr}}
\newcommand{\bth}[1]{\begin{theo}#1\end{theo}}
\newcommand{\blm}[1]{\begin{lm}#1\end{lm}}
\newcommand{\bco}[1]{\begin{co}#1\end{co}}
\newcommand{\bdef}[1]{\begin{de}#1\end{de}}
\newcommand{\brem}[1]{\begin{rem}#1\end{rem}}
\newcommand{\bex}[1]{\begin{exam}#1\end{exam}}

\newcommand{\beq}{\begin{equation}}
\newcommand{\eeq}{\end{equation}}
\newcommand{\beqa}{\begin{eqnarray*}}
\newcommand{\eeqa}{\end{eqnarray*}}
\newcommand{\beqan}{\begin{eqnarray}}
\newcommand{\eeqan}{\end{eqnarray}}
\newcommand{\bitem}{\begin{itemize}}
\newcommand{\eitem}{\end{itemize}}

\newcommand{\dem}{\noindent {\sc Proof. }}
\newcommand{\qed}{
\relax\ifmmode\quad\hbox{\rlap{$\sqcap$}$\sqcup$}\else
    {\unskip\nobreak\hfil\penalty50\hskip1em\null\nobreak\hfil
    \quad\hbox{\rlap{$\sqcap$}$\sqcup$}
    \parfillskip=0em\finalhyphendemerits=0\endgraf\fi}
}
\newcommand{\fin}{\ \rule{0.5em}{0.5em}}
\newcommand{\refp}[1]{(\ref{#1})}

\numberwithin{equation}{section}

\renewcommand{\a}{\alpha}
\renewcommand{\b}{\beta}
\newcommand{\g}{\gamma}
\newcommand{\G}{\Gamma}
\renewcommand{\d}{\delta}
\newcommand{\D}{\Delta}
\newcommand{\e}{\varepsilon}

\newcommand{\z}{\zeta}
\renewcommand{\k}{\kappa}
\renewcommand{\l}{\lambda}

\newcommand{\s}{\sigma}
\renewcommand{\t}{\theta}

\newcommand{\f}{\varphi}

\renewcommand{\ge}{\geqslant}
\renewcommand{\le}{\leqslant}

\newcommand{\C}{\ensuremath{\mathbb{C}}}
\newcommand{\E}{\ensuremath{\mathbb{E}}}
\newcommand{\N}{\ensuremath{\mathbb{N}}}
\renewcommand{\P}{\ensuremath{\mathbb{P}}}

\newcommand{\R}{\ensuremath{\mathbb{R}}}

\newcommand{\cC}{\mathcal{C}}
\newcommand{\cD}{\mathcal{D}}
\newcommand{\cJ}{\mathcal{J}}
\newcommand{\cK}{\mathcal{K}}
\newcommand{\cS}{\mathcal{S}}

\newcommand{\card}{\#}
\newcommand{\dis}{\displaystyle}
\newcommand{\eg}{{e.g.,}}
\newcommand{\Eo}{\E_{\,0}}
\newcommand{\Ex}{\E_{\,x}}
\newcommand{\egallaw}{\stackrel{dist}{=}}

\newcommand{\lqn}[1]{\noalign{\noindent $\displaystyle{#1}$}}

\newcommand{\ind}{1\hspace{-.27em}\mbox{\rm l}}

\newcommand{\fnk}{\f_{n,k}}
\newcommand{\Fnk}{F_{n,k}}
\newcommand{\Gnk}{G_{n,k}}

\newcommand{\tnun}{t_{n,1}}
\newcommand{\tnj}{t_{n,j}}
\newcommand{\tnk}{t_{n,k}}
\newcommand{\tnkun}{t_{n,k+1}}
\newcommand{\tauan}{\tau_{a,n}}

\begin{document}
\maketitle


\begin{abstract}
Consider the high-order heat-type equation
$\partial u/\partial t=\pm\partial^N u/\partial x^N$
for an integer $N>2$ and introduce the related Markov pseudo-process
$(X(t))_{t\ge 0}$. In this paper, we study several functionals related
to $(X(t))_{t\ge 0}$: the maximum $M(t)$ and minimum $m(t)$ up to time $t$;
the hitting times $\tau_a^+$ and $\tau_a^-$ of the half lines $(a,+\infty)$
and $(-\infty,a)$ respectively.
We provide explicit expressions for the distributions of
the vectors $(X(t),M(t))$ and $(X(t),m(t))$, as well as those of the vectors
$(\tau_a^+,X(\tau_a^+))$ and $(\tau_a^-,X(\tau_a^-))$.
\end{abstract}

\begin{footnotesize}\sc
\noindent AMS 2000 subject classifications:
\begin{tabular}[t]{l}
primary 60G20;\\ secondary 60J25. \end{tabular}  
\end{footnotesize}

{\small\tableofcontents}


\section{Introduction}

Let $N$ be an integer greater than $2$ and consider the high-order heat-type
equation
\beq\label{EDP}
\frac{\partial u}{\partial t}=\k_N\frac{\partial^{N} u}{\partial x^{N}}
\eeq
where $\k_N=(-1)^{1+N/2}$ if $N$ is even and $\k_N=\pm 1$ if $N$ is odd.
Let $p(t;z)$ be the fundamental solution of Eq.~\refp{EDP} and put
$$
p(t;x,y)=p(t;x-y).
$$
The function $p$ is characterized by its Fourier transform
\beq\label{fourier}
\int_{-\infty}^{+\infty} e^{iu\xi} p(t;\xi)\, d\xi = e^{\k_N t(-iu)^N}.
\eeq
With Eq.~\refp{EDP} one associates a Markov pseudo-process $(X(t))_{t\ge 0}$
defined on the real line and governed by a signed measure $\P$,
which is \textbf{not} a probability measure,
according to the usual rules of ordinary stochastic processes:
$$
\P_x\{X(t)\in dy\}=p(t;x,y)\,dy
$$
and for $0=t_0<t_1<\cdots<t_n$, $x_0=x$,
$$
\P_x\{X(t_1)\in dx_1,\ldots,X(t_n)\in dx_n\}=
\prod_{i=1}^{n} p(t_i-t_{i-1} ; x_{i-1}-x_i)\,dx_i.
$$
Relation~\refp{fourier} reads, by means of the expectation associated with $\P$,
$$
\Ex\!\!\left(e^{i uX(t)}\right)= e^{iux+\k_N t(iu)^N}.
$$
Such pseudo-processes have been considered by several authors,
especially in the particulary cases $N=3$ and $N=4$. The case $N=4$
is related to the biharmonic operator $\partial^4/\partial x^4$.
Few results are known in the case $N>4$.
Let us mention that for $N=2$, the pseudo-process considered here is a genuine
stochastic process (i.e., driven by a genuine probability measure),
this is the most well-known Brownian motion.

The following problems have been tackled:
\begin{itemize}
\item
Analytical study of the sample paths of that pseudo-process: Hochberg~\cite{hoch}
defined a stochastic integral (see also Motoo~\cite{motoo} in higher dimension)
and proposed an It\^o formula based on the correspondence $dx^4=dt$, he obtained
a formula for the distribution of the maximum over $[0,t]$ in the case $N=4$
with an extension to the even-order case. Noteworthy, the sample paths do not
seem to be continuous in the case $N=4$;

\item
Study of the sojourn time spent on the positive half-line up to time $t$,
$T(t)=\mbox{meas}\{s\in[0,t]:X(s)>0\}=\int_0^t\ind_{\{X(s)>0\}}\, ds$:
Krylov~\cite{kry}, Orsingher~\cite{ors}, Hochberg and Orsingher~\cite{hoch-ors1},
Nikitin and Orsingher~\cite{nik-ors}, Lachal~\cite{arcsine}
explicitly obtained the distribution of $T(t)$
(with possible conditioning on the events $\{X(t)>\mbox{(or $=$, or $<$)$0$}\}$).
Sojourn time is useful for defining local times related to the pseudo-process
$X$, see Beghin and Orsingher~\cite{bh};

\item
Study of the maximum and the minimum functionals
$$
M(t)=\max_{0\le s\le t} X(s)\quad\mbox{and}\quad m(t)=\min_{0\le s\le t} X(s):
$$
Hochberg~\cite{hoch}, Beghin et al.~\cite{bho,bor}, Lachal~\cite{arcsine}
explicitly derived the distribution of $M(t)$ and that of $m(t)$
(with possible conditioning on some values of $X(t)$);

\item
Study of the couple $(X(t),M(t))$:
Beghin et al.~\cite{ors} wrote out several formulas for the joint distribution
of $X(t)$ and $M(t)$ in the cases $N=3$ and $N=4$;

\item
Study of the first time the pseudo-process $(X(t))_{t\ge 0}$ overshoots the
level $a>0$, $\tau_a^+=\inf\{t\ge0:X(t)>a\}$:
Nishioka~\cite{nish1,nish2}, Nakajima and Sato~\cite{nak-sat} adopt a
distributional approach (in the sense of Schwartz distributions) and explicitly
obtained the joint distribution of $\tau_a^+$ and $X(\tau_a^+)$
(with possible drift) in the case $N=4$. The quantity $X(\tau_a^+)$ is the
first hitting place of the half-line $[a,+\infty)$. Nishioka~\cite{nish3}
then studied killing, reflecting and absorbing pseudo-processes;

\item
Study of the last time before becoming definitively negative
up to time $t$, $O(t)=\sup\{s\in[0,t]:X(s)>0\}$: Lachal~\cite{arcsine}
derived the distribution of $O(t)$;

\item
Study of Equation~\refp{EDP} in the case $N=4$ under other points
of view: Funaki~\cite{fun}, and next Hochberg and Orsingher~\cite{hoch-ors2} exhibited
relationships with compound processes, namely iterated Brownian motion,
Benachour et al.~\cite{vallois} provided other probabilistic interpretations.
See also the references therein.
\end{itemize}

This aim of this paper is to study the problem of the first times straddling
a fixed level $a$ (or the first hitting times of the half-lines $(a,+\infty)$
and $(-\infty,a)$):
$$
\tau_a^+=\inf\{t\ge0:X(t)>a\},\quad
\tau_a^-=\inf\{t\ge0:X(t)<a\}
$$
with the convention $\inf(\emptyset)=+\infty$.
In the spirit of the method developed by Nishioka in the case $N=4$,
we explicitly compute the joint ``signed-distributions'' (we simply shall call
``distributions'' throughout the paper for short) of the vectors
$(X(t),M(t))$ and $(X(t),m(t))$ from which we deduce those of the vectors
$(\tau_a^+,X(\tau_a^+))$ and $(\tau_a^-,X(\tau_a^-))$. The method consists of
several steps:
\bitem
\item
Defining a step-process by sampling the pseudo-process $(X(t))_{t\ge 0}$ on
dyadic times $\tnk=k/2^n$, $k\in\N$;
\item
Observing that the classical Spitzer identity holds for any signed measure,
provided the total mass equals one, and then
using this identity for deriving the distribution of
$(X(\tnk),\max_{0\le j\le k}X(\tnj))$ through its Laplace-Fourier transform
by means of that of $X(\tnk)^+$ where $x^+=\max(x,0)$;
\item
Expressing time $\tau_a^+$ (for instance) related to the sampled process
$(X(\tnk))_{k\in\N}$ by means of $(X(\tnk),\max_{0\le j\le k}X(\tnj))$;
\item
Passing to the limit when $n\to +\infty$.
\eitem

Meaningfully, we have obtained that the distributions of the hitting places
$X(\tau_a^+)$ and $X(\tau_a^-)$ are linear combinations of the
successive derivatives of the Dirac distribution $\d_a$.
In the case $N=4$, Nishioka~\cite{nish1} already found a linear combination of
$\d_a$ and $\d_a'$ and called each corresponding part ``monopole'' and
``dipole'' respectively, considering that an electric dipole having two
opposite charges $\d_{a+\e}$ and $\d_{a-\e}$ with a distance $\e$ tending to 0
may be viewed as one monopole with charge $\d_a'$. In the general case,
we shall speak of ``multipoles''.

Nishioka~\cite{nish2} used precise estimates for carrying out the rigorous
analysis of the pseudo-process corresponding to the case $N=4$.
The most important fact for providing such estimates is that
the integral of the density $p$ is absolutely convergent.
Actually, this fact holds for any even integer $N$.
When $N$ is an odd integer, the integral of $p$ is \textbf{not} absolutely
convergent and then similar estimates may not be obtained;
this makes the study of $X$ very much harder in this case.
Nevertheless, we have found, formally at least,
remarkable formulas which agree with those of Beghin et al.~\cite{bho,bor}
in the case $N=3$. They obtained them by using
a Feynman-Kac approach and solving differential equations.
We also mention some similar differential equations for any $N$.
So, we guess our formulas should hold for any odd integer $N\ge 3$.
Perhaps a distributional definition (in the sense of Schwartz distributions
since the heat-kernel is locally integrable) of the
pseudo-process $X$ might provide a properly justification to comfirm
our results. We shall not tackle this question here.

The paper is organized as follows:
in Section~\ref{sect-settings}, we write down general notations
and recall some known results.
In Section~\ref{sect-step-process}, we construct the step-process deduced from
$(X(t))_{t\ge 0}$ by sampling this latter on dyadic times.
Section~\ref{sect-distrib-(X,M)} is devoted to the distributions of the vectors
$(X(t),M(t))$ and $(X(t),m(t))$ with the aid of Spitzer identity.
Section~\ref{sect-distrib-(tau,Xtau)} deals with the distributions of the
vectors $(\tau_a^+,X(\tau_a^+))$ and $(\tau_a^-,X(\tau_a^-))$ which can be
expressed by means of those of $(X(t),M(t))$ and $(X(t),m(t))$.
Each section is completed by an illustration of the displayed results therein
to the particular cases $N\in\{2,3,4\}$.

We finally mention that the most important results have been announced,
without details, in a short Note~\cite{multipole-cras}.

\section{Settings}\label{sect-settings}


The relation $\int_{-\infty}^{+\infty} p(t;\xi)\, d\xi = 1$ holds for
all integers $N$. Moreover, if $N$ is even, the integral is absolutely
convergent (see~\cite{arcsine}) and we put
$$
\rho=\int_{-\infty}^{+\infty} |p(t;\xi)|\, d\xi>1.
$$
Notice that $\rho$ does not depend on $t$ since $p(t;\xi)=t^{-1/N}p(1;\xi/t^{1/N})$.
For odd integer $N$, the integral of $p$ is not absolutely convergent;
in this case $\rho =+\infty$.

\subsection{$N^{\mbox{\scriptsize th}}$ roots of $\k_N$}

We shall have to consider the $N^{\mbox{\scriptsize th}}$
roots of $\k_N$ ($\t_l$ for $0\le l\le N-1$ say) and distinguish the indices
$l$ such that $\Re\t_l<0$ and  $\Re\t_l>0$ (one never
has $\Re\t_l=0$). So, let us introduce the following set of indices
\beqa
J &=&  \{l\in\{0,\ldots,N-1\}:\Re\t_l>0\}, \\
K &=&  \{l\in\{0,\ldots,N-1\}:\Re\t_l<0\}.
\eeqa
We clearly have $J\cup K=\{0,\ldots,N-1\}, \quad J\cap K=\emptyset$ and
\beq\label{somme-card}
\card J+\card K=N.
\eeq

If $N=2p$, then $\k_N=(-1)^{p+1}$, $\t_l=e^{i[(2l+p+1)\pi /N]}$,
$$
J=\{p,\ldots,2p-1\}
\quad\mbox{and}\quad
K=\{0,\ldots,p-1\}.
$$
The numbers of elements of the sets $J$ and $K$ are
$$
\card J=\card K=p.
$$

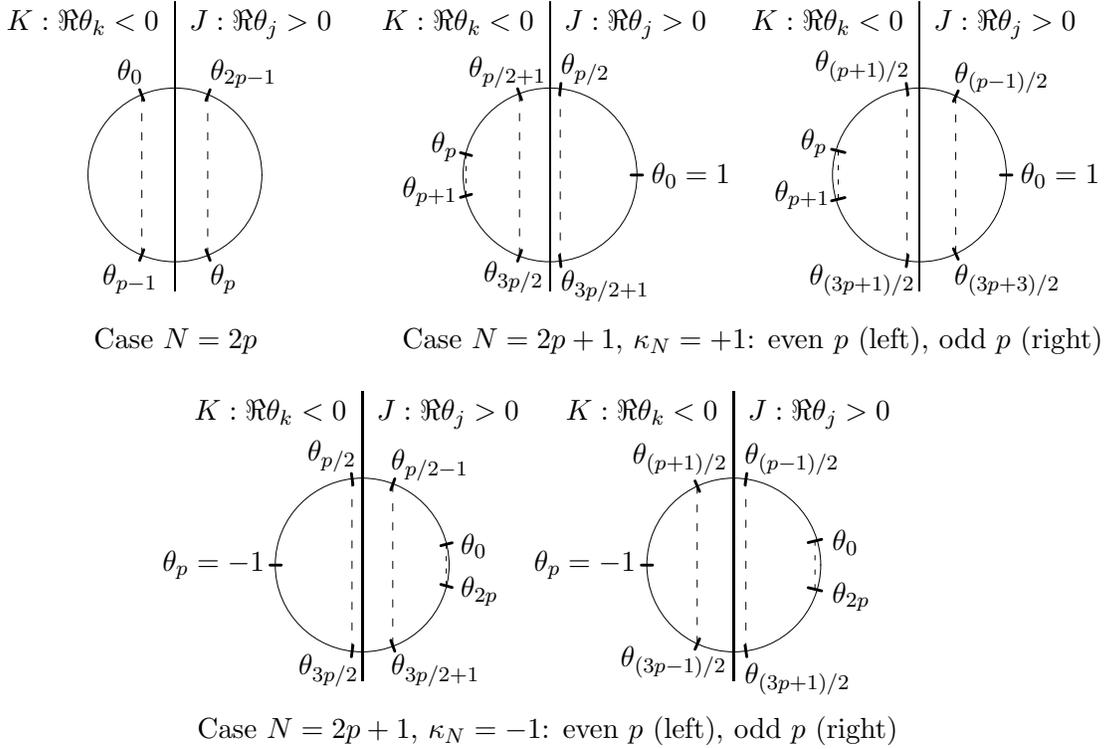
\begin{figure}[h]
\setlength{\unitlength}{1em}
\hspace*{3em}
\begin{picture}(12.7,10)
\put(3,0){\line(0,1){10}}
\put(3,4){\circle{6}}
\dashline{.3}(4.14,6.79)(4.14,1.21)
\dashline{.3}(1.86,6.79)(1.86,1.21)
\Thicklines
\drawline(4.07,6.58)(4.22,6.95)\put(4.22,7.3){$\t_{2p-1}$}
\drawline(1.93,6.58)(1.78,6.95)\put(1.,7.3){$\t_0$}
\drawline(4.07,1.42)(4.22,1.05)\put(4.22,.2){$\t_p$}
\drawline(1.93,1.42)(1.78,1.05)\put(0.5,.2){$\t_{p-1}$}
\put(-2.8,9){$K:\Re\t_k<0$}\put(3.5,9){$J:\Re\t_j>0$}
\put(0.2,-2){Case $N=2p$}
\end{picture}
\begin{picture}(12.5,10)
\put(3,0){\line(0,1){10}}
\put(3,4){\circle{6}}
\dashline{.3}(3.355,6.97)(3.355,1.03)
\dashline{.3}(1.95,6.79)(1.95,1.21)
\dashline{.2}(0.095,4.72)(0.095,3.28)
\Thicklines
\drawline(5.8,4)(6.2,4)\put(6.5,3.7){$\t_0=1$}
\drawline(3.33,6.77)(3.38,7.17)\put(3.38,7.5){$\t_{p/2}$}
\drawline(2.02,6.6)(1.88,6.98)\put(.2,7.4){$\t_{p/2+1}$}
\drawline(0.29,4.67)(-0.1,4.77)\put(-1.1,4.77){$\t_p$}
\drawline(0.29,3.33)(-0.1,3.23)\put(-2.1,3.23){$\t_{p+1}$}
\drawline(2.02,1.4)(1.88,1.02)\put(.7,.3){$\t_{3p/2}$}
\drawline(3.33,1.23)(3.38,0.83)\put(3.38,0){$\t_{3p/2+1}$}
\put(-2.8,9){$K:\Re\t_k<0$}\put(3.5,9){$J:\Re\t_j>0$}
\put(-2.1,-2){Case $N=2p+1$, $\k_N=+1$: even $p$ (left), odd $p$ (right)}
\end{picture}
\begin{picture}(10,10)
\put(3,0){\line(0,1){10}}
\put(3,4){\circle{6}}
\dashline{.3}(4.26,6.73)(4.26,1.27)
\dashline{.3}(2.58,6.97)(2.58,1.03)
\dashline{.2}(0.195,4.835)(0.195,3.165)
\Thicklines
\drawline(5.8,4)(6.2,4)\put(6.5,3.7){$\t_0=1$}
\drawline(4.18,6.55)(4.34,6.91)\put(4.2,7.2){$\t_{(p-1)/2}$}
\drawline(2.61,6.77)(2.55,7.17)\put(-.6,7.6){$\t_{(p+1)/2}$}
\drawline(0.43,4.78)(-0.04,4.89)\put(-1.1,4.9){$\t_p$}
\drawline(0.43,3.22)(-0.04,3.11)\put(-2,3.1){$\t_{p+1}$}
\drawline(2.61,1.23)(2.55,0.83)\put(-0.8,0.2){$\t_{(3p+1)/2}$}
\drawline(4.18,1.45)(4.34,1.09)\put(4.2,0.2){$\t_{(3p+3)/2}$}
\put(-2.8,9){$K:\Re\t_k<0$}\put(3.5,9){$J:\Re\t_j>0$}
\end{picture}
\vspace{5ex}

\hspace*{9.4em}
\begin{picture}(12.5,10)
\put(3,0){\line(0,1){10}}
\put(3,4){\circle{6}}
\dashline{.3}(2.645,6.97)(2.645,1.03)
\dashline{.3}(4.05,6.79)(4.05,1.21)
\dashline{.2}(5.905,4.72)(5.905,3.28)
\Thicklines
\drawline(5.71,4.67)(6.1,4.77)\put(6.4,4.4){$\t_0$}
\drawline(3.98,6.6)(4.12,6.98)\put(4,7.3){$\t_{p/2-1}$}
\drawline(2.67,6.77)(2.62,7.17)\put(1.1,7.6){$\t_{p/2}$}
\drawline(.2,4)(-.2,4)\put(-3.9,3.8){$\t_p=-1$}
\drawline(2.67,1.23)(2.62,0.83)\put(0.8,0.2){$\t_{3p/2}$}
\drawline(3.98,1.4)(4.12,1.02)\put(4,0.2){$\t_{3p/2+1}$}
\drawline(5.71,3.33)(6.1,3.23)\put(6.4,2.9){$\t_{2p}$}
\put(-2.8,9){$K:\Re\t_k<0$}\put(3.5,9){$J:\Re\t_j>0$}
\put(-2.7,-2){Case $N=2p+1$, $\k_N=-1$: even $p$ (left), odd $p$ (right)}
\end{picture}
\begin{picture}(10,10)
\put(3,0){\line(0,1){10}}
\put(3,4){\circle{6}}
\dashline{.3}(1.74,6.73)(1.74,1.27)
\dashline{.3}(3.42,6.97)(3.42,1.03)
\dashline{.2}(5.805,4.835)(5.805,3.165)
\Thicklines
\drawline(5.57,4.78)(6.04,4.89)\put(6.4,4.5){$\t_0$}
\drawline(3.39,6.77)(3.45,7.17)\put(3.38,7.5){$\t_{(p-1)/2}$}
\drawline(1.82,6.55)(1.66,6.91)\put(-0.5,7.5){$\t_{(p+1)/2}$}
\drawline(.2,4)(-.2,4)\put(-3.9,3.8){$\t_p=-1$}
\drawline(1.82,1.45)(1.66,1.09)\put(-0.9,0.5){$\t_{(3p-1)/2}$}
\drawline(3.39,1.23)(3.45,0.83)\put(3.38,-.1){$\t_{(3p+1)/2}$}
\drawline(5.57,3.22)(6.04,3.11)\put(6.4,2.7){$\t_{2p}$}
\put(-2.8,9){$K:\Re\t_k<0$}\put(3.5,9){$J:\Re\t_j>0$}
\end{picture}
\vspace{5ex}
\caption{The $N^{\mbox{\scriptsize th}}$ roots of $\k_N$}
\label{Nth-roots}
\end{figure}

If $N=2p+1$, two cases must be considered:
\bitem
\item
    For $\k_N=+1$: $\t_l=e^{i[2l\pi /N]}$ and
$$
\hspace*{-3em}\begin{array}{llll}
\dis J = \Big\{0,\ldots,\frac{p}{2}\Big\}\cup
\Big\{\frac{3p}{2}+1,\ldots,2p\Big\}
&\hspace{-.5em}\mbox{and}\hspace{-.5em}&
\dis K= \Big\{\frac{p}{2}+1,\ldots,\frac{3p}{2}\Big\}
&\mbox{if $p$ is even,}
\\[2ex]
\dis J = \Big\{0,\ldots,\frac{p-1}{2}\Big\}\cup
\Big\{\frac{3p+3}{2},\ldots,2p\Big\}
&\hspace{-.5em}\mbox{and}\hspace{-.5em}&
\dis K= \Big\{\frac{p+1}{2},\ldots,\frac{3p+1}{2}\Big\}
&\mbox{if $p$ is odd.}
\end{array}
$$
The numbers of elements of the sets $J$ and $K$ are
$$
\begin{array}{llll}
\card J=p+1   &\mbox{and}& \card K=p &\mbox{if $p$ is even,} \\
\card J=p &\mbox{and}& \card K=p+1  &\mbox{if $p$ is odd;}
\end{array}
$$
\item
For $\k_N=-1$: $\t_l=e^{i[(2l+1)\pi /N]}$ and
$$
\hspace*{-3em}\begin{array}{llll}
\dis J = \Big\{0,\ldots,\frac{p}{2}-1\Big\}\cup
\Big\{\frac{3p}{2}+1,\ldots,2p\Big\}
&\hspace{-.5em}\mbox{and}\hspace{-.5em}&
\dis K= \Big\{\frac{p}{2},\ldots,\frac{3p}{2}\Big\}
&\mbox{if $p$ is even,}
\\[2ex]
\dis J = \Big\{0,\ldots,\frac{p-1}{2}\Big\}\cup
\Big\{\frac{3p+1}{2},\ldots,2p\Big\}
&\hspace{-.5em}\mbox{and}\hspace{-.5em}&
\dis K= \Big\{\frac{p+1}{2},\ldots,\frac{3p-1}{2}\Big\}
&\mbox{if $p$ is odd.}
\end{array}
$$
The numbers of elements of the sets $J$ and $K$ are
$$
\begin{array}{llll}
\card J=p &\mbox{and}& \card K=p+1   &\mbox{if $p$ is even,} \\
\card J=p+1   &\mbox{and}& \card K=p &\mbox{if $p$ is odd.}
\end{array}
$$
\eitem
Figure~\ref{Nth-roots} illustrates the different cases.


\subsection{Recalling some known results}

We recall from~\cite{arcsine} the expressions of the kernel $p(t;\xi)$
\beq\label{fourier-inverse}
p(t;\xi)=\frac{1}{2\pi} \int_{-\infty}^{+\infty}  e^{-i\xi u+\k_N t (-iu)^N}\,du
\eeq
together with its Laplace transform (the so-called $\l$-potential of the
pseudo-process $(X(t))_{t\ge 0}$), for $\l>0$,
\beq\label{potential}
\Phi(\l;\xi)=\int_0^{+\infty} e^{-\l t}\,p(t;\xi)\,dt=
\left\{\begin{array}{ll}
\dis-\frac1N \,\l^{1/N-1} \sum_{k\in K} \t_k\,e^{\t_k\!\!\sqrt[N]{\l}\,\xi}
&\mbox{for $\xi\ge 0$,}\\[2ex]
\dis\frac1N \,\l^{1/N-1} \sum_{j\in J} \t_j\,e^{\t_j\!\!\sqrt[N]{\l}\,\xi}
&\mbox{for $\xi\le 0$.}
\end{array}\right.
\eeq
Notice that
$$
\Phi(\l;\xi)=\int_0^{+\infty} e^{-\l t}\,dt\,\P\{X(t)\in -d\xi\}/d\xi.
$$
We also recall (see the proof of Proposition~4 of~\cite{arcsine}):
\beq\label{potential-bis}
\Psi(\l;\xi)=\int_0^{+\infty} e^{-\l t}\,\P\{X(t)\le -\xi\}\,dt=
\left\{\begin{array}{ll}
\dis \frac{1}{N\l} \sum_{k\in K} e^{\t_k\!\!\sqrt[N]{\l}\,\xi}
&\mbox{for $\xi\ge 0$,}\\[3ex]
\dis\frac{1}{\l}\bigg[1-\frac1N\sum_{j\in J}
e^{\t_j\!\!\sqrt[N]{\l}\,\xi}\bigg]
&\mbox{for $\xi\le 0$.}
\end{array}\right.
\eeq
We recall the expressions of the distributions of $M(t)$ and $m(t)$ below.
\bitem
\item
Concerning the densities:
\beq\label{densities-Mm}
\begin{array}{rcl}
\dis\int_0^{+\infty} e^{-\l t}\,dt\,\P_x\{M(t)\in dz\}/dz
&=&
\dis\frac{1}{\l}\,\f_{\l}(x-z)\quad\mbox{for $x\le z$},
\\[2ex]
\dis\int_0^{+\infty} e^{-\l t}\,dt\,\P_x\{m(t)\in dz\}/dz
&=&
\dis\frac{1}{\l}\,\psi_{\l}(x-z)\quad\mbox{for $x\ge z$},
\end{array}
\eeq
with
\beq\label{phi-and-psi}
\f_{\l}(\xi)=\sqrt[N]{\l} \,\sum_{j\in J} \t_jA_j
\,e^{\t_j\!\!\sqrt[N]{\l}\,\xi},
\quad
\psi_{\l}(\xi)=-\sqrt[N]{\l} \,\sum_{k\in K} \t_kB_k
\,e^{\t_k\!\!\sqrt[N]{\l}\,\xi}
\eeq
and
$$
A_j=\prod_{l\in J\setminus \{j\}} \frac{\t_l}{\t_l-\t_j}
\quad\mbox{for $j\in J$,}\quad
B_k=\prod_{l\in K\setminus \{k\}} \frac{\t_l}{\t_l-\t_k}
\quad\mbox{for $k\in K$.}
$$
\item
Concerning the distribution functions:
\beq\label{dist-func-Mm}
\begin{array}{rcl}
\dis\int_0^{+\infty} e^{-\l t} \,\P_x(M(t)\le z) \,dt
&=&
\dis\frac{1}{\l}\Bigg[1-\sum_{j\in J} A_j
\, e^{\t_j\!\!\sqrt[N]{\l}\,(x-z)}\Bigg] \quad\mbox{for } x\le z,
\\[3ex]
\dis\int_0^{+\infty} e^{-\l t} \,\P_x(m(t)\ge z) \,dt
&=&
\dis\frac{1}{\l}\Bigg[1-\sum_{k\in K} B_k
\, e^{\t_k\!\!\sqrt[N]{\l}\,(x-z)}\Bigg] \quad\mbox{for } x\ge z.
\end{array}
\eeq
\eitem

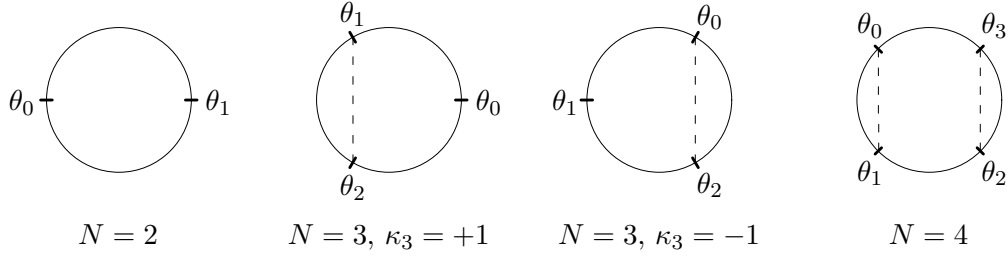
\begin{figure}[h]
\setlength{\unitlength}{1em}
\hspace*{2.2em}
\begin{picture}(9,10)
\put(3,4){\circle{5}}
\Thicklines
\drawline(.3,4)(.7,4)\put(-.8,3.65){$\t_0$}
\drawline(5.7,4)(5.3,4)\put(6,3.65){$\t_1$}
\put(1.6,-1){$N=2$}
\end{picture}
\begin{picture}(9,10)
\put(3,4){\circle{5}}
\dashline{.3}(1.75,1.83)(1.75,6.17)
\Thicklines
\drawline(5.7,4)(5.3,4)\put(6,3.65){$\t_0$}
\drawline(1.65,6.34)(1.85,6)\put(1.3,6.7){$\t_1$}
\drawline(1.65,1.66)(1.85,2)\put(1.3,.7){$\t_2$}
\put(-.5,-1){$N=3$, $\k_3=+1$}
\end{picture}
\begin{picture}(9,10)
\put(3,4){\circle{5}}
\dashline{.3}(4.25,1.83)(4.25,6.17)
\Thicklines
\drawline(4.35,6.34)(4.15,6)\put(4.3,6.6){$\t_0$}
\drawline(.3,4)(.7,4)\put(-.7,3.65){$\t_1$}
\drawline(4.35,1.66)(4.15,2)\put(4.3,.8){$\t_2$}
\put(-.5,-1){$N=3$, $\k_3=-1$}
\end{picture}
\begin{picture}(8,10)
\put(3,4){\circle{5}}
\dashline{.3}(4.76,2.24)(4.76,5.76)
\dashline{.3}(1.246,2.24)(1.24,5.76)
\Thicklines
\drawline(4.65,2.35)(4.87,2.13)\put(0.5,6.3){$\t_0$}
\drawline(1.35,2.35)(1.13,2.13)\put(0.5,1.2){$\t_1$}
\drawline(1.35,5.65)(1.13,5.87)\put(4.8,1.2){$\t_2$}
\drawline(4.65,5.65)(4.87,5.87)\put(4.8,6.3){$\t_3$}
\put(1.6,-1){$N=4$}
\end{picture}
\vspace{5ex}
\caption{The $N^{\mbox{\scriptsize th}}$ roots of $\k_N$ in the cases $N\in\{1,2,3\}$}
\label{Nth-roots-part}
\end{figure}

We explicitly write out the settings in the particular cases $N\in\{2,3,4\}$
(see~Fig.~\ref{Nth-roots-part}).

\bex{\label{case-2}
\textsl{Case $N=2$:} we have
$\k_2=+1$, $\t_0=-1,\t_1=1$, $J=\{1\},K=\{0\}$, $A_1=1,B_0=1$.
}
\bex{\label{case-3}
\textsl{Case $N=3$:} we split this (odd) case into two subcases:
\bitem
\item
for $\k_3=+1$, we have
$\t_0=1,\t_1=e^{i\,2\pi/3},\t_2=e^{-i\,2\pi/3}$,
$J=\{0\},K=\{1,2\}$,
$A_0=1,B_1=\frac{1}{1-e^{-i\,2\pi/3}}=\frac{1}{\sqrt3}\,e^{-i\,\pi/6},
B_2=\bar{B}_1=\frac{1}{\sqrt3}\,e^{i\,\pi/6}$;

\item
for $\k_3=-1$, we have
$\t_0=e^{i\,\pi/3},\t_1=-1,\t_2=e^{-i\,\pi/3}$,
$J=\{0,2\},K=\{1\}$,
$A_0=\frac{1}{1-e^{-i\,4\pi/3}}=\frac{1}{\sqrt3}\,e^{i\,\pi/6},
A_2=\bar{A}_0=\frac{1}{\sqrt3}\,e^{-i\,\pi/6}, B_1=1$.
\eitem
}
\bex{\label{case-4}
\textsl{Case $N=4$:} we have $\k_4=-1$,
$\t_0=e^{i\,3\pi/4},\t_1=e^{-i\,3\pi/4},\t_2=e^{-i\,\pi/4},\t_3=e^{i\,\pi/4}$,
$J=\{2,3\},K=\{0,1\}$,
$A_2=B_0=\frac{1}{1-e^{-i\,\pi/2}}=\frac{1}{\sqrt2}\,e^{-i\,\pi/4},
A_3=B_1=\bar{A}_2=\frac{1}{\sqrt2}\,e^{i\,\pi/4}$.
}


\subsection{Some elementary properties}\label{subsect-elem-prop}

Let us mention some elementary properties: the relation
$\prod_{l=1}^{N-1}(1-e^{i(2l\pi/N)})=N$ entails
\beq\label{prod-elem}
\prod_{l=0,l\neq m}^{N-1} \frac{\t_l}{\t_m-\t_l}=\frac 1N
\quad\mbox{for $0\le m\le N-1$.}
\eeq
The following result will be used further:
expanding into partial fractions yields, for any polynomial $P$ of
degree $\deg P\le \card J$,
\beq\label{decompo-frac}
\frac{P(x)}{\displaystyle\mathop{\textstyle\prod}_{j\in J}(1-x/\t_j)}
=\left\{\begin{array}{ll}
\dis \sum_{j\in J} \frac{A_j P(\t_j)}{1-x/\t_j}
&\mbox{if $\deg P\le \card J-1$,}
\\
\dis \sum_{j\in J} \frac{A_j P(\t_j)}{1-x/\t_j}+(-1)^{\card J}\prod_{j\in J}\t_j
&\parbox{13em}{if $\deg P= \card J$ and the highest degree coefficient of
$P$ is 1.}\end{array}\right.
\eeq
\bitem
\item
Applying~\refp{decompo-frac} to $x=0$ and $P=1$ gives
$\sum_{j\in J} A_j=\sum_{k\in K} B_k=1$. Actually, the $A_j$'s
and $B_k$'s are solutions of a Vandermonde system (see~\cite{arcsine}).

\item
Applying~\refp{decompo-frac} to $x=\t_k,k\in K$, and $P=1$ gives
\beqa\label{relation-sumBkbis}
\sum_{j\in J} \frac{\t_jA_j}{\t_j-\t_k}
=\sum_{j\in J} \frac{A_j}{1-\t_k/\t_j}
=\bigg[\prod_{j\in J} (1-\t_k/\t_j)\bigg]^{-1}
=\frac{{\displaystyle\mathop{\textstyle\prod}_{l=0,l\neq k}^{N-1}}
\frac{\t_l}{\t_l-\t_k}}
{{\displaystyle\mathop{\textstyle\prod}_{l\in K\setminus\{k\}}}\frac{\t_l}{\t_l-\t_k}}
\eeqa
which simplifies, by~\refp{prod-elem}, into (and also for the $B_k$'s)
\beq\label{somme-partial-frac}
\sum_{j\in J} \frac{\t_jA_j}{\t_j-\t_k}
=\frac{1}{NB_k}\quad\mbox{for $k\in K$ and}\quad
\sum_{k\in K} \frac{\t_kB_k}{\t_k-\t_j}=\frac{1}{NA_j}
\quad\mbox{for $j\in J$}.
\eeq

\item
Applying~\refp{decompo-frac} to $P=x^p$, $p\le \card J$, gives,
by observing that $1/\t_j=\bar{\t}_j$,
\beq\label{expansion}
\sum_{j\in J} \frac{\t_j^p A_j}{\vphantom{\bar{\bar{\t}}}1-\bar{\t}_jx}
=\left\{\begin{array}{ll}
\dis \frac{x^p}{{\displaystyle\mathop{\textstyle\prod}_{j\in J}}
\vphantom{\bar{\bar{\t}}}(1-\bar{\t}_jx)}
&\mbox{if $p\le \card J-1$,}
\\
\dis \frac{x^p}{{\displaystyle\mathop{\textstyle\prod}_{j\in J}}
\vphantom{\bar{\bar{\t}}}(1-\bar{\t}_jx)}+(-1)^{\card J-1}\prod_{j\in J}\t_j
&\mbox{if $p= \card J$.}
\end{array}\right.
\eeq
\eitem

\section{Step-process}\label{sect-step-process}

In this part, we proceed to sampling the pseudo-process $X=(X(t))_{t\ge 0}$
on the dyadic times $\tnk=k/2^n$, $k,n\in\N$ and we introduce the corresponding
step-process $X_n=(X_n(t))_{t\ge 0}$ defined for any $n\in\N$ by
$$
X_n(t)=\sum_{k=0}^{\infty} X(\tnk)\ind_{[\tnk,\tnkun)}(t).
$$
The quantity $X_n$ is a function of discrete observations of $X$ at times
$\tnk$, $k\in\N$.

For the convenience of the reader, we recall the definitions of tame functions,
functions of discrete observations, and admissible functions introduced by
Nishioka~\cite{nish2} in the case $N=4$.
%
\bdef{\label{def1}
Fix $n\in\N$. A tame function is a function of a finite number of observations of the
pseudo-process $X$ at times $\tnj$, $1\le j\le k$, that is a quantity of
the form $\Fnk=F(X(\tnun),\ldots,X(\tnk))$
for a certain $k$ and a certain bounded Borel function $F:\R^k\longrightarrow \C$.
The ``expectation'' of $\Fnk$ is defined as
$$
\Ex(\Fnk)=\int\cdots\int_{\R^k} F(x_1,\ldots,x_k)\,p(1/2^n;x,x_1)\cdots
p(1/2^n;x_{k-1},x_k)\,dx_1\cdots dx_k.
$$
}
%
We plainly have the inequality
$$
|\Ex(\Fnk)|\le \rho^k\sup_{\R^k}|F|.
$$
%
\bdef{\label{def2}
Fix $n\in\N$. A function of the discrete observations of $X$ at times $\tnk$, $k\ge 1$,
is a convergent series of tame functions:
$F_{X_n}=\sum_{k=1}^{\infty} \Fnk$ where $\Fnk$ is a tame function for all $k\ge 1$.
Assuming the series $\sum_{k=1}^{\infty} |\Ex(\Fnk)|$ convergent,
the ``expectation'' of $F_{X_n}$ is defined as
$$
\Ex(F_{X_n})=\sum_{k=1}^{\infty} \Ex(\Fnk).
$$
}
%
The definition of the expectation is consistent in the sense that it does not
depend on the representation of $F_{X_n}$ as a series (see~\cite{nish2}):
if $\sum_{k=1}^{\infty} \Fnk=\sum_{k=1}^{\infty} \Gnk$
and if the series
$\sum_{k=1}^{\infty} |\Ex(\Fnk)|$ and $\sum_{k=1}^{\infty} |\Ex(\Gnk)|$
are convergent, then
$\sum_{k=1}^{\infty} \Ex(\Fnk)=\sum_{k=1}^{\infty} \Ex(\Gnk).$
%
\bdef{\label{def3}
An admissible function is a functional $F_X$ of the pseudo-process $X$ which is
the limit of a sequence $(F_{X_n})_{n\in\N}$ of functions of discrete observations of $X$:
$$
F_X=\lim_{n\to \infty} F_{X_n},
$$
such that the sequence $(\Ex(F_{X_n}))_{n\in\N}$ is convergent.
The ``expectation'' of $F_X$ is defined as
$$
\Ex(F_X)=\lim_{n\to \infty}\Ex(F_{X_n}).
$$
}
%
This definition eludes the difficulty due to the lack of $\s$-additivity
of the signed measure~$\P$. On the other hand,
any bounded Borel function of a finite number of observations of $X$ at any times
(not necessarily dyadic) $t_1<\cdots<t_k$ is admissible and it can be seen
that, according to Definitions~\ref{def1}, \ref{def2} and \ref{def3},
\beqa
\Ex[F(X(t_1),\ldots,X(t_k))]
&=&
\int\cdots\int_{\R^k} F(x_1,\ldots,x_k)
\,p(t_1;x,x_1)\,p(t_2-t_1;x_1,x_2)\cdots
\\
&&
\hphantom{\int\cdots\int_{\R^k}}
\times p(t_k-t_{k-1};x_{k-1},x_k) \,dx_1\cdots dx_k
\eeqa
as expected in the usual sense.

For any pseudo-process $Z=(Z(t))_{t\ge 0}$, consider the functional
defined for $\l\in\C$ such that $\Re(\l)>0$, $\mu\in\R,$ $\nu>0$ by
\beq\label{functional}
F_Z(\l,\mu,\nu)=\int_0^{+\infty} e^{-\l t+i\mu H_Z(t)-\nu K_Z(t)} I_Z(t)\,dt
\eeq
where $H_Z,K_Z,I_Z$ are functionals of $Z$ defined on $[0,+\infty)$,
$K_Z$ being non negative and $I_Z$ bounded; we suppose that, for all $t\ge 0$,
$H_Z(t),K_Z(t),I_Z(t)$ are functions of the continuous observations $Z(s)$, $0\le s\le t$
(that is the observations of $Z$ up to time $t$). For $Z=X_n$, we have
\beqan
F_{X_n}(\l,\mu,\nu)
&=&
\sum_{k=0}^{\infty} \int_{\tnk}^{\tnkun}
e^{-\l t+i\mu H_{X_n}(\tnk)-\nu K_{X_n}(\tnk)} I_{X_n}(\tnk)\,dt
\nonumber\\
&=&
\sum_{k=0}^{\infty} \bigg(\int_{\tnk}^{\tnkun} e^{-\l t}\,dt\bigg)
e^{i\mu H_{X_n}(\tnk)-\nu K_{X_n}(\tnk)} I_{X_n}(\tnk)
\nonumber\\
&=&
\frac{1-e^{-\l/2^n}}{\l} \sum_{k=0}^{\infty}
e^{-\l\tnk+i\mu H_{X_n}(\tnk)-\nu K_{X_n}(\tnk)} I_{X_n}(\tnk).
\label{step-functional}
\eeqan
Since $H_{X_n}(\tnk),K_{X_n}(\tnk),I_{X_n}(\tnk)$ are functions of
$X_n(\tnj)=X(\tnj)$, $0\le j\le k$, the quantity
$e^{i\mu H_{X_n}(\tnk)-\nu K_{X_n}(\tnk)} I_{X_n}(\tnk)$
is a tame function and the series in~\refp{step-functional}
is a function of discrete observations of $X$. If the series
$$
\sum_{k=0}^{\infty} \left|\Ex\!\!\left[
e^{-\l\tnk+i\mu H_{X_n}(\tnk)-\nu K_{X_n}(\tnk)} I_{X_n}(\tnk)\right]\right|
$$
converges, the expectation of $F_{X_n}(\l,\mu,\nu)$ is defined,
according to Definition~\ref{def2}, as
$$
\Ex[F_{X_n}(\l,\mu,\nu)]=\frac{1-e^{-\l/2^n}}{\l} \sum_{k=0}^{\infty}
\Ex\!\!\left[e^{-\l\tnk+i\mu H_{X_n}(\tnk)-\nu K_{X_n}(\tnk)} I_{X_n}(\tnk)\right]\!.
$$
Finally, if $\lim_{n\to+\infty} F_{X_n}(\l,\mu,\nu)=F_X(\l,\mu,\nu)$ and
if the limit of $\Ex[F_{X_n}(\l,\mu,\nu)]$ exists as $n$ goes to $\infty$,
$F_X(\l,\mu,\nu)$ is an admissible function and its expectation
is defined, according to Definition~\ref{def3}, as
$$
\Ex[F_X(\l,\mu,\nu)]=\lim_{n\to+\infty} \Ex[F_{X_n}(\l,\mu,\nu)].
$$


\section{Distributions of $(X(t),M(t))$ and $(X(t),m(t))$}\label{sect-distrib-(X,M)}

We assume that $N$ is even.
In this section, we derive the Laplace-Fourier transforms of the vectors
$(X(t),M(t))$ and $(X(t),m(t))$ by using Spitzer identity
(Subsection~\ref{subsect-LFT-XMm}), from which we deduce the densities of
these vectors by successively inverting---three times---the Laplace-Fourier
transforms (Subsection~\ref{subsect-inverting-LFT-XMm}).
Next, we write out the formulas corresponding to the particular cases
$N\in\{2,3,4\}$ (Subsection~\ref{subsect-particular-cases}).
Finally, we compute the distribution functions of the vectors $(X(t),m(t))$
and $(X(t),M(t))$ (Subsection~\ref{subsect-dist-func-XMm}) and write out the
formulas associated with $N\in\{2,3,4\}$
(Subsection~\ref{subsect-dist-func-XMm-part}).
Although $N$ is assumed to be even, all the formulas obtained in this part when
replacing $N$ by $3$ lead to some well-known formulas in the literature.
\enlargethispage{1\baselineskip}

\subsection{Laplace-Fourier transforms}\label{subsect-LFT-XMm}

\bth{\label{theorem-LFT-XMm}
The Laplace-Fourier transform of the vectors $(X(t),M(t))$ and $(X(t),m(t))$
are given, for $\Re(\l)>0,\mu\in\R,\nu>0$, by
\beqan
\Ex\!\!\left[\int_0^{+\infty} e^{-\l t+i\mu X(t)-\nu M(t)}\,dt\right]\!\!
&=&
\frac{e^{(i\mu-\nu)x}}{\dis\mathop{\textstyle{\prod}}_{j\in J}
(\!\sqrt[N]{\l}-(i\mu-\nu)\t_j)
\dis\mathop{\textstyle{\prod}}_{k\in K} (\!\sqrt[N]{\l}-i\mu\t_k)},
\nonumber\\[-1ex]
\label{LFT-Mm}\\[-1ex]
\Ex\!\!\left[\int_0^{+\infty} e^{-\l t+i\mu X(t)+\nu m(t)}\,dt\right]\!\!
&=&
\frac{e^{(i\mu+\nu)x}}{\dis\mathop{\textstyle{\prod}}_{j\in J}
(\!\sqrt[N]{\l}-i\mu\t_j)
\dis\mathop{\textstyle{\prod}}_{k\in K} (\!\sqrt[N]{\l}-(i\mu+\nu)\t_k)}.
\nonumber
\eeqan
}
\dem
We divide the proof of Theorem~\ref{theorem-LFT-XMm} into four parts.

\noindent\textbf{$\bullet$ Step 1}

Write functionals~\refp{functional} with
$H_X(t)=X(t)$, $K_X(t)=M(t)$ or $K_X(t)=-m(t)$ and \mbox{$I_X(t)=1$:}
$$
F_X^+(\l,\mu,\nu)=\int_0^{+\infty} e^{-\l t+i\mu X(t)-\nu M(t)}\,dt
\quad\mbox{and}\quad
F_X^-(\l,\mu,\nu)=\int_0^{+\infty} e^{-\l t+i\mu X(t)+\nu m(t)}\,dt.
$$
So, putting $X_{n,k}=X(\tnk)$, $M_n(t)=\max_{0\le s\le t} X_n(s)
=\max_{0\le j\le \lfloor 2^nt\rfloor} X_{n,j}$ where $\lfloor . \rfloor$
denotes the floor function, and next $M_{n,k}=M_n(\tnk)
=\max_{0\le j\le k} X_{n,j}$, \refp{step-functional} yields, \eg\ for
$F_{X_n}^+$,
$$
F_{X_n}^+(\l,\mu,\nu)=\frac{1-e^{-\l/2^n}}{\l} \sum_{k=0}^{\infty}
e^{-\l\tnk+i\mu X_{n,k}-\nu M_{n,k}}.
$$
The functional $F_{X_n}^+(\l,\mu,\nu)$ is a function of discrete observations of $X$.
Our aim is to compute its expectation, that is to compute the expectation of the
above series and next to take the limit as $n$ goes to infinity.
For this, we observe that, using the Markov property,
\beqa
\lqn{
\left|\Ex\!\!\left[e^{-\l \tnk+i\mu X_{n,k}-\nu M_{n,k}}\right]\right|
=|e^{-\l \tnk}| \Bigg|\sum_{j=0}^k \Ex\!\!\left[e^{i\mu X_{n,k}-\nu X_{n,j}}
\ind_{\{X_{n,1}\le X_{n,j},\ldots,X_{n,k}\le X_{n,j}\}}\right]\Bigg|
}
\\
&\le&
\!\!\!(e^{-\Re(\l)/2^n})^k \Bigg|\sum_{j=0}^k\int\!\ldots\!
\int_{\{x_1\le x_j,\ldots,x_k\le x_j\}}
\hspace{-5.2em}
e^{i\mu x_k-\nu x_j} p(1/2^n;x-x_1)\cdots p(1/2^n;x_{k-1}-x_k)\,dx_1\cdots dx_k\Bigg|
\\
&\le&
\!\!\!(k+1)(\rho\,e^{-\Re(\l)/2^n})^k.
\eeqa
So, if $\Re(\l)>2^n\ln\rho$, the series
$\sum\Ex\!\!\left[e^{-\l \tnk+i\mu X_{n,k}-\nu M_{n,k}}\right]$
is absolutely convergent and then we can write the expectation of
$F_{X_n}^+(\l,\mu,\nu)$:
\beq\label{expectation-Fplus}
\Ex\!\!\left[F_{X_n}^+(\l,\mu,\nu)\right]\!=\frac{1-e^{-\l/2^n}}{\l}
\sum_{k=0}^{\infty} e^{-\l\tnk}\,\fnk^+(\mu,\nu;x)\quad\mbox{for $\Re(\l)>2^n\ln\rho$}
\eeq
with
$$
\fnk^+(\mu,\nu;x)=\Ex\!\!\left[e^{i\mu X_{n,k}-\nu M_{n,k}}\right]
=e^{(i\mu-\nu)x} \,\Eo\!\!\left[e^{-(\nu-i\mu)
M_{n,k}-i\mu(M_{n,k}-X_{n,k})}\right]\!.
$$
However, because of the domain of validity of~\refp{expectation-Fplus}, we
cannot take directly the limit as $n$ tends to infinity.
Actually, we shall see that this difficulty can be circumvented by
using sharp results on Dirichlet series.

\noindent\textbf{$\bullet$ Step 2}

Putting $z=e^{-\l/2^n}$ and noticing that $e^{-\l\tnk}=z^k$, \refp{expectation-Fplus} writes
$$
\Ex[F_{X_n}^+(\l,\mu,\nu)]=\frac{1-z}{\l} \sum_{k=0}^{\infty} \fnk^+(\mu,\nu;x)\,z^k.
$$
The generating function appearing in the last displayed equality can be
evaluated thanks to an extension of Spitzer identity, which we recall below.
%
\blm{\label{spitzer1}
Let $(\xi_k)_{k\ge 1}$ be a sequence of ``i.i.d. random variables''
and set $X_0=0$, $X_k=\sum_{j=1}^k \xi_j$ for $k\ge 1$, and
$M_k=\max_{0\le j\le k} X_j$ for $k\ge 0$.
The following relationship holds for $|z|<1$:
$$
\sum_{k=0}^{\infty} \E\!\left[e^{i\mu X_k-\nu M_k}\right]z^k
=\exp\left[\sum_{k=1}^{\infty} \E\!\left[e^{i\mu X_k-\nu X_k^+}\right]
\frac{z^k}{k} \right]\!.
$$
}
%
Observing that $1-z=\exp[\log(1-z)]=\exp[-\sum_{k=1}^{\infty} z^k/k]$,
Lemma~\ref{spitzer1} yields, for $\xi_k=X_{n,k}-X_{n,k-1}$:
\beq
\Ex[F_{X_n}^+(\l,\mu,\nu)]
=\frac{1}{\l}\,e^{(i\mu-\nu)x}\, \exp\left[\frac{1}{2^n}\sum_{k=1}^{\infty}
\frac{e^{-\l\tnk}}{\tnk}\,\psi^+(\mu,\nu;\tnk)\right]
\label{Spitzer-discret}
\eeq
where
\beqan
\psi^+(\mu,\nu;t)
&=&
\Eo\!\!\left[e^{i\mu X(t)-\nu X(t)^+}\right]\!-1
\nonumber\\
&=&
\Eo\!\!\left[\left(e^{i\mu X(t)}-1\right)\ind_{\{X(t)<0\}}\right]\!+
\Eo\!\!\left[\left(e^{(i\mu-\nu)X(t)}-1\right)\ind_{\{X(t)\ge 0\}}\right]
\nonumber\\
&=&
\int_{-\infty}^0 \left(e^{i\mu \xi}-1\right) p(t;-\xi)\,d\xi
+\int_0^{+\infty} \left(e^{(i\mu-\nu)\xi}-1\right) p(t;-\xi)\,d\xi.
\label{def-psi}
\eeqan
We plainly have $|\psi^+(\mu,\nu;t)|\le 2\rho$, and then the series
in~\refp{Spitzer-discret} defines an analytical function on the half-plane
$\{\l\in\C:\Re(\l)>0\}$. It is the analytical continuation of the function
$\l\longmapsto \Ex[F_{X_n}^+(\l,\mu,\nu)]$ which was \textit{a priori} defined
on the half-plane $\{\l\in\C:\Re(\l)>2^n\ln \rho\}$. As a byproduct,
we shall use the same notation $\Ex[F_{X_n}^+(\l,\mu,\nu)]$ for $\Re(\l)>0$.
We emphasize that the rhs of~\refp{Spitzer-discret} involves only \textbf{one}
observation of the pseudo-processus $X$ (while the lhs involves \textbf{several}
discrete observations). This important feature of Spitzer identity entails the
convergence of the series lying in~\refp{expectation-Fplus} with a lighter
constraint on the domain of validity for $\l$.

\noindent\textbf{$\bullet$ Step 3}

In order to prove that the functional $F_X^+(\l,\mu,\nu)$ is admissible, we show
that the series $\sum\Ex\!\!\left[e^{-\l \tnk+i\mu X_{n,k}-\nu M_{n,k}}\right]$
is absolutely convergent for $\Re(\l)>0$.
For this, we invoke a lemma of Bohr concerning Dirichlet series (\cite{bohr}).
Let $\sum \a_k e^{-\b_k\l}$ be a Dirichlet series of the complex variable $\l$,
where $(\a_k)_{k\in\N}$ is a sequence of complex numbers
and $(\b_k)_{k\in\N}$ is an increasing sequence of positive numbers tending to
infinity. Let us denote $\s_c$ its abscissa of convergence, $\s_a$ its abscissa
of absolute convergence and $\s_b$ the abscissa of boundedness of the analytical
continuation of its sum. If the condition $\limsup_{k\to\infty}\ln(k)/\b_k=0$
is fulfilled, then $\s_c=\s_a=\s_b$.

In our situation, we will show that the function of the variable $\l$ lying
in the rhs in~\refp{Spitzer-discret} is bounded on each half-plane $\Re(\l)\ge \e$
for any $\e>0$. We write it as
\beqa
\exp\left[\sum_{k=1}^{\infty} \psi^+(\mu,\nu;\tnk)\,\frac{e^{-\l \tnk}}{k}\right]
&=&
\exp\left[\sum_{k=1}^{\infty} \,\frac{e^{-\l \tnk}}{k}\,
\Eo\!\!\left[\left(e^{i\mu X_{n,k}}-1\right)\ind_{\{X_{n,k}<0\}}\right]\right]
\\
&&
\times\exp\left[\sum_{k=1}^{\infty} \,\frac{e^{-\l \tnk}}{k}\,
\Eo\!\!\left[\left(e^{(i\mu-\nu)X_{n,k}}-1\right)\ind_{\{X_{n,k}\ge 0\}}\right]\right]\!.
\eeqa
For any $\a\in\C$ such that $\Re(\a)\le 0$, we have
\beqa
\left|\Eo\!\!\left[\left(e^{\a X(t)}-1\right)\ind_{\{X(t)\ge 0\}}\right]\right|
&=&
\left|\Eo\!\!\left[\left(e^{\a t^{1/N}X(1)}-1\right)\ind_{\{X(1)\ge 0\}}\right]\right|
\\
&\le&
\int_0^{+\infty} \left|1-e^{\a t^{1/N}\xi}\right|\times|p(1;-\xi)|\,d\xi
\\
&\le&
2\varrho|\a| t^{1/N}
\eeqa
where we set $\varrho=\int_0^{+\infty}\xi\,|p(1;-\xi)|\,d\xi$ ($\varrho<+\infty$)
and we used the elementary inequality
$|e^{\z}-1|\le 2|\z|$ which holds for any $\z\in\C$ such that $\Re(\z)\le 0$.
Similarly,
$$
\left|\Eo\!\!\left[\left(e^{\a X(t)}-1\right)\ind_{\{X(t)< 0\}}\right]\right|
\le 2\varrho|\a| t^{1/N}.
$$
Therefore,
\beqan
\lqn{
\left|\sum_{k=1}^{\infty} \frac{e^{-\l \tnk}}{k}\,
\Eo\!\!\left[\left(e^{(\a X_{n,k}}-1\right)
\ind_{\{X_{n,k}\ge 0\mbox{ \scriptsize (or $<0$)}\}}\right]\right|
}\nonumber\\
&\le&
2\varrho|\a| \sum_{k=1}^{\infty} \frac{e^{-\Re(\l) \tnk}}{k}\,\tnk^{1/N}
=\frac{2\varrho|\a|}{2^n} \sum_{k=1}^{\infty} \frac{e^{-\Re(\l) \tnk}}{\tnk^{1-1/N}}
\nonumber\\
&\le&
2\varrho|\a| \sum_{k=1}^{\infty}
\int_{\tnk}^{\tnkun}\frac{e^{-\Re(\l)t}}{t^{1-1/N}}\,dt
\le 2\varrho|\a| \int_0^{+\infty}\frac{e^{-\Re(\l)t}}{t^{1-1/N}}\,dt
\nonumber\\
&\le&
\frac{2\,\G(1/N)\varrho|\a|}{\Re(\l)^{1/N}}.
\label{majorationXM}
\eeqan
This proves that the rhs of~\refp{Spitzer-discret} is bounded on each
half-plane $\Re(\l)\ge \e$ for any $\e>0$. So, the convergence of the series
lying in~\refp{expectation-Fplus} holds in the domain $\Re(\l)>0$ and the
functional $F_X^+(\l,\mu,\nu)$ is admissible.

\noindent\textbf{$\bullet$ Step 4}

Now, we can pass to the limit when $n\to +\infty$ in~\refp{Spitzer-discret} and we obtain
\beq\label{esp-FX}
\Ex[F_X^+(\l,\mu,\nu)] = \frac{1}{\l}\,e^{(i\mu-\nu)x}\,
\exp\left[\int_0^{+\infty} e^{-\l t}\,\psi^+(\mu,\nu;t)\,\frac{dt}{t}\right]
\quad\mbox{for $\Re(\l)>0$.}
\eeq
A similar formula holds for $F_X^-$.

From~\refp{def-psi}, we see that we need to evaluate integrals of the form
$$
\int_0^{+\infty} e^{-\l t}\,\frac{dt}{t}\int_0^{+\infty}
(e^{\a\xi}-1) p(t;-\xi)\,d\xi \quad\mbox{for $\Re(\a)\le 0$}
$$
and
$$
\int_0^{+\infty} e^{-\l t}\,\frac{dt}{t}\int_{-\infty}^0
(e^{\a\xi}-1) p(t;-\xi)\,d\xi \quad\mbox{for $\Re(\a)\ge 0$.}
$$
We have, for $\Re(\a)\le 0$,
\beqan
\lqn{
\int_0^{+\infty} e^{-\l t}\,\frac{dt}{t}\int_0^{+\infty}
(e^{\a \xi}-1) p(t;-\xi)\,d\xi}\nonumber\\
&=&
\int_0^{+\infty} dt \int_{\l}^{+\infty} e^{-ts}\,ds
\int_0^{+\infty} (e^{\a \xi}-1) p(t;-\xi)\,d\xi
\nonumber\\
&=&
\int_{\l}^{+\infty} ds \int_0^{+\infty} (e^{\a \xi}-1) \,d\xi
\int_0^{+\infty} e^{-ts}\,p(t;-\xi)\,dt
\nonumber\\
&=&
\int_{\!\!\sqrt[N]{\l}}^{+\infty} d\s \int_0^{+\infty} (e^{\a \xi}-1)
\Bigg(\sum_{j\in J} \t_j\,e^{-\t_j\s \xi}\Bigg)d\xi
\qquad\mbox{(by putting $\s=\sqrt[N]{s}$)}
\nonumber\\
&=&
\sum_{j\in J} \int_{\!\!\sqrt[N]{\l}}^{+\infty} d\s
\left[\t_j\int_0^{+\infty} \left(e^{-(\t_j\s-\a)\xi}-
e^{-\t_j\s \xi}\right) d\xi \right]
\nonumber\\
&=&
\sum_{j\in J} \int_{\!\!\sqrt[N]{\l}}^{+\infty}
\left(\frac{\t_j}{\t_j\s-\a}-\frac{1}{\s}\right)d\s
\;=\;\sum_{j\in J} \log\frac{\sqrt[N]{\l}}{\sqrt[N]{\l}-\a\t_j}.
\label{TLp-plus}
\eeqan
In the last step, we used the fact that the set $\{\t_j,j\in J\}$ is
invariant by conjugating.

In the same way, for $\Re(\a)\ge 0$,
\beq\label{TLp-moins}
\int_0^{+\infty} e^{-\l t}\,\frac{dt}{t}\int_{-\infty}^0
(e^{\a \xi}-1) p(t;-\xi)\,d\xi
=\sum_{k\in K}\frac{\sqrt[N]{\l}}{\sqrt[N]{\l}-\a\t_k}.
\eeq
Consequently, by choosing $\a=i\mu$ in~\refp{TLp-plus} and $\a=i\mu-\nu$
in~\refp{TLp-moins}, and using~\refp{somme-card}, it comes from~\refp{def-psi}:
$$
\exp\left[\int_0^{+\infty} e^{-\l t}\,\psi^+(\mu,\nu;t)\,\frac{dt}{t}\right]
=\frac{\l}{\displaystyle\mathop{\textstyle\prod}_{j\in J} (\!\sqrt[N]{\l}-(i\mu-\nu)\t_j)
\mathop{\textstyle\prod}_{k\in K} (\!\sqrt[N]{\l}-i\mu\t_k)}.
$$
From this and~\refp{esp-FX}, we derive the Laplace-Fourier transform
of the vector $(X(t),M(t))$. In a similar manner, we can obtain
that of $(X(t),m(t))$. The proof of Theorem~\ref{theorem-LFT-XMm} is now completed.
\fin

\brem{\label{remark-exchange-formulas}
Any of both formulas~\refp{LFT-Mm} can be deduced from the other one by using a
symmetry argument.
\bitem
\item
For even integers $N$, the obvious symmetry property $X\egallaw -X$ holds and
entails
\beqa
\Eo\!\!\left[e^{i\mu X(t)+\nu \min_{0\le s\le t}X(s)}\right]\!\!
&=&
\Eo\!\!\left[e^{-i\mu X(t)+\nu \min_{0\le s\le t}(-X(s))}\right]
\\
&=&
\Eo\!\!\left[e^{-i\mu X(t)-\nu \max_{0\le s\le t}X(s)}\right]\!.
\eeqa
Observing that in this case $\{\t_k,k\in K\}=\{-\t_j,j\in J\}$, we have
\beqa
\prod_{j\in J}\frac{\sqrt[N]{\l}}{\sqrt[N]{\l}-i\mu\t_j}
&=&
\prod_{k\in K}\frac{\sqrt[N]{\l}}{\sqrt[N]{\l}+i\mu\t_k}
\\
\noalign{\hspace{2.5em}and}
\prod_{j\in J}\frac{\sqrt[N]{\l}}{\sqrt[N]{\l}+(i\mu+\nu)\t_j}
&=&
\prod_{k\in K}\frac{\sqrt[N]{\l}}{\sqrt[N]{\l}-(i\mu+\nu)\t_k},
\eeqa
which confirms the simple relationship between both expectations~\refp{LFT-Mm}.

\item
If $N$ is odd, although this case is not recovered by~\refp{LFT-Mm},
it is interesting to note the asymmetry property $X^+\egallaw -X^-$ and
$X^-\egallaw -X^+$ where $X^+$ and $X^-$ are the pseudo-processes
respectively associated with $\k_N=+1$ and $\k_N=-1$. This would give
\beqa
\Eo\!\!\left[e^{i\mu X^+(t)+\nu \min_{0\le s\le t}X^+(s)}\right]\!\!
&=&
\Eo\!\!\left[e^{-i\mu X^-(t)+\nu \min_{0\le s\le t}(-X^-(s))}\right]
\\
&=&
\Eo\!\!\left[e^{-i\mu X^-(t)-\nu \max_{0\le s\le t}X^-(s)}\right]\!.
\eeqa
Observing that now, with similar notations,
$\{\t_j^+,j\in J^+\}=\{-\t_k^-,k\in K^-\}$ and
$\{\t_k^+,k\in K^+\}=\{-\t_j^-,j\in J^-\}$, the following relations hold:
\beqa
\prod_{j\in J^+}\frac{\sqrt[N]{\l}}{\sqrt[N]{\l}-i\mu\t_j^+}
&=&
\prod_{k\in K^-}\frac{\sqrt[N]{\l}}{\sqrt[N]{\l}+i\mu\t_k^-}
\\
\noalign{\hspace{2.5em}and}
\prod_{j\in J^-}\frac{\sqrt[N]{\l}}{\sqrt[N]{\l}+(i\mu+\nu)\t_j^-}
&=&
\prod_{k\in K^+}\frac{\sqrt[N]{\l}}{\sqrt[N]{\l}-(i\mu+\nu)\t_k^+}.
\eeqa
Hence $(X^+(t),m^+(t))$ and $(X^-(t),-M^-(t))$ should have identical
distributions, which would explain the relationship between both
expectations~\refp{LFT-Mm} in this case.
\eitem
}
\enlargethispage{\baselineskip}
%
\brem{
By choosing $\nu=0$ in~\refp{LFT-Mm}, we obtain the Fourier transform
of the \mbox{$\l$-potential} of the kernel $p$. In fact, remarking that
$$
\prod_{j\in J}(\!\sqrt[N]{\l}-i\mu\t_j)\prod_{k\in K}(\!\sqrt[N]{\l}-i\mu\t_k)
=\prod_{l=0}^{N-1} (\!\sqrt[N]{\l}-i\mu\t_l)=\l-\k_N(i\mu)^N,
$$
\refp{LFT-Mm} yields
$$
\Ex\!\!\left[\int_0^{+\infty} e^{-\l t+i\mu X(t)}\,dt\right]\!=
\frac{e^{i\mu x}}{\l-\k_N(i\mu)^N}
$$
which can be directly checked according as
$$
\int_0^{+\infty} e^{-\l t} \Ex\!\!\left[e^{i\mu X(t)}\right]dt
=\int_0^{+\infty} e^{i\mu x-\left(\l-\k_N(i\mu)^N\right)t}\,dt.
$$
}

\subsection{Density functions}\label{subsect-inverting-LFT-XMm}

We are able to invert the Laplace-Fourier transforms~\refp{LFT-Mm}
with respect to~$\mu$ and~$\nu$.

\subsubsection{Inverting with respect to $\nu$}

\bpr{
We have, for $z\ge x$,
\beqan
\int_0^{+\infty} e^{-\l t}\,dt\,\Ex\!\!\left[e^{i\mu X(t)},M(t)\in dz\right]/dz
&=&
\frac{\l^{(1-\card J)/N}e^{i\mu x}}{\dis\mathop{\textstyle{\prod}}_{k\in K}
(\!\sqrt[N]{\l}-i\mu\t_k)}
\sum_{j\in J} \t_jA_j\,e^{(i\mu-\t_j\!\!\sqrt[N]{\l}\,)(z-x)},
\nonumber\\[-3ex]
\noalign{\noindent and, for $z\le x$,}
\nonumber\\[-2\baselineskip]
\label{inverting-nu-XMm}\\
\int_0^{+\infty} e^{-\l t}\,dt\,\Ex\!\!\left[e^{i\mu X(t)},m(t)\in dz\right]/dz
&=&
-\frac{\l^{(1-\card K)/N}e^{i\mu x}}{\dis\mathop{\textstyle{\prod}}_{j\in J}
(\!\sqrt[N]{\l}-i\mu\t_j)}
\sum_{k\in K} \t_kB_k\,e^{(i\mu-\t_k\!\!\sqrt[N]{\l}\,)(z-x)}.
\nonumber
\eeqan
}
%
\dem
Observing that $\{\t_j,j\in J\}=\{\bar{\t}_j,j\in J\}=\{1/\t_j,j\in J\}$,
we have
$$
\frac{1}{\dis\mathop{\textstyle{\prod}}_{j\in J} (\!\sqrt[N]{\l}-(i\mu-\nu)\t_j)}
=\frac{1}{\dis\mathop{\textstyle{\prod}}_{j\in J}
\left(\!\sqrt[N]{\l}-\frac{i\mu-\nu}{\t_j}\right)}
=\frac{\l^{-\card J/N}}{\dis\mathop{\textstyle{\prod}}_{j\in J}
\left(1-\frac{i\mu-\nu}{\t_j\!\sqrt[N]{\l}}\right)}.
$$
Applying expansion~\refp{expansion} to $x=(i\mu-\nu)/\sqrt[N]{\l}$ yields:
\beq\label{expansion-bis}
\frac{1}{\dis\mathop{\textstyle{\prod}}_{j\in J} (\!\sqrt[N]{\l}-(i\mu-\nu)\t_j)}
=\l^{-\card J/N} \sum_{j\in J}\frac{A_j}{1-\frac{i\mu-\nu}{\t_j\!\sqrt[N]{\l}}}
=\l^{(1-\card J)/N} \sum_{j\in J}\frac{\t_jA_j}{\nu-i\mu+\t_j\!\sqrt[N]{\l}}.
\eeq
Writing now
$$
\frac{e^{-\nu x}}{\nu-i\mu+\t_j\!\sqrt[N]{\l}}
=\int_x^{+\infty} e^{-\nu z}\,e^{(i\mu-\t_j\!\!\sqrt[N]{\l}\,)(z-x)}\,dz,
$$
we find that
\beqa
\lqn{
\int_0^{+\infty} e^{-\l t}\,\Ex\!\!\left[e^{i\mu X(t)-\nu M(t)}\right]dt
}\\
&=&
\frac{\l^{(1-\card J)/N}e^{i\mu x}}{\dis\mathop{\textstyle{\prod}}_{k\in K}
(\!\sqrt[N]{\l}-i\mu\t_k)} \int_x^{+\infty} e^{-\nu z}
\Bigg(\sum_{j\in J} \t_jA_j\,e^{(i\mu-\t_j\!\!\sqrt[N]{\l}\,)(x-z)}\Bigg) dz.
\eeqa
We can therefore invert the foregoing Laplace transform with respect to $\nu$
and we get the formula~\refp{inverting-nu-XMm} corresponding the case of
the maximum functional. That corresponding in the case of the minimum functional
is obtained is a similar way.
\fin

Formulas~\refp{inverting-nu-XMm} will be used further when determining the
distributions of $(\tau_a^+,X(\tau_a^+))$ and $(\tau_a^-,X(\tau_a^-))$.

\subsubsection{Inverting with respect to $\mu$}

\bth{
The Laplace transforms with respect to time $t$ of the joint density of
$X(t)$ and, respectively, $M(t)$ and $m(t)$, are given, for $z\ge x\vee y$,
by
\beqan
\int_0^{+\infty} e^{-\l t}\,dt\,\P_x\{X(t)\in dy,M(t)\in dz\}/dy\,dz
&=&
\frac{1}{\l}\,\f_{\l}(x-z)\,\psi_{\l}(z-y),
\nonumber\\[-1ex]
\noalign{\noindent and, for $z\le x\wedge y$,}
\nonumber\\[-2\baselineskip]
\label{densities-XMm}\\
\int_0^{+\infty} e^{-\l t}\,dt\,\P_x\{X(t)\in dy,m(t)\in dz\}/dy\,dz
&=&
\frac{1}{\l}\,\psi_{\l}(x-z)\,\f_{\l}(z-y),
\nonumber
\eeqan
where the functions $\f_{\l}$ and $\psi_{\l}$ are defined by~\refp{phi-and-psi}.
}
%
\dem
Let us write the following equality, as in the previous subsubsection
(see~\refp{expansion-bis}):
$$
\frac{1}{\dis\mathop{\textstyle{\prod}}_{K\in K} (\!\sqrt[N]{\l}-i\mu\t_k)}
=-\l^{(1-\card K)/N} \sum_{k\in K}\frac{\t_kB_k}{i\mu-\t_k\!\sqrt[N]{\l}}.
$$
Set
$$
G(\l,\mu;x,z)=\int_0^{+\infty} e^{-\l t}\,dt\,
\Ex\!\!\left[e^{i\mu X(t)},M(t)\in dz\right]/dz
$$
We get, by~\refp{inverting-nu-XMm} and~\refp{somme-card}, for $z\ge x$,
\beqa
G(\l,\mu;x,z)
&=&
\frac{\l^{(1-\card J)/N}e^{i\mu x}}{\dis\mathop{\textstyle{\prod}}_{k\in K}
(\!\sqrt[N]{\l}-i\mu\t_k)}
\sum_{j\in J} \t_jA_j\,e^{(i\mu-\t_j\!\!\sqrt[N]{\l}\,)(z-x)}
\\
&=&
-\l^{(2-\card J-\card K)/N} e^{i\mu x}
\sum_{k\in K}\frac{\t_kB_k}{i\mu-\t_k\!\sqrt[N]{\l}}
\sum_{j\in J} \t_jA_j\,e^{(i\mu-\t_j\!\!\sqrt[N]{\l}\,)(z-x)}
\\
&=&
-\l^{2/N-1} \sum_{j\in J,k\in K} e^{\t_j\!\!\sqrt[N]{\l}\,x} \t_jA_j\,\t_kB_k\,
\frac{e^{(i\mu-\t_j\!\!\sqrt[N]{\l}\,)\,z}}{i\mu-\t_k\!\sqrt[N]{\l}}.
\eeqa
Writing now
$$
\frac{e^{(i\mu-\t_j\!\!\sqrt[N]{\l}\,)\,z}}{i\mu-\t_k\!\sqrt[N]{\l}}
=e^{(\t_k-\t_j)\!\sqrt[N]{\l}\,z}
\int_{-\infty}^z e^{(i\mu-\t_k\!\!\sqrt[N]{\l}\,)\,y}\,dy
$$
gives
$$
G(\l,\mu;x,z)=-\l^{2/N-1} \ind_{\{z\ge x\}}
\int_{-\infty}^z e^{i\mu y}\Bigg[\sum_{j\in J,k\in K} \t_jA_j\,\t_kB_k \,
e^{\!\!\sqrt[N]{\l}\,(\t_j x-\t_k y+(\t_k-\t_j)z)}\Bigg] dy
$$
and then
\beqa
\lqn{
\int_0^{+\infty} e^{-\l t}\,dt\,\P_x\{X(t)\in dy,M(t)\in dz\}/dy\,dz
}\\
&=&
-\l^{2/N-1} \sum_{j\in J,k\in K} \t_jA_j\,\t_kB_k \,
e^{\!\!\sqrt[N]{\l}\,(\t_j x-\t_k y+(\t_k-\t_j)z)}\,\ind_{\{z\ge x\vee y\}}.
\eeqa
This proves~\refp{densities-XMm} in the case of the maximum functional
and the formula corresponding to the minimum functional can be proved in
a same manner.
\fin

%
\brem{
Formulas~\refp{densities-XMm} contain in particular the Laplace transforms
of $X(t)$, $M(t)$ and $m(t)$ separately. As a verification, we
integrate~\refp{densities-XMm} with respect to $y$ and $z$ separately.
\bitem
\item
By integrating with respect to $y$ on $[z,+\infty)$ for $z\le x$, we get
\beqa
\lqn{\hspace{2.5em}
\int_0^{+\infty} e^{-\l t}\,dt\,\P_x\{m(t)\in dz\}/dz
}\\
&=&
-\l^{2/N-1}\sum_{j\in J} \t_jA_j \int_z^{+\infty}
e^{-\t_j\!\!\sqrt[N]{\l}(y-z)}\,dy
\sum_{k\in K} \t_kB_k \,e^{\t_k\!\!\sqrt[N]{\l}(x-z)}
\\
&=&
-\l^{1/N-1}\sum_{j\in J} A_j\sum_{k\in K} \t_kB_k
\,e^{\t_k\!\!\sqrt[N]{\l}(x-z)}
\\
&=&
-\l^{1/N-1} \sum_{k\in K} \t_kB_k \,e^{\t_k\!\!\sqrt[N]{\l}(x-z)}
\;=\;\frac{1}{\l}\,\psi_{\l}(x-z).
\eeqa
We used the relation $\sum_{j\in J} A_j=1$; see
Subsection~\ref{subsect-elem-prop}. We retrieve the Laplace
transform~\refp{densities-Mm} of the distribution of $m(t)$.

\item
Suppose for instance that $x\le y$. Let us integrate~\refp{densities-XMm}
now with respect to $z$ on $(-\infty,x]$. This gives
\beqa
\lqn{\hspace{2.5em}
\int_0^{+\infty} e^{-\l t}\,dt\,\P_x\{X(t)\in dy\}/dy
}\\
&=&
-\l^{2/N-1} \sum_{j\in J,k\in K} \t_jA_j\,\t_kB_k
\,e^{\t_k\!\!\sqrt[N]{\l}\,x-\t_j\!\!\sqrt[N]{\l}\,y}
\int_{-\infty}^{x} e^{(\t_j-\t_k)\!\sqrt[N]{\l}\,z}\,dz
\\
&=&
\l^{1/N-1} \sum_{j\in J,k\in K} \frac{\t_jA_j\,\t_kB_k}{\t_k-\t_j}
\,e^{\t_j\!\!\sqrt[N]{\l}(x-y)}
\\
&=&
\l^{1/N-1} \sum_{j\in J} \left(\sum_{k\in K} \frac{\t_kB_k}{\t_k-\t_j}\right)
\t_jA_j\,e^{\t_j\!\!\sqrt[N]{\l}(x-y)}
\\
&=&
\frac1N\,\l^{1/N-1} \sum_{j\in J} \t_j\,e^{\t_j\!\!\sqrt[N]{\l}(x-y)},
\eeqa
where we used~\refp{somme-partial-frac} in the last step. We retrieve the
$\l$-potential~\refp{potential} of the pseudo-process $(X(t))_{t\ge0}$ since
$$
\int_0^{+\infty} e^{-\l t}\,dt\,\P_x\{X(t)\in dy\}/dy
=\int_0^{+\infty} e^{-\l t}\,p(t;x-y)\,dt.
$$
\eitem}

\brem{
Consider the reflected process at its maximum $(M(t)-X(t))_{t\ge 0}$.
The joint distribution of $(M(t),M(t)-X(t))$ writes in terms of
the joint distribution of $(X(t),M(t))$, for $x=0$ (set $\P=\P_0$ for short)
and $z,\z\ge 0$, as:
$$
\P\{M(t)\in dz,M(t)-X(t)\in d\z\}
=\P\{X(t)\in z-d\z,M(t)\in dz\}.
$$
Formula~\refp{densities-XMm} writes
\beqan
\lqn{
\int_0^{+\infty} \l\,e^{-\l t}\,dt\,\P\{M(t)\in dz,M(t)-X(t)\in d\z\}/dz\,d\z
\;=\;\f_{\l}(z)\psi_{\l}(-\z)
}\nonumber\\
&=&
\int_0^{+\infty} \l\,e^{-\l t}\,dt\,\P\{M(t)\in dz\}/dz\times
\int_0^{+\infty} \l\,e^{-\l t}\,dt\,\P\{-m(t)\in d\z\}/d\z.
\label{factorization}
\eeqan
By introducing an exponentially distributed time $T_{\l}$ with parameter $\l$
which is independent of $(X(t))_{t\ge 0}$, \refp{factorization} reads
$$
\P\{M(T_{\l})\in dz,M(T_{\l})-X(T_{\l})\in d\z\}
=\P\{M(T_{\l})\in dz\}\,\P\{-m(T_{\l})\in d\z\}.
$$
This relationship may be interpreted by saying that
$-m(T_{\l})$ and $M(T_{\l})-X(T_{\l})$ admit the same distribution and
$M(T_{\l})$ and $M(T_{\l})-X(T_{\l})$ are independent.
}

\brem{
The similarity between both formulas~\refp{densities-XMm} may be explained by
invoking a ``duality'' argument.
In effect, the dual pseudo-process
$(X^*(t))_{t\ge 0}$ of $(X(t))_{t\ge 0}$ is defined by $X^*(t)=-X(t)$ for all
$t\ge 0$ and we have the following equality related to the inversion of
the extremities (see~\cite{arcsine}):
\beqa
\P_x\{X(t)\in dy,M(t)\in dz\}/dy\,dz
&=&
\P_y\{X^*(t)\in dx,m^*(t)\in dz\}/dx\,dz
\\
&=&
\P_{-y}\{X(t)\in d(-x),m(t)\in d(-z)\}/dx\,dz.
\eeqa
}
\brem{
Let us expand the function $\f_{\l}$ as $\l\to 0^+$:
\beqa
\f_{\l}(\xi)
&=&
\sqrt[N]{\l}\, \sum_{j\in J} \t_j A_j\left[ \sum_{l=0}^{\card J-1}
\frac{[\t_j\sqrt[N]{\l}\,\xi]^l}{l!} +o\Big(\l^{(\card J-1)/N}\Big)\right]
\\
&=&
\sum_{l=0}^{\card J-1} \Bigg(\sum_{j\in J} \t_j^{l+1} A_j\Bigg)
\frac{\l^{(l+1)/N}\xi^l}{l!} +o\Big(\l^{\card J/N}\Big).
\eeqa
We have by~\refp{expansion} (for $x=0$) $\sum_{j\in J} \t_j^{l+1} A_j=0$
for $0\le l\le \card J-2$ and $\sum_{j\in J} \t_j^{\card J} A_j=
(-1)^{\card J-1}\prod_{j\in J}\t_j.$ Hence
\beq\label{equivalent-phi}
\f_{\l}(\xi)\underset{\l \to 0^+}{\overset{}{\sim}}
(-1)^{\card J-1}\prod_{j\in J}\t_j\,
\frac{\xi^{\card J-1}}{(\card J-1)!}\,\l^{\card J/N}.
\eeq
Similarly
\beq\label{equivalent-psi}
\psi_{\l}(\xi)\underset{\l \to 0^+}{\overset{}{\sim}}
(-1)^{\card K}\prod_{k\in K}\t_k\,
\frac{\xi^{\card K-1}}{(\card K-1)!}\,\l^{\card K/N}.
\eeq
As a result, putting~\refp{equivalent-phi} and~\refp{equivalent-psi}
into~\refp{densities-XMm} and using~\refp{somme-card}
and $\prod_{l=0}^{N-1} \t_l=(-1)^{N-1}\k_N$ lead to
$$
\int_0^{+\infty} e^{-\l t}\,dt\,\P_x\{X(t)\in dy,M(t)\in dz\}/dy\,dz
\underset{\l \to 0^+}{\overset{}{\sim}}
\k_N\,\frac{(x-z)^{\card J-1}(z-y)^{\card K-1}}{(\card J-1)!\,(\card K-1)!}.
$$
By integrating this asymptotic with respect to $z$, we derive the value
of the so-called $0$-potential of the absorbed pseudo-process
(see~\cite{nish3} for the definition of several kinds of absorbed or killed
pseudo-processes):
$$
\int_0^{+\infty} \P_x\{X(t)\in dy,M(t)\le a\}/dy
=(-1)^{\card J-1}\k_N \int_{x\vee y}^a
\frac{(z-x)^{\card J-1}(z-y)^{\card K-1}}{(\card J-1)!\,(\card K-1)!}\,dz.
$$
}

\subsubsection{Inverting with respect to $\l$}

Formulas~\refp{densities-XMm} may be inverted with respect to $\l$
and an expression by means of the successive derivatives of the kernel $p$
may be obtained for the densities of $(X(t),M(t))$ and $(X(t),m(t))$.
However, the computations and the results are cumbersome and
we prefer to perform them in the case of the distribution functions.
They are exhibited in Subsection~\ref{subsect-dist-func-XMm}.

\subsection{Density functions: particular cases}\label{subsect-particular-cases}

In this subsection, we pay attention to the cases $N\in\{2,3,4\}$.
Although our results are not justified when $N$ is odd, we nevertheless retrieve
well-known results in the literature related to the case $N=3$. In order to
lighten the notations, we set for, $\Re(\l)>0$,
\beqa
\Phi_{\l}(x,y,z)
&=&
\int_0^{+\infty} e^{-\l t}\,dt\,\P_x\{X(t)\in dy,M(t)\in dz\}/dy\,dz,
\\[1ex]
\Psi_{\l}(x,y,z)
&=&
\int_0^{+\infty} e^{-\l t}\,dt\,\P_x\{X(t)\in dy,m(t)\in dz\}/dy\,dz.
\eeqa
\bex{\textsl{Case $N=2$:} using the numerical results of Example~\ref{case-2}
gives
$$
\f_{\l}(\xi) =\sqrt{\l}\,e^{\sqrt{\l}\,\xi}
\quad\mbox{and}\quad
\psi_{\l}(\xi)=\sqrt{\l}\,e^{-\sqrt{\l}\,\xi},
$$
and then
$$
\Phi_{\l}(x,y,z)=e^{\sqrt{\l}(x+y-2z)}\,\ind_{\{z\ge x\vee y\}}
\quad\mbox{and}\quad
\Psi_{\l}(x,y,z)=e^{\sqrt{\l}(2z-x-y)}\,\ind_{\{z\le x\wedge y\}}.
$$
}

\bex{\textsl{Case $N=3$:} referring to Example~\ref{case-3}, we have
\bitem
\item
for $\k_3=+1$:
\beqa
\f_{\l}(\xi)
&=&
\sqrt[3]{\l}\,e^{\sqrt[3]{\l}\,\xi},
\\
\psi_{\l}(\xi)
&=&
-\frac{i\,\sqrt[3]{\l}}{\sqrt3} \left(e^{e^{i\,2\pi/3}\sqrt[3]{\l}\,\xi}
-\,e^{e^{-i\,2\pi/3}\sqrt[3]{\l}\,\xi}\right)
\;=\;\frac{2\,\sqrt[3]{\l}}{\sqrt3}\,e^{-\frac{\sqrt[3]{\l}}{2}\,\xi}
\,\sin\bigg(\frac{\sqrt3}{2}\,\sqrt[3]{\l}\,\xi\bigg),
\eeqa
which gives
\beqa
\Phi_{\l}(x,y,z)
&=&
\frac{2}{\sqrt3\sqrt[3]{\l}}\,e^{\sqrt[3]{\l}(x+\frac12\,y-\frac32\,z)}
\sin\bigg(\frac{\sqrt3}{2}\,\sqrt[3]{\l}(z-y)\bigg)\,\ind_{\{z\ge x\vee y\}},
\\
\Psi_{\l}(x,y,z)
&=&
\frac{2}{\sqrt3\sqrt[3]{\l}}\,e^{\sqrt[3]{\l}(\frac32\,z-\frac12\,x-y)}
\sin\bigg(\frac{\sqrt3}{2}\,\sqrt[3]{\l}(x-z)\bigg)\,\ind_{\{z\le x\wedge y\}};
\eeqa

\item
for $\k_3=-1$,
\beqa
\f_{\l}(\xi)
&=&
\frac{i\,\sqrt[3]{\l}}{\sqrt3} \left(e^{e^{i\,\pi/3}\sqrt[3]{\l}\,\xi}
-e^{e^{-i\,\pi/3}\sqrt[3]{\l}\,\xi}\right)
\;=\;-\frac{2\,\sqrt[3]{\l}}{\sqrt3}\,e^{\frac{\sqrt[3]{\l}}{2}\,\xi}
\,\sin\bigg(\frac{\sqrt3}{2}\,\sqrt[3]{\l}\,\xi\bigg),
\\
\psi_{\l}(\xi)
&=&
\sqrt[3]{\l}\,e^{-\sqrt[3]{\l}\,\xi},
\eeqa
which gives
\beqa
\Phi_{\l}(x,y,z)
&=&
\frac{2}{\sqrt3\sqrt[3]{\l}}\,e^{\sqrt[3]{\l}(\frac12\,x+y-\frac32\,z)}
\sin\bigg(\frac{\sqrt3}{2}\,\sqrt[3]{\l}(z-x)\bigg)\,\ind_{\{z\ge x\vee y\}},
\\
\Psi_{\l}(x,y,z)
&=&
\frac{2}{\sqrt3\sqrt[3]{\l}}\,e^{\sqrt[3]{\l}(\frac32\,z-x-\frac12\,y)}
\sin\bigg(\frac{\sqrt3}{2}\,\sqrt[3]{\l}(y-z)\bigg)\,\ind_{\{z\le x\wedge y\}}.
\eeqa
\eitem
}
\bex{\textsl{Case $N=4$:} the numerical results of Example~\ref{case-4} yield
\beqa
\f_{\l}(\xi)
&=&
-\frac{i\,\sqrt[4]{\l}}{\sqrt2}\left(e^{e^{-i\,\pi/4}\sqrt[4]{\l}\,\xi}
-e^{e^{i\,\pi/4}\sqrt[4]{\l}\,\xi}\right)
\;=\;-\sqrt2\,\sqrt[4]{\l}\,e^{\frac{\sqrt[4]{\l}}{\sqrt2}\,\xi}
\,\sin\bigg(\frac{\sqrt[4]{\l}}{2}\,\xi\bigg),
\\
\psi_{\l}(\xi)
&=&
\frac{i\,\sqrt[4]{\l}}{\sqrt2} \left(e^{e^{-i\,3\pi/4}\sqrt[4]{\l}\,\xi}
-e^{e^{i\,3\pi/4}\sqrt[4]{\l}\,\xi}\right)
\;=\;\sqrt2\,\sqrt[4]{\l}\,e^{-\frac{\sqrt[4]{\l}}{\sqrt2}\,\xi}
\,\sin\bigg(\frac{\sqrt[4]{\l}}{2}\,\xi\bigg),
\eeqa
which gives
\beqa
\Phi_{\l}(x,y,z)
&=&
\frac{1}{\sqrt{\l}}\,e^{\frac{\sqrt[4]{\l}}{\sqrt2}(x+y-2z)}
\bigg[\cos\bigg(\frac{\sqrt[4]{\l}}{\sqrt2}\,(x-y)\bigg)
-\cos\bigg(\frac{\sqrt[4]{\l}}{\sqrt2}\,(x+y-2z)\bigg)\bigg]
\,\ind_{\{z\ge x\vee y\}},
\\
\Psi_{\l}(x,y,z)
&=&
\frac{1}{\sqrt{\l}}\,e^{\frac{\sqrt[4]{\l}}{\sqrt2}(2z-x-y)}
\bigg[\cos\bigg(\frac{\sqrt[4]{\l}}{\sqrt2}\,(x-y)\bigg)
-\cos\bigg(\frac{\sqrt[4]{\l}}{\sqrt2}\,(x+y-2z)\bigg)\bigg]
\,\ind_{\{z\le x\wedge y\}}.
\eeqa
}

\subsection{Distribution functions}\label{subsect-dist-func-XMm}

In this part, we integrate~\refp{densities-XMm} in view to get the
distribution function of the vector $(X(t),M(t))$:
$\P_x\{X(t)\le y,M(t)\le z\}$.
Obviously, if $x>z$, this quantity vanishes. Suppose now $x\le z$. We must
consider the cases $y\le z$ and $z\le y$. In the latter, we have
$\P_x\{X(t)\le y,M(t)\le z\}=\P\{M(t)\le z\}$ and this quantity is given
by~\refp{dist-func-Mm}. So, we assume that $z\ge x\vee y$. Actually,
the quantity $\P_x\{X(t)\le y,M(t)\ge z\}$ is easier to derive.

\subsubsection{Laplace transform}

Put for $\Re(\l)>0$:
\beqa
F_{\l}(x,y,z)
&=&
\int_0^{+\infty} e^{-\l t}\,\P_x\{X(t)\le y,M(t)\le z\}\,dt
\\
\tilde{F}_{\l}(x,y,z)
&=&
\int_0^{+\infty} e^{-\l t}\,\P_x\{X(t)\le y,M(t)\ge z\}\,dt.
\eeqa
The functions $F_{\l}$ and $\tilde{F}_{\l}$ are related together through
\beq\label{relation-F-Ftilde}
F_{\l}(x,y,z)+\tilde{F}_{\l}(x,y,z)=\Psi(\l;x-y)
\eeq
where $\Psi$ is given by \refp{potential-bis}.
Using~\refp{densities-XMm}, we get
\beqa
\lqn{\tilde{F}_{\l}(x,y,z)}
\\[-2ex]
&=&
\int_{-\infty}^y \int_z^{+\infty} \!\!\!\int_0^{+\infty}
e^{-\l t}\,\P_x\{X(t)\in d\xi,M(t)\in d\z\}\,dt
\\
&=&
-\l^{2/N-1} \sum_{j\in J,k\in K} \t_jA_j\,\t_kB_k \,e^{\t_j\!\!\sqrt[N]{\l}\,x}
\int_{-\infty}^y e^{-\t_k\!\!\sqrt[N]{\l}\,\xi}\,d\xi
\int_z^{+\infty} e^{(\t_k-\t_j)\!\sqrt[N]{\l}\,\z}
\ind_{\{\z\ge x\vee \xi\}}\,d\z.
\eeqa
We plainly have $\z\ge z\ge x\vee y\ge x\vee \xi$ over the integration set
$(-\infty,y]\times[z,+\infty)$. So, the indicator $\ind_{\{\z\ge x\vee \xi\}}$
is useless and we obtain the following expression for $\tilde{F}_{\l}$.
%
\bpr{\label{prop-dist-func-inter}
We have for $z\ge x\vee y$ and $\Re(\l)>0$:
\beqa
\int_0^{+\infty} e^{-\l t}\,\P_x\{X(t)\le y\le z\le M(t)\}\,dt
&=&
\frac{1}{\l} \sum_{j\in J,k\in K} \frac{\t_jA_jB_k}{\t_j-\t_k}
\,e^{\t_j\!\!\sqrt[N]{\l}(x-z)+\t_k\!\!\sqrt[N]{\l}(z-y)}
\\[-2ex]
\noalign{\noindent and for $z\le x\wedge y$:}
\\[-1ex]
\int_0^{+\infty} e^{-\l t}\,\P_x\{X(t)\ge y\ge z\ge m(t)\}\,dt
&=&
\frac{1}{\l} \sum_{j\in J,k\in K} \frac{A_j\t_kB_k}{\t_k-\t_j}
\,e^{\t_j\!\!\sqrt[N]{\l}(z-y)+\t_k\!\!\sqrt[N]{\l}(x-z)}.
\eeqa
}
%
As a result, combining the above formulas with~\refp{relation-F-Ftilde}, the
distribution function of the couple $(X(t),M(t))$ emerges and that of
$(X(t),m(t))$ is obtained in a similar way.
%
\bth{
The distribution functions of $(X(t),M(t))$ and $(X(t),m(t))$ are respectively
determined through their Laplace transforms with respect to $t$ by
\beqan
\lqn{
\int_0^{+\infty} e^{-\l t}\,\P_x\{X(t)\le y,M(t)\le z\}\,dt
}\nonumber\\
&=&
\left\{\begin{array}{ll}
\dis\frac{1}{\l} \sum_{j\in J,k\in K} \frac{\t_jA_jB_k}{\t_k-\t_j}
\,e^{\t_j\!\!\sqrt[N]{\l}(x-z)+\t_k\!\!\sqrt[N]{\l}(z-y)}
+\frac{1}{N\l} \sum_{k\in K} e^{\t_k\!\!\sqrt[N]{\l}(x-y)}
&\mbox{\hspace{-.5em}if $y\le x\le z$, \hspace{-1em}}
\\[4ex]
\dis\frac{1}{\l}\Bigg[1 -\frac 1N\sum_{j\in J} e^{\t_j\!\!\sqrt[N]{\l}(x-y)}
+\sum_{j\in J,k\in K}\frac{\t_jA_jB_k}{\t_k-\t_j}
\,e^{\t_j\!\!\sqrt[N]{\l}(x-z)+\t_k\!\!\sqrt[N]{\l}(z-y)}\Bigg]
&\mbox{\hspace{-.5em}if $x\le y\le z$, \hspace{-1em}}
\end{array}\right.
\nonumber\\[-10.5ex]
&&\label{distXM}\\[7ex]
\noalign{\noindent and}
\lqn{
\int_0^{+\infty} e^{-\l t}\,\P_x\{X(t)\ge y,m(t)\ge z\}\,dt
}\nonumber\\
&=&
\left\{\begin{array}{ll}
\dis\frac{1}{\l} \sum_{j\in J,k\in K} \frac{A_j\t_kB_k}{\t_j-\t_k}
\,e^{\t_j\!\!\sqrt[N]{\l}(x-z)+\t_k\!\!\sqrt[N]{\l}(z-y)}
+\frac{1}{N\l} \sum_{j\in J} e^{\t_j\!\!\sqrt[N]{\l}(x-y)}
&\mbox{\hspace{-.5em}if $z\le x\le y$, \hspace{-1em}}
\\[4ex]
\dis\frac{1}{\l}\Bigg[1 -\frac 1N\sum_{k\in K} e^{\t_k\!\!\sqrt[N]{\l}(x-y)}
+\sum_{j\in J,k\in K}\frac{A_j\t_kB_k}{\t_j-\t_k}
\,e^{\t_j\!\!\sqrt[N]{\l}(x-z)+\t_k\!\!\sqrt[N]{\l}(z-y)}\Bigg]
&\mbox{\hspace{-.5em}if $z\le y\le x$. \hspace{-1em}}
\end{array}\right.
\nonumber
\eeqan
}

\subsubsection{Inverting the Laplace transform}

%
\bth{\label{th-representations}
The distribution function of $(X(t),M(t))$ admits the following representation:
\beq\label{dist-func-XM}
\P_x\{X(t)\le y\le z\le M(t)\}=\!\!\!\!\!\!\sum_{k\in K\atop 0\le m\le\card J-1}
\!\!\!\!\!\!a_{km}\int_0^t\!\!\int_0^s\frac{\partial^m p}{\partial x^m}(\s;x-z)
\,\frac{I_{k0}(s-\s;z-y)}{(t-s)^{1-(m+1)/N}}\,ds\,d\s
\eeq
where $I_{k0}$ is given by~\refp{Ik0} and
$$
a_{km}=\frac{NB_k}{\G(\frac{m+1}{N})}\sum_{j\in J} \frac{A_j\a_{jm}}{\t_j-\t_k},
$$
the $\a_{jm}$'s being some coefficients given by~\refp{coefficients}.
}
%
\dem
We intend to invert the Laplace transform~\refp{distXM}.
For this, we interpret both exponentials
$e^{\t_j\!\!\sqrt[N]{\l}(x-z)}$ and $e^{\t_k\!\!\sqrt[N]{\l}(z-y)}$
as Laplace transforms in two different manners:
one is the Laplace transform of a combination of the successive
derivatives of the kernel $p$, the other one is the Laplace
transform of a function which is closely related to the density
of some stable distribution. More explicitly, we proceed as follows.
\bitem
\item
On one hand, we start from the $\l$-potential~\refp{potential} that we shall
call~$\Phi$:
$$
\Phi(\l;\xi)=\frac{1}{N\l^{1-1/N}}\sum_{j\in J}
\t_j e^{\t_j\!\!\sqrt[N]{\l}\,\xi} \quad\mbox{for $\xi\le 0$.}
$$
Differentiating this potential $(\card J-1)$ times with respect to $\xi$
leads to the Vandermonde system
of $\card J$ equations where the exponentials $e^{\t_j\!\!\sqrt[N]{\l}\,\xi}$
are unknown:
$$
\sum_{j\in J} \t_j^{l+1} e^{\t_j\!\!\sqrt[N]{\l}\,\xi}=
N\l^{1-(l+1)/N} \,\frac{\partial^l \Phi}{\partial x^l}(\l;\xi)
\quad\mbox{for $0\le l\le \card J-1$.}
$$
Introducing the solutions $\a_{jm}$ of the $\card J$ elementary Vandermonde systems
(indexed by $m$ varying from $0$ to $\card J-1$):
$$
\sum_{j\in J} \t_j^l\a_{jm}=\d_{lm},\quad 0\le l\le\card J-1,
$$
we extract
\beqa
\frac{\t_j}{\l}\,e^{\t_j\!\!\sqrt[N]{\l}\,\xi}
&=&
N\sum_{m=0}^{\card J-1} \frac{\a_{jm}}{\l^{(m+1)/N}}\,
\frac{\partial^m \Phi}{\partial x^m}(\l;\xi) 
\\
&=&
\int_0^{+\infty} e^{-\l t}\,dt
\int_0^t \sum_{m=0}^{\card J-1} \frac{N\a_{jm}}{\G(\frac{m+1}{N})}
\frac{\partial^m p}{\partial x^m}(s;\xi)\,\frac{ds}{(t-s)^{1-(m+1)/N}}.
\eeqa
The explicit expression of $\a_{jm}$ is
\beq\label{coefficients}
\a_{jm}=(-1)^m \frac{\s_{\card J-1-m}(\t_l,l\in J\setminus\{j\})}
{\prod_{l\in J\setminus\{j\}}(\t_l-\t_j)}
=\frac{(-1)^m}{\prod_{l\in J}\t_l}\,c_{j,\card J-1-m}\t_jA_j
\eeq
where the coefficients $c_{jq}$, $0\le q\le\card J-1$, are the elementary
symmetric functions of the $\t_l$'s, $l\in J\setminus\{j\}$, that is
$c_{j0}=1$ and for $1\le q\le \card J-1$,
$$
c_{jq}=\s_q\left(\t_l,l\in J\setminus\{j\}\right)
=\sum_{l_1,\ldots,l_q\in J\setminus\{j\}\atop l_1<\cdots<l_q}\t_{l_1}\cdots\t_{l_q}.
$$

\item
On the other hand, using the Bromwich formula, the function
$\xi \longmapsto e^{\t_k\!\!\sqrt[N]{\l}\,\xi}$ can be written as a Laplace transform.
Indeed, referring to Section~\ref{subsubsect-invert-lambda-tau-Xtau}, we have for
$k\in K$ and $\xi\ge 0$,
$$
e^{\t_k\!\!\sqrt[N]{\l}\,\xi}=\int_0^{+\infty} e^{-\l t} I_{k0}(t;\xi)\,dt
$$
where $I_{k0}$ is given by~\refp{Ik0}.
\eitem
Consequently, the sum lying in Proposition~\ref{prop-dist-func-inter} may be
written as a Laplace transform which gives the
representation~\refp{dist-func-XM} for the
the distribution function of $(X(t),M(t))$.
\fin

%
\brem{
A similar expression obtained by exchanging the roles of the indices $j$ and $k$
in the above discussion and slightly changing the coefficient $a_{km}$
into another $b_{jn}$ may be derived:
\beq\label{dist-func-XM-bis}
\P_x\{X(t)\le y\le z\le M(t)\}=\!\!\!\!\!\!
\sum_{j\in J\atop 0\le n\le\card K-1} \!\!\!\!\!\!b_{jn}
\int_0^t\!\!\int_0^s\frac{\partial^n p}{\partial x^n}(\s;z-y)
\,\frac{I_{j0}(s-\s;x-z)}{(t-s)^{1-(n+1)/N}}\,ds\,d\s
\eeq
where
$$
b_{jn}=\frac{N\t_jA_j}{\G(\frac{n+1}{N})}\sum_{k\in K}
\frac{\bar{\t}_kB_k\b_{kn}}{\t_k-\t_j}.
$$
However, the foregoing result involves the same number of integrals as
that displayed in Theorem~\ref{th-representations}.
}


\subsection{Distribution functions: particular cases}\label{subsect-dist-func-XMm-part}

Here, we write out~\refp{distXM} and~\refp{dist-func-XM} or~\refp{dist-func-XM-bis}
in the cases $N\in\{2,3,4\}$ with the same remark about the case $N=3$
already mentioned at the beginning of Subsection~\ref{subsect-particular-cases}.
The expressions are rather simple and remarkable.

\bex{\textsl{Case $N=2$:}
the double sum lying in~\refp{distXM} reads
$$
\sum_{j\in J,k\in K} \frac{\t_jA_jB_k}{\t_k-\t_j}
\,e^{\t_j\sqrt[4]{\l}(x-z)+\t_k\sqrt[4]{\l}(z-y)}
=\frac{\t_1A_1B_0}{\t_0-\t_1}\,e^{\t_1\sqrt{\l}(x-z)+\t_0\sqrt{\l}(z-y)}
$$
with $\frac{\t_1A_1B_0}{\t_0-\t_1}=-\frac12$, and then
\beqa
F_{\l}(x,y,z)
&=&
\left\{\begin{array}{ll}
\dis\frac{1}{2\l}\,\Big[e^{-\sqrt{\l}(x-y)}-e^{\sqrt{\l}(x+y-2z)}\Big]
&\mbox{if $y\le x\le z$,}
\\[2ex]
\dis\frac{1}{\l}-\frac{1}{2\l}\,\Big[e^{\sqrt{\l}(x-y)}+e^{\sqrt{\l}(x+y-2z)}\Big]
&\mbox{if $x\le y\le z$.}
\end{array}\right.
\eeqa
Formula~\refp{dist-func-XM} writes
$$
\P_x\{X(t)\le y\le z\le M(t)\}=
a_{00}\int_0^t\!\!\int_0^s p(\s;x-z)\,I_{00}(s-\s;z-y)\,\frac{ds\,d\s}{\sqrt{t-s}}
$$
with
$$
p(t;\xi)=\frac{1}{\sqrt{\pi t}}\,e^{-\frac{\xi^2}{4t}}.
$$
The reciprocal relations, which are valid for $\xi\le 0$,
$$
\Phi(\l;\xi)=\frac{e^{\sqrt{\l}\,\xi}}{2\sqrt{\l}}
\quad\mbox{and}\quad
e^{\sqrt{\l}\,\xi}=2\sqrt{\l}\,\a_{10}\Phi(\l;\xi)
$$
imply that $\a_{10}=1$. Then
$
a_{00}=\frac{2B_0}{\G(1/2)}\,\frac{A_1\a_{10}}{\t_1-\t_0}=\frac{1}{\sqrt{\pi}}.
$
On the other hand, we have for $\xi\ge 0$, by~\refp{Ik0},
\beqa
I_{00}(t;\xi)
&=&
\frac{i\t_0\xi}{2\pi t}
\left[i\int_0^{+\infty} e^{-t\l^2+i\t_0\xi\l} \,d\l
+i\int_0^{+\infty} e^{-t\l^2-i\t_0\xi\l} \,d\l\right]
\\
&=&
\frac{\xi}{2\pi t}
\left[\int_0^{+\infty} e^{-t\l^2+i\xi\l} \,d\l
+\int_0^{+\infty} e^{-t\l^2-i\xi\l} \,d\l\right]
\\
&=&
\frac{\xi}{2\pi t} \int_{-\infty}^{+\infty} e^{-t\l^2+i\xi\l} \,d\l
=\frac{\xi}{2\sqrt{\pi}\,t^{3/2}}\,e^{-\frac{\xi^2}{4t}}.
\eeqa
Consequently,
$$
\P_x\{X(t)\le y\le z\le M(t)\}
=\frac{z-y}{4\pi^{3/2}} \int_0^t\!\!\int_0^s
\frac{e^{-\frac{(x-z)^2}{4\s}-\frac{(z-y)^2}{4(s-\s)}}}{\sqrt{\s}(s-\s)^{3/2}\sqrt{t-s}}\,ds\,d\s.
$$
Using the substitution $\s=\frac{s^2}{u+s}$ together with a well-known integral related
to the modified Bessel function $K_{1/2}$, we get
\beqa
\int_0^s \frac{e^{-\frac{(x-z)^2}{4\s}-\frac{(z-y)^2}{4(s-\s)}}}{\sqrt{\s}(s-\s)^{3/2}}\,d\s
&=&
\frac{e^{-\frac{(x-z)^2+(z-y)^2}{4s}}}{\sqrt s}
\int_0^{\infty} e^{-\frac{(x-z)^2}{s^2}\,u-\frac{(z-y)^2}{4u}}\, \frac{du}{u^{3/2}}
\\
&=&
\frac{\sqrt{\pi}}{\sqrt{(z-y)s}}\,e^{-\frac{(2z-x-y)^2}{4s}}.
\eeqa
Then
$$
\P_x\{X(t)\le y\le z\le M(t)\}=
\frac{1}{2\pi}\int_0^t \frac{e^{-\frac{(2z-x-y)^2}{4s}}}{\sqrt{\s(t-s)}}\,ds.
$$
Finally, it can be easily checked, by using the Laplace transform, that
$$
\int_0^t \frac{e^{-\frac{(2z-x-y)^2}{4s}}}{\sqrt{\s(t-s)}}\,ds
=\sqrt{\pi} \int_{2z-x-y}^{\infty} \frac{e^{-\frac{\xi^2}{4t}}}{\sqrt{t}}\,dt
=2\pi \int_{2z-x-y}^{\infty} p(t;-\xi)\,dt.
$$
As a result, we retrieve the famous reflection principle for Brownian motion:
$$
\P_x\{X(t)\le y\le z\le M(t)\}=\P\{X(t)\ge 2z-x-y\}.
$$
}

\bex{\textsl{Case $N=3$:}
we have to cases to distinguish.
\bitem
\item
Case $\k_3=+1$: the sum of interest in~\refp{distXM} reads here
$$
\frac{\t_0A_0B_1}{\t_1-\t_0}\,e^{\t_0\sqrt[3]{\l}(x-z)+\t_1\sqrt[3]{\l}(z-y)}
+\frac{\t_0A_0B_2}{\t_2-\t_0}\,e^{\t_0\sqrt[3]{\l}(x-z)+\t_2\sqrt[3]{\l}(z-y)}
$$
with $\frac{B_1}{\t_1-\t_0}=\frac{B_2}{\t_2-\t_0}=-\frac13$, and then
\beqa
F_{\l}(x,y,z)
&=&
\left\{\begin{array}{ll}
\dis\frac{2}{3\l}\bigg[e^{-\frac{\sqrt[3]{\l}}{2}(x-y)}
\cos\bigg(\frac{\sqrt3}{2}\,\sqrt[3]{\l}(x-y)\bigg)&
\\[2ex]
\dis -\,e^{\sqrt[3]{\l}(x+\frac12\,y-\frac32\,z)}
\cos\bigg(\frac{\sqrt3}{2}\,\sqrt[3]{\l}(z-y)\bigg)\bigg]
&\mbox{if $y\le x\le z$,}
\\[3ex]
\dis\frac{1}{\l}-\frac{1}{3\l}\bigg[e^{\sqrt[3]{\l}(x-y)}&
\\[2ex]
\dis +2\,e^{\sqrt[3]{\l}(x+\frac12\,y-\frac32\,z)}
\cos\bigg(\frac{\sqrt3}{2}\,\sqrt[3]{\l}(z-y)\bigg)\bigg]
&\mbox{if $x\le y\le z$.}
\end{array}\right.
\eeqa
We retrieve the results~(2.2) of~\cite{bor}.
Now, formula~\refp{dist-func-XM} writes
$$
\P_x\{X(t)\le y\le z\le M(t)\}
=\int_0^t\!\!\int_0^s p(\s;x-z)\,(a_{10}I_{10}+a_{20}I_{20})(s-\s;z-y)\,
\frac{ds\,d\s}{(t-s)^{2/3}}
$$
where, by~\refp{fourier-inverse},
$$
p(t;\xi)=\frac{1}{\pi} \int_0^{+\infty} \cos(\xi\l-t\l^3)\,d\l.
$$
The reciprocal relations, for $\xi\le 0$,
$$
\Phi(\l;\xi)=\frac{e^{\sqrt[3]{\l}\,\xi}}{3\l^{2/3}}
\quad\mbox{and}\quad
e^{\sqrt[3]{\l}\,\xi}=3\l^{2/3}\a_{00}\Phi(\l;\xi)
$$
imply that $\a_{00}=1$. Then
\beqa
a_{10}
&=&
\frac{3B_1}{\G(1/3)}\,\frac{A_0\a_{00}}{\t_0-\t_1}=\frac{1}{\G(1/3)},
\\
a_{20}
&=&
\frac{3B_2}{\G(1/3)}\,\frac{A_0\a_{00}}{\t_0-\t_2}=\frac{1}{\G(1/3)}.
\eeqa
Consequently,
$$
\P_x\{X(t)\le y\le z\le M(t)\}
=\frac{1}{\G(1/3)}\int_0^t\!\!\int_0^s p(\s;x-z)\,q(s-\s;z-y)\,
\frac{ds\,d\s}{(t-s)^{2/3}}
$$
with, for $\xi\ge 0$, by~\refp{Ik0},
\beqa
q(t;\xi)&=&(I_{10}+I_{20})(t;\xi)
\\
&=&
\frac{i\xi}{2\pi t}
\left[\t_1e^{\frac{i\pi}{3}} \int_0^{+\infty}
e^{-t\l^3+\t_1 e^{\frac{i\pi}{3}}\xi\l} \,d\l
-\t_1e^{-\frac{i\pi}{3}} \int_0^{+\infty}
e^{-t\l^3+\t_1 e^{-\frac{i\pi}{3}}\xi\l} \,d\l\right.
\\
&&
+\left.\t_2e^{\frac{i\pi}{3}} \int_0^{+\infty}
e^{-t\l^3+\t_2 e^{\frac{i\pi}{3}}\xi\l} \,d\l
-\t_2e^{-\frac{i\pi}{3}} \int_0^{+\infty}
e^{-t\l^3+\t_2 e^{-\frac{i\pi}{3}}\xi\l} \,d\l\right]
\\
&=&
-\frac{i\xi}{2\pi t}
\left[e^{\frac{i\pi}{3}} \int_0^{+\infty}
e^{-t\l^3+e^{\frac{i\pi}{3}}\xi\l} \,d\l
-e^{-\frac{i\pi}{3}} \int_0^{+\infty}
e^{-t\l^3+e^{-\frac{i\pi}{3}}\xi\l} \,d\l\right]
\\
&=&
\frac{\xi}{\pi t} \int_0^{+\infty} e^{-t\l^3+\frac12\,\xi\l}
\sin\bigg(\frac{\sqrt3}{2}\,\xi\l+\frac{\pi}{3}\bigg)\,d\l.
\eeqa

\item
Case $\k_3=-1$: the sum of interest in~\refp{distXM} reads here
$$
\frac{\t_0A_0B_1}{\t_1-\t_0}\,e^{\t_0\sqrt[3]{\l}(x-z)+\t_1\sqrt[3]{\l}(z-y)}
+\frac{\t_2A_2B_1}{\t_1-\t_2}\,e^{\t_2\sqrt[3]{\l}(x-z)+\t_1\sqrt[3]{\l}(z-y)}
$$ with $\frac{\t_0A_0}{\t_1-\t_0}=-\frac13\,e^{i\,\pi/3}$ and
$\frac{\t_2A_2}{\t_1-\t_2}=-\frac13\,e^{-i\,\pi/3}$, and then
\beqa
F_{\l}(x,y,z)
&=&
\left\{\begin{array}{l}
\dis\frac{1}{3\l}\bigg[e^{-\sqrt[3]{\l}(x-y)}-
2\,e^{\sqrt[3]{\l}(\frac12\,x+y-\frac32\,z)}
\cos\bigg(\frac{\sqrt3}{2}\,\sqrt[3]{\l}(x-z)+\frac{\pi}{3}\bigg)\bigg]
\\
\mbox{if $y\le x\le z$,}
\\[3ex]
\dis\frac{1}{\l}-\frac{1}{3\l}\bigg[2\,e^{\frac{\sqrt[3]{\l}}{2}(x-y)}
\cos\bigg(\frac{\sqrt3}{2}\,\sqrt[3]{\l}(x-y)\bigg)
\\[2ex]
\dis+\,e^{\frac{\sqrt[3]{\l}}{2}(\frac12\,x+y-\frac32\,z)}
\cos\bigg(\frac{\sqrt3}{2}\,\sqrt[3]{\l}(x-z)+\frac{\pi}{3}\bigg)\bigg]
\\
\mbox{if $x\le y\le z$.}
\end{array}\right.
\eeqa
We retrieve the results~(2.2) of~\cite{bor}.
Next, formula~\refp{dist-func-XM-bis} writes
\beqa
\lqn{\hspace{2.6em}
\P_x\{X(t)\le y\le z\le M(t)\}
}\\[-2ex]
&=&
\sum_{j\in \{0,2\}} b_{j0}\int_0^t\!\!\int_0^s
p(\s;z-y)\,I_{j0}(s-\s;x-z)\,\frac{ds\,d\s}{(t-s)^{2/3}}
\\
&=&
\int_0^t\!\!\int_0^s p(\s;z-y)\,(b_{00}I_{00}+b_{20}I_{20})(s-\s;x-z)\,
\frac{ds\,d\s}{(t-s)^{2/3}}
\eeqa
where, by~\refp{fourier-inverse},
$$
p(t;\xi)=\frac{1}{\pi} \int_0^{+\infty} \cos(\xi\l+t\l^3)\,d\l.
$$
From the reciprocal relations, which are valid for $\xi\ge 0$,
$$
\Phi(\l;\xi)=\frac{e^{-\sqrt[3]{\l}\,\xi}}{3\l^{2/3}}
\quad\mbox{and}\quad
-e^{-\sqrt[3]{\l}\,\xi}=3\l^{2/3}\b_{10}\Phi(\l;\xi)
$$
we extract the value $\b_{10}=-1$. Therefore,
\beqa
b_{00}
&=&
\frac{3\t_0A_0}{\G(1/3)}\,\frac{\bar{\t}_1B_1\b_{10}}{\t_1-\t_0}
=\frac{e^{\frac{2i\pi}{3}}}{\G(1/3)},
\\
b_{20}
&=&
\frac{3\t_2A_2}{\G(1/3)}\,\frac{\bar{\t}_1B_1\b_{10}}{\t_1-\t_2}
=\frac{e^{-\frac{2i\pi}{3}}}{\G(1/3)}.
\eeqa
Consequently,
$$
\P_x\{X(t)\le y\le z\le M(t)\}=\frac{1}{\G(1/3)}\int_0^t\!\!\int_0^s
p(\s;z-y)\,q(s-\s;x-z)\, \frac{ds\,d\s}{(t-s)^{2/3}}
$$
where, for $\xi\le 0$, by~\refp{Ik0},
\beqa
q(t;\xi)&=&(e^{\frac{2i\pi}{3}}I_{00}+e^{-\frac{2i\pi}{3}}I_{20})(t;\xi)
\\
&=&
\frac{i\xi}{2\pi t}
\left[-\t_0 \int_0^{+\infty}
e^{-t\l^3+\t_0 e^{\frac{i\pi}{3}}\xi\l} \,d\l
-\t_0e^{\frac{i\pi}{3}} \int_0^{+\infty}
e^{-t\l^3+\t_0 e^{-\frac{i\pi}{3}}\xi\l} \,d\l\right.
\\
&&
+\left.\t_2 e^{-\frac{i\pi}{3}} \int_0^{+\infty}
e^{-t\l^3+\t_2 e^{\frac{i\pi}{3}}\xi\l} \,d\l
+\t_2 \int_0^{+\infty}
e^{-t\l^3+\t_2 e^{-\frac{i\pi}{3}}\xi\l} \,d\l\right]
\\
&=&
\frac{\xi}{\pi t} \left[\sqrt3 \int_0^{+\infty} e^{-t\l^3+\xi\l}\,d\l
+\int_0^{+\infty} e^{-t\l^3-\frac12\,\xi\l}
\sin\bigg(\frac{\sqrt3}{2}\,\xi\l+\frac{\pi}{3}\bigg)\,d\l\right]\!.
\eeqa

\eitem
}

\bex{\textsl{Case $N=4$:}
in this case, we have
\beqa
\lqn{
\sum_{j\in J,k\in K} \frac{\t_jA_jB_k}{\t_k-\t_j}
\,e^{\t_j\sqrt[4]{\l}(x-z)+\t_k\sqrt[4]{\l}(z-y)}
}\\
&=&
\frac{\t_2A_2B_0}{\t_0-\t_2}\,e^{\t_2\sqrt[4]{\l}(x-z)+\t_0\sqrt[4]{\l}(z-y)}
+\frac{\t_2A_2B_1}{\t_1-\t_2}\,e^{\t_2\sqrt[4]{\l}(x-z)+\t_1\sqrt[4]{\l}(z-y)}
\\
&&
+\frac{\t_3A_3B_0}{\t_0-\t_3}\,e^{\t_3\sqrt[4]{\l}(x-z)+\t_0\sqrt[4]{\l}(z-y)}
+\frac{\t_3A_3B_1}{\t_1-\t_3}\,e^{\t_3\sqrt[4]{\l}(x-z)+\t_1\sqrt[4]{\l}(z-y)}
\eeqa
with
$$
\frac{\t_2A_2B_0}{\t_0-\t_2}=\frac i4,\quad
\frac{\t_3A_3B_1}{\t_1-\t_3}=-\frac i4,\quad
\frac{\t_2A_2B_1}{\t_1-\t_2}=-\frac{e^{-i\,\pi/4}}{2\sqrt2},\quad
\frac{\t_3A_3B_0}{\t_0-\t_3}=-\frac{e^{i\,\pi/4}}{2\sqrt2}.
$$
Hence,
\beqa
F_{\l}(x,y,z)
&=&
\left\{\begin{array}{l}
\dis\frac{1}{2\l}\bigg[e^{-\frac{\sqrt[4]{\l}}{\sqrt2}(x-y)}
\,\cos\bigg(\frac{\sqrt[4]{\l}}{\sqrt2}(x-y)\bigg)\!
-e^{\frac{\sqrt[4]{\l}}{\sqrt2}(x+y-2z)}
\\[2ex]
\dis\times\bigg(\cos\bigg(\frac{\sqrt[4]{\l}}{\sqrt2}(x-y)\bigg)\!
-\sin\bigg(\frac{\sqrt[4]{\l}}{\sqrt2}(x-y)\bigg)\!
-\sin\bigg(\frac{\sqrt[4]{\l}}{\sqrt2}(x+y-2z)\bigg)\bigg]
\\[2ex]
\mbox{if $y\le x\le z$,}
\\[3ex]
\dis\frac{1}{\l}-\frac{1}{2\l}\bigg[e^{\frac{\sqrt[4]{\l}}{\sqrt2}(x-y)}
\,\cos\bigg(\frac{\sqrt[4]{\l}}{\sqrt2}(x-y)\bigg)\!+
e^{\frac{\sqrt[4]{\l}}{\sqrt2}(x+y-2z)}
\\[2ex]
\dis\times\bigg(\cos\bigg(\frac{\sqrt[4]{\l}}{\sqrt2}(x-y)\bigg)\!
-\sin\bigg(\frac{\sqrt[4]{\l}}{\sqrt2}(x-y)\bigg)\!
-\sin\bigg(\frac{\sqrt[4]{\l}}{\sqrt2}(x+y-2z)\bigg)\bigg]
\\[2ex]
\mbox{if $x\le y\le z$.}
\end{array}\right.
\eeqa
We retrieve the results~(3.2) of~\cite{bor}.
Now, formula~\refp{dist-func-XM} writes
\beqa
\P_x\{X(t)\le y\le z\le M(t)\}
&\!\!\!=\!\!\!&\!
\int_0^t\!\!\int_0^s p(\s;x-z)\,(a_{00}I_{00}+a_{10}I_{10})(s-\s;z-y)\,
\frac{ds\,d\s}{(t-s)^{3/4}}
\\
&&
\!+\int_0^t\!\!\int_0^s \frac{\partial p}{\partial x}(\s;x-z)\,
(a_{01}I_{00}+a_{11}I_{10})(s-\s;z-y)\,\frac{ds\,d\s}{\sqrt{t-s}}.
\eeqa
where, by~\refp{fourier-inverse},
$$
p(t;\xi)=\frac{1}{\pi} \int_0^{+\infty} e^{-t \l^4}\cos(\xi\l) \,d\l.
$$
Let us consider the system
$$
\left\{\begin{array}{rcl}
\t_2\,e^{\t_2\sqrt[4]{\l}\,\xi}+\t_3\,e^{\t_3\sqrt[4]{\l}\,\xi}
&=&
4\l^{3/4}\Phi(\l;\xi)
\\[1ex]
\t_2^2\,e^{\t_2\sqrt[4]{\l}\,\xi}+\t_3^2\,e^{\t_3\sqrt[4]{\l}\,\xi}
&=&
\dis 4\sqrt{\l}\,\frac{\partial\Phi}{\partial \xi}(\l;\xi)
\end{array}
\right.
$$
which can be conversely written
$$
\left\{\begin{array}{rcl}
\t_2\,e^{\t_2\sqrt[4]{\l}\,\xi}
&=&
\dis\frac{4}{\t_3-\t_2}\bigg(\t_3\l^{3/4}\Phi(\l;\xi)-\sqrt{\l}
\,\frac{\partial\Phi}{\partial \xi}(\l;\xi)\bigg)
\\[2ex]
\t_3\,e^{\t_3\sqrt[4]{\l}\,\xi}
&=&
\dis\frac{4}{\t_3-\t_2}\bigg(-\t_2\l^{3/4}\Phi(\l;\xi)+\sqrt{\l}
\,\frac{\partial\Phi}{\partial \xi}(\l;\xi)\bigg)
\end{array}
\right.
$$
or, by means of the coefficients $\a_{20},\a_{21},\a_{30},\a_{31}$,
$$
\left\{\begin{array}{rcl}
\t_2\,e^{\t_2\sqrt[4]{\l}\,\xi}
&=&
\dis 4\bigg(\a_{20}\l^{3/4}\Phi(\l;\xi)+\a_{21}\sqrt{\l}\,
\frac{\partial\Phi}{\partial \xi}(\l;\xi)\bigg),
\\[2ex]
\t_3\,e^{\t_3\sqrt[4]{\l}\,\xi}
&=&
\dis 4\bigg(\a_{30}\l^{3/4}\Phi(\l;\xi)+\a_{31}\sqrt{\l}\,
\frac{\partial\Phi}{\partial \xi}(\l;\xi)\bigg).
\end{array}
\right.
$$
Identifying the two above systems yields the coefficients we are looking for:
\beqa
\a_{20}
&=&
\frac{\t_3}{\t_3-\t_2}=A_2=\frac{e^{-\frac{i\pi}{4}}}{\sqrt 2},
\\
\a_{30}
&=&
-\frac{\t_2}{\t_3-\t_2}=A_3=\frac{e^{\frac{i\pi}{4}}}{\sqrt 2},
\\
\a_{21}
&=&
-\frac{1}{\t_3-\t_2}=-\t_2A_2=\frac{i}{\sqrt 2},
\\
\a_{31}
&=&
\frac{1}{\t_3-\t_2}=-\t_3A_3=-\frac{i}{\sqrt 2},
\eeqa
and next:
\beqa
a_{00}
&=&
\frac{4B_0}{\G(1/4)}\left[\frac{A_2\a_{20}}{\t_2-\t_0}+
\frac{A_3\a_{30}}{\t_3-\t_0}\right]=\frac{1}{\sqrt2\,\G(1/4)},
\\
a_{10}
&=&
\frac{4B_1}{\G(1/4)}\left[\frac{A_2\a_{20}}{\t_2-\t_1}+
\frac{A_3\a_{30}}{\t_3-\t_1}\right]=\frac{1}{\sqrt2\,\G(1/4)},
\\
a_{01}
&=&
\frac{4B_0}{\G(1/2)}\left[\frac{A_2\a_{21}}{\t_2-\t_0}+
\frac{A_3\a_{31}}{\t_3-\t_0}\right]=\frac{e^{-\frac{i\pi}{4}}}{\sqrt{2\pi}},
\\
a_{11}
&=&
\frac{4B_1}{\G(1/2)}\left[\frac{A_2\a_{21}}{\t_2-\t_1}+
\frac{A_3\a_{31}}{\t_3-\t_1}\right]=\frac{e^{\frac{i\pi}{4}}}{\sqrt{2\pi}}.
\eeqa
Consequently,
\beqa
\P_x\{X(t)\le y\le z\le M(t)\}
&=&
\int_0^t\!\!\int_0^s p(\s;x-z)\,q_1(s-\s;z-y)\, \frac{ds\,d\s}{(t-s)^{3/4}}
\\
&&
+\int_0^t\!\!\int_0^s \frac{\partial p}{\partial x}(\s;x-z)\,q_2(s-\s;z-y)\,
\frac{ds\,d\s}{\sqrt{t-s}}
\eeqa
with, for $\xi\ge 0$, by~\refp{Ik0},
$$
$$
$$
q_1(t;\xi)=(a_{00}I_{00}+a_{10}I_{10})(t;\xi)
=\frac{1}{\sqrt2\,\G(1/4)}\,(I_{00}+I_{10})(t;\xi)
$$
and
\beqa
q_2(t;\xi)&=&(a_{01}I_{00}+a_{11}I_{10})(t;\xi)
=\frac{1}{\sqrt{2\pi}}\,(e^{-\frac{i\pi}{4}} I_{00}+ e^{\frac{i\pi}{4}} I_{10})(t;\xi)
\\
&=&
\frac{1}{2\pi}\,(I_{00}+I_{10})(t;\xi)
-\frac{i}{2\pi}\,(I_{00}-I_{10})(t;\xi).
\eeqa
Let us evaluate the intermediate quantities $(I_{00}\pm I_{10})(t;\xi)$:
\beqa
(I_{00}+I_{10})(t;\xi)
&=&
\frac{i\xi}{2\pi t}
\left[\t_0 e^{\frac{i\pi}{4}} \int_0^{+\infty}
e^{-t\l^4+\t_0 e^{\frac{i\pi}{4}}\xi\l} \,d\l
-\t_0 e^{-\frac{i\pi}{4}} \int_0^{+\infty}
e^{-t\l^4+\t_0 e^{-\frac{i\pi}{4}}\xi\l} \,d\l\right.
\\
&&
+\left.\t_1 e^{\frac{i\pi}{4}} \int_0^{+\infty}
e^{-t\l^4+\t_1 e^{\frac{i\pi}{4}}\xi\l} \,d\l
-\t_1 e^{-\frac{i\pi}{4}} \int_0^{+\infty}
e^{-t\l^4+\t_1 e^{-\frac{i\pi}{4}}\xi\l} \,d\l\right]
\\
&=&
\frac{\xi}{\pi t} \int_0^{+\infty} e^{-t\l^4}\,\cos(\xi\l)\,d\l
\eeqa
and
\beqa
(I_{00}-I_{10})(t;\xi)
&=&
-\frac{i\xi}{\pi t} \left[ \int_0^{+\infty} e^{-t\l^4-\xi\l}\,d\l
-\int_0^{+\infty} e^{-t\l^4}\sin(\xi\l)\,d\l\right]\!.
\eeqa
Then
\beqa
q_1(t;\xi)
&=&
\frac{\xi}{\pi\sqrt2\,\G(1/4)\,t} \int_0^{+\infty} e^{-t\l^4}\,\cos(\xi\l)\,d\l,
\\
q_2(t;\xi)
&=&
\frac{\xi}{2\pi^2\,t} \int_0^{+\infty} e^{-t\l^4}\,
\left(\cos(\xi\l)+\sin(\xi\l)-e^{-\xi\l}\right)\,d\l.
\eeqa
}


\subsection{Boundary value problem}

In this part, we show that the function $x\longmapsto F_{\l}(x,y,z)$ solves a
boundary value problem
related to the differential operator $\cD_x=\k_N\,\frac{d^N}{d x^N}$.
Fix $y<z$ and set $F(x)=F_{\l}(x,y,z)$ for $x\in(-\infty,z]$.
%
\bpr{
The function $F$ satisfies the differential equation
\beq\label{ODE-F}
\cD_xF(x)=\left\{\begin{array}{ll}
\l F(x)-1 &\mbox{for $x\in(-\infty,y)$,}
\\[1ex]
\l F(x) &\mbox{for $x\in(y,z)$,}
\end{array}\right.
\eeq
together with the conditions
\beqan
F^{(l)}(z^-)&=0&\mbox{for $0\le l\le \card J-1$},
\label{boundary-cond1-F}
\\[1ex]
F^{(l)}(y^+)-F^{(l)}(y^-)&=0&\mbox{for $0\le l \le N-1$.}
\label{boundary-cond2-F}
\eeqan
}
%
\dem
The differential equation~\refp{ODE-F} is readily obtained by
differentiating~\refp{distXM} with
respect to $x$. Let us derive the boundary condition~\refp{boundary-cond1-F}:
\beqa
F^{(l)}(z^-)
&=&
\frac{\l^{l/N}}{\l} \sum_{j\in J,k\in K} \frac{\t_j^{l+1}A_jB_k}{\t_k-\t_j}
\,e^{\t_k\!\!\sqrt[N]{\l}(z-y)}
+\frac{\l^{l/N}}{N\l} \sum_{k\in K} \t_k^{l+1}\,e^{\t_k\!\!\sqrt[N]{\l}(z-y)}
\\
&=&
\l^{l/N-1} \sum_{k\in K}\Bigg[\sum_{j\in J} \frac{\t_j^{l+1}A_j}{\t_k-\t_j}
\,B_k+\frac{\t_k^{l+1}}{N}\Bigg]\,e^{\t_k\!\!\sqrt[N]{\l}(z-y)}=0
\eeqa
where the last equality comes from~\refp{expansion} with $x=\t_k$ which yields
$\dis \sum_{j\in J}\frac{\t_j^{l+1}A_j}{\t_k-\t_j}=-\frac{\t_k^{l+1}}{NB_k}.$
Condition~\refp{boundary-cond2-F} is quite easy to check.
\fin
%
\brem{
Condition~\refp{boundary-cond2-F} says that the function $F$ is regular up to
the order $n-1$. It can also be easily seen that
$F^{(N)}(y^+)-F^{(N)}(y^-)=\k_N$ which says that the function $F^{(N)}$
has a jump at point $y$.
On the other hand, the boundary value
problem~\refp{ODE-F}--\refp{boundary-cond1-F}--\refp{boundary-cond2-F}
(the differential equation together with the $N+\card J$ conditions)
augmented of a boundedness condition on $(-\infty,y)$
may be directly solved by using Vandermonde determinants.
}


\section{Distributions of $(\tau_a^+,X(\tau_a^+))$ and $(\tau_a^-,X(\tau_a^-))$}
\label{sect-distrib-(tau,Xtau)}

The integer $N$ is again assumed to be even.
Recall we set $\tau_a^+=\inf\{t\ge0:X(t)>a\}$ and $\tau_a^-=\inf\{t\ge0:X(t)<a\}$.
The aim of this section is to derive the distributions of the vectors
$(\tau_a^+,X(\tau_a^+))$ and $(\tau_a^-,X(\tau_a^-))$.
For this, we proceed in three steps: we first compute the Laplace-Fourier
transform of, \eg\ $(\tau_a^+,X(\tau_a^+))$
(Subsection~\ref{subsect-LFT-tau-Xtau});
we next invert the Fourier transform (with respect to $\mu$,
Subsubsection~\ref{subsubsect-invert-mu-tau-Xtau}) and we finally invert the
Laplace transform (with respect to $\l$,
Subsubsection~\ref{subsubsect-invert-lambda-tau-Xtau}).
We have especially obtained a remarkable formula for the densities
of $X(\tau_a^+)$ and $X(\tau_a^-)$ by means of multipoles
(Subsection~\ref{subsect-FT-Xtau}).

\subsection{Laplace-Fourier transforms}\label{subsect-LFT-tau-Xtau}

We have a relationship between the distributions of $(\tau_a^+,X(\tau_a^+))$
and $(X(t),M(t))$, and between those of $(\tau_a^-,X(\tau_a^-))$ and $(X(t),m(t))$.
%
\blm{\label{lemma-relation}
The Laplace-Fourier transforms of the vectors $(\tau_a^+,X(\tau_a^+))$ and
$(\tau_a^-,X(\tau_a^-))$ are related to the distributions of the vectors
$(X(t),M(t))$ and $(X(t),m(t))$ according as, for $\Re(\l)>0$ and $\mu\in\R$,
\beqan
\Ex\!\!\left[e^{-\l\tau_a^++i\mu X(\tau_a^+)}\right]\!\!
&=\;
\dis\left(\l-\k_N(i\mu)^N\right) \int_0^{+\infty}
e^{-\l t}\,\Ex\!\!\left[e^{i\mu X(t)},M(t)>a\right]dt
&\mbox{for $x\le a$,}
\nonumber\\[-1ex]
\label{LFT-tau-Xtau-inter}\\[-1ex]
\Ex\!\!\left[e^{-\l\tau_a^-+i\mu X(\tau_a^-)}\right]\!\!
&=\;
\dis\left(\l-\k_N(i\mu)^N\right) \int_0^{+\infty}
e^{-\l t}\,\Ex\!\!\left[e^{i\mu X(t)},m(t)<a\right]dt
&\mbox{for $x\ge a$.}
\nonumber
\eeqan
}
%
\dem
We divide the proof of Lemma~\ref{lemma-relation} into five steps.

\noindent\textbf{$\bullet$ Step 1}

For the step-process $(X_n(t))_{t\ge 0}$, the corresponding first hitting time
$\tauan^+$ is the instant $\tnk$ with $k$ such that
$X(\tnj)\le a$ for all $j\in\{0,\ldots,k-1\}$ and $X(\tnk)> a$, or, equivalently,
such that $M_{n,k-1}\le a$ and $M_{n,k}> a$ where $M_{n,k}=\max_{0\le j\le k}X_{n,j}$
and $X_{n,k}=X(\tnk)$ for $k\ge 0$ and $M_{n,-1}=-\infty$. We have, for $x\le a$,
\beqan
e^{-\l\tauan^++i\mu X_n(\tauan^+)}
&=&
\sum_{k=0}^{\infty} e^{-\l \tnk+i\mu X_{n,k}}
\ind_{\{M_{n,k-1}\le a< M_{n,k}\}}
\nonumber\\
&=&
\sum_{k=0}^{\infty} e^{-\l \tnk+i\mu X_{n,k}}\!
\left[\ind_{\{M_{n,k}>a\}}-\ind_{\{M_{n,k-1}>a\}}\right]\!.
\label{sum-Abel}
\eeqan
Let us apply classical Abel's identity to sum~\refp{sum-Abel}.
This yields, since $\ind_{\{M_{n,-1}>a\}}=0$ and
$\lim_{k\to +\infty} e^{-\l \tnk+i\mu X_{n,k}}\ind_{\{M_{n,k}>a\}}=0$,
for $\Re(\l)>0$:
$$
e^{-\l\tauan^++i\mu X_n(\tauan^+)}
=\sum_{k=0}^{\infty} \left[e^{-\l \tnk+i\mu X_{n,k}}
-e^{-\l \tnkun+i\mu X_{n,k+1}}\right] \!\ind_{\{M_{n,k}>a\}}.
$$
The functional $e^{-\l\tauan^++i\mu X_n(\tauan^+)}$ is a function
of discrete observations of $X$.

\noindent\textbf{$\bullet$ Step 2}

In order to evaluate the expectation of the foregoing functional,
we need to check that the series
$$
\sum_{k=0}^{\infty} \Ex\!\!\left[\left(e^{-\l \tnk+i\mu X_{n,k}}
-e^{-\l \tnkun+i\mu X_{n,k+1}}\right)\!\ind_{\{M_{n,k}>a\}}\right]
$$
is absolutely convergent. For this, we use the Markov property
and derive the following estimate:
\beqa
\lqn{
\left|\Ex\!\!\left[e^{-\l \tnk+i\mu X_{n,k}}\ind_{\{M_{n,k}\le a\}}\right]\right|
=\left|e^{-\l \tnk} \Ex\!\!\left[e^{i\mu X_{n,k}}
\ind_{\{X_{n,1}\le a,\ldots,X_{n,k}\le a\}}\right]\right|
}
\\
&\le&
|e^{-\l/2^n}|^k \left|\int_{-\infty}^a\ldots\int_{-\infty}^a e^{i\mu x_k}
p(1/2^n;x-x_1)\cdots p(1/2^n;x_{k-1}-x_k)\,dx_1\cdots dx_k\right|
\\
&\le&
(\rho\,e^{-\Re(\l)/2^n})^k.
\eeqa
We recall that in the last inequality
$\rho=\int_{-\infty}^{+\infty} |p(t;z)|\,dz<+\infty$.
Similar computations yield the inequality
$\left|\Ex\!\!\left[e^{-\l \tnk+i\mu X_{n,k}}\right]\right|\le
(\rho\,e^{-\Re(\l)/2^n})^k.$
Because of the identity $\ind_{\{M_{n,k}>a\}}=1-\ind_{\{M_{n,k}\le a\}}$, we
plainly get
$$
\left|\Ex\!\!\left[e^{-\l \tnk+i\mu X_{n,k}} \ind_{\{M_{n,k}>a\}}\right]\right|
\le 2(\rho\,e^{-\Re(\l)/2^n})^k.
$$
Upon adding one integral more in the above discussion, it is easily seen that
$$
\left|\Ex\!\!\left[e^{-\l \tnk+i\mu X_{n,k+1}} \ind_{\{M_{n,k}>a\}}\right]\right|
\le 2(\rho\,e^{-\Re(\l)/2^n})^{k+1}.
$$
As a result, when choosing $\l$ such that $\Re(\l)>2^n \ln \rho$, we have
$$
\sum_{k=0}^{\infty} \left|\Ex\!\!\left[\left(e^{-\l \tnk+i\mu X_{n,k}}
-e^{-\l \tnkun+i\mu X_{n,k+1}}\right)\!\ind_{\{M_{n,k}>a\}}\right]\right|
\le \frac{2(1+\rho\,e^{-\Re(\l)/2^n})}{1-\rho\,e^{-\Re(\l)/2^n}}<+\infty.
$$
\noindent\textbf{$\bullet$ Step 3}

Therefore, we can evaluate the expectation of
$e^{-\l\tauan^++i\mu X_n(\tauan^+)}$.
By the Markov property we get, for $\Re(\l)>2^n \ln \rho$,
\beqa
\lqn{
\Ex\!\!\left[e^{-\l\tauan^++i\mu X_n(\tauan^+)}\right]
}\\[-2ex]
&=&
\sum_{k=0}^{\infty} e^{-\l \tnk} \Ex\!\!\left[e^{i\mu X_{n,k}}\ind_{\{M_{n,k}>a\}}
\left(1-e^{-\l/2^n}\,e^{i\mu (X_{n,k+1}-X_{n,k})}\right)\right]
\\
&=&
\sum_{k=0}^{\infty} e^{-\l \tnk} \Ex\!\!\left[e^{i\mu X_{n,k}} \ind_{\{M_{n,k}>a\}}
\left(1-e^{-\l/2^n}\,\E_{X_{n,k}}\!\!\left(e^{i\mu (X_{n,1}-X_{n,0})}\right)\right) \right]\!.
\eeqa
Since
$\E_{X_{n,k}}\!\!\left(e^{i\mu (X_{n,1}-X_{n,0})}\right)=e^{\k_N(i\mu)^N/2^n}$
we obtain, for $\Re(\l)>2^n \ln \rho$,
\beqan
\Ex\!\!\left[e^{-\l\tauan^++i\mu X_n(\tauan^+)}\right]
&=&
2^n\!\left(1-e^{-\left(\l-\k_N(i\mu)^N\right)/2^n}\right)
\nonumber\\
&&
\times\frac{1}{2^n}\sum_{k=0}^{\infty} e^{-\l \tnk}
\Ex\!\!\left[e^{i\mu X_{n,k}}\ind_{\{M_{n,k}>a\}}\right]\!.
\label{LFT-tauplus-Xtauplus}
\eeqan
\noindent\textbf{$\bullet$ Step 4}

In order to take the limit of~\refp{LFT-tauplus-Xtauplus} as $n$ tends to
$\infty$, we have to check the validity of~\refp{LFT-tauplus-Xtauplus} for any
$\l$ such that $\Re(\l)>0$. For this, we first consider its Laplace transform
with respect to $a$:
\beqa
\lqn{
\int_x^{+\infty} e^{-\nu a}
\Ex\!\!\left[e^{-\l\tauan^++i\mu X_n(\tauan^+)}\right]da
}
&=&
\left(1-e^{-\left(\l-\k_N(i\mu)^N\right)/2^n}\right)
\sum_{k=0}^{\infty} e^{-\l \tnk}
\int_x^{+\infty} e^{-\nu a} \Ex\!\!\left[e^{i\mu X_{n,k}}
\ind_{\{M_{n,k}>a\}}\right]da.
\eeqa
The sum in the above equality writes, using the settings of
Subsection~\ref{subsect-LFT-XMm},
\beqa
\lqn{
\sum_{k=0}^{\infty} e^{-\l \tnk}
\int_x^{+\infty} e^{-\nu a} \Ex\!\!\left[e^{i\mu X_{n,k}}
\ind_{\{M_{n,k}>a\}}\right]da
}
&=&
\sum_{k=0}^{\infty} e^{-\l \tnk}
\Ex\!\!\left[e^{i\mu X_{n,k}}\int_x^{M_{n,k}} e^{-\nu a}\,da\right]
\\
&=&
\sum_{k=0}^{\infty} \frac{e^{-\l \tnk}}{\nu}
\left[e^{-\nu x}\Ex\!\!\left(e^{i\mu X_{n,k}}\right)-
\Ex\!\! \left(e^{i\mu X_{n,k}-\nu M_{n,k}} \right)\right]
\\
&=&
\frac{1}{\nu}\Bigg[e^{(i\mu-\nu)x}\sum_{k=0}^{\infty} e^{-(\l-\k_N(i\mu)^N) \tnk}
-\sum_{k=0}^{\infty} e^{-\l \tnk}
\Ex\!\! \left(e^{i\mu X_{n,k}-\nu M_{n,k}} \right)\Bigg]
\\
&=&
\frac{1}{\nu}\Bigg[\frac{e^{(i\mu-\nu)x}}{1-e^{-(\l-\k_N(i\mu)^N)/2^n}}
-\frac{\l}{1-e^{-\l/2^n}} \,\Ex\!\!\left(F_{X_n}^+(\l,\mu,\nu)\right)\Bigg]
\\
&=&
\frac{e^{(i\mu-\nu)x}}{\nu}\left[\frac{1}{1-e^{-(\l-\k_N(i\mu)^N)/2^n}}
-\frac{1}{1-e^{-\l/2^n}} \,\exp\left(\frac{1}{2^n}\sum_{k=1}^{\infty}
\frac{e^{-\l\tnk}}{\tnk}\,\psi^+(\mu,\nu;\tnk)\right)
\right]\!.
\eeqa
We then obtain
\beqa
\lqn{
\int_x^{+\infty} e^{-\nu a}\Ex\!\!\left[e^{-\l\tauan^++i\mu X_n(\tauan^+)}\right]da
}
&=&
\frac{e^{(i\mu-\nu)x}}{\nu}\left[1-\frac{1-e^{-(\l-\k_N(i\mu)^N)/2^n}}{1-e^{-\l/2^n}}
\,\exp\left(\frac{1}{2^n}\sum_{k=1}^{\infty}
\frac{e^{-\l\tnk}}{\tnk}\,\psi^+(\mu,\nu;\tnk)\right)
\right]\!.
\eeqa
Inverting the Laplace transform yields, noting that the function
$a\longmapsto \Ex\!\!\left[e^{-\l\tauan^++i\mu X_n(\tauan^+)}\right]$
is right-continuous,
\beqa
\lqn{
\Ex\!\!\left[e^{-\l\tauan^++i\mu X_n(\tauan^+)}\right]
=\frac{1}{2i\pi} \lim_{\e\to 0^+} \int_{c-i\infty}^{c+i\infty}
e^{(i\mu-\nu)x+\nu (a+\e)}}\\[1ex]
&&
\times\left[1-\frac{1-e^{-(\l-\k_N(i\mu)^N)/2^n}}{1-e^{-\l/2^n}}
\,\exp\left(\frac{1}{2^n}\sum_{k=1}^{\infty}
\frac{e^{-\l\tnk}}{\tnk}\,\psi^+(\mu,\nu;\tnk)\right)
\right]\!\frac{d\nu}{\nu}.
\eeqa
Putting
$$
\psi^+(\mu,\nu;t)=\psi_1(i\mu;t)+\psi_2(i\mu-\nu;t)
$$
with
\beqa
\psi_1(\a;t)=\Eo\!\!\left[\left(e^{\a X(t)}-1\right)\ind_{\{X(t)<0\}}\right],
\quad \psi_2(\a;t)=\Eo\!\!\left[\left(e^{\a X(t)}-1\right)\ind_{\{X(t)\ge 0\}}\right]\!,
\eeqa
the exponential within the last displayed integral writes
\beqa
\lqn{\exp\left(\frac{1}{2^n}\sum_{k=1}^{\infty}
\frac{e^{-\l\tnk}}{\tnk}\,\psi^+(\mu,\nu;\tnk)\right)
}
&=&
\exp\left(\frac{1}{2^n}\sum_{k=1}^{\infty}
\frac{e^{-\l\tnk}}{\tnk}\,\psi_1(i\mu;\tnk)\right)
\exp\left(\frac{1}{2^n}\sum_{k=1}^{\infty}
\frac{e^{-\l\tnk}}{\tnk}\,\psi_2(i\mu-\nu;\tnk)\right)\!.
\eeqa
Noticing that
$$
\frac{1}{2i\pi}\int_{c-i\infty}^{c+i\infty}
e^{(i\mu-\nu)x+\nu a}\frac{d\nu}{\nu}=e^{i\mu x},
$$
we get
\beqa
\Ex\!\!\left[e^{-\l\tauan^++i\mu X_n(\tauan^+)}\right]
&\!\!\!=\!\!&
e^{i\mu x}\left[1-\frac{1-e^{-(\l-\k_N(i\mu)^N)/2^n}}{1-e^{-\l/2^n}}
\exp\left(\frac{1}{2^n}\sum_{k=1}^{\infty}
\frac{e^{-\l\tnk}}{\tnk}\,\psi_1(i\mu;\tnk)\right)\right.
\\
&&
\left.\times\frac{1}{2i\pi}\int_{c-i\infty}^{c+i\infty} e^{(a-x)\nu}
\exp\left(\frac{1}{2^n}\sum_{k=1}^{\infty}
\frac{e^{-\l\tnk}}{\tnk}\,\psi_2(i\mu-\nu;\tnk)\right)
\!\frac{d\nu}{\nu}\right]\!.
\eeqa

By imitating the method used by Nishioka (Appendix in~\cite{nish2})
for deriving subtil extimates, it may be seen that this last expression
is bounded over the half-plane $\Re(\l)\ge \e$ for any $\e>0$.
Hence, as in the proof of the validity of~\refp{expectation-Fplus} for
$\Re(\l)>0$, we see that~\refp{LFT-tauplus-Xtauplus} is also valid for $\Re(\l)>0$.
It follows that the functional
$e^{-\l\tau_a^++i\mu X(\tau_a^+)}$ is admissible.

\noindent\textbf{$\bullet$ Step 5}

Now, we can let $n$ tend to $+\infty$ in~\refp{LFT-tauplus-Xtauplus}.
For $\Re(\l)>0$, we obviously have
$$
\lim_{n\to +\infty} 2^n\!\left(1-e^{-\left(\l-\k_N(i\mu)^N\right)/2^n}\right)
=\l-\k_N(i\mu)^N,
$$
and we finally obtain the relationship~\refp{LFT-tau-Xtau-inter} corresponding
to $\tau_a^+$. The proof of that corresponding to $\tau_a^-$ is quite similar.
\fin

\bth{
The Laplace-Fourier transforms of the vectors $(\tau_a^+,X(\tau_a^+))$ and
$(\tau_a^-,X(\tau_a^-))$ are determined, for $\Re(\l)>0$ and $\mu\in\R$, by
\beqan
\Ex\!\!\left[e^{-\l\tau_a^++i\mu X(\tau_a^+)}\right]\!\!
&=\;
\dis\sum_{j\in J} A_j\prod_{l\in J\setminus\{j\}}
\left(1-\frac{i\mu}{\sqrt[N]{\l}}\,\bar{\t}_l\right)
e^{\t_j\!\!\sqrt[N]{\l}\,(x-a)} \,e^{i\mu a}
&\mbox{for $x\le a$,}
\nonumber\\[-1ex]
\label{LFT-tau-Xtau}\\[-1ex]
\Ex\!\!\left[e^{-\l\tau_a^-+i\mu X(\tau_a^-)}\right]\!\!
&=\;
\dis\sum_{k\in K} B_k\prod_{l\in K\setminus\{k\}}
\left(1-\frac{i\mu}{\sqrt[N]{\l}}\,\bar{\t}_l\right)
e^{\t_k\!\!\sqrt[N]{\l}\,(x-a)} \,e^{i\mu a}
&\mbox{for $x\ge a$.}
\nonumber
\eeqan
}
%
\dem
Using~\refp{inverting-nu-XMm} gives
\beqan
\lqn{
\int_0^{+\infty} e^{-\l t}\,\Ex\!\!\left[e^{i\mu X(t)},M(t)>a\right]dt
}&=&
\frac{\l^{(1-\card J)/N}\,e^{i\mu x}}
{\dis\mathop{\textstyle{\prod}}_{k\in K} (\!\sqrt[N]{\l}-i\mu\t_k)}
\sum_{j\in J} \t_jA_j \int_a^{+\infty} e^{(i\mu-\t_j\!\!\sqrt[N]{\l}\,)(z-x)}\,dz.
\label{integ-M>a}
\eeqan
Plugging the following equality
$$
\l-\k_N(i\mu)^N=\prod_{l=0}^{N-1} (\!\sqrt[N]{\l}-i\mu\t_l)
=\prod_{j\in J} (\!\sqrt[N]{\l}-i\mu\t_j)\times
\prod_{k\in K} (\!\sqrt[N]{\l}-i\mu\t_k)
$$
into~\refp{integ-M>a} and remarking that the set $\{\t_j,j\in J\}$
is invariant by conjugating yield
\beqan
\lqn{
\int_0^{+\infty} e^{-\l t}\,\Ex\!\!\left[e^{i\mu X(t)},M(t)>a\right]dt
}&=&
\l^{(1-\card J)/N} e^{i\mu x}\,\frac{\dis\mathop{\textstyle{\prod}}_{j\in J}
(\!\sqrt[N]{\l}-i\mu\t_j)}{\l-\k_N(i\mu)^N}
\sum_{j\in J} \frac{A_j}{\!\sqrt[N]{\l}-i\mu\bar{\t}_j} \,
e^{-(i\mu-\t_j\!\!\sqrt[N]{\l}\,)(x-a)}
\nonumber\\
&=&
\frac{e^{i\mu x}}{\l-\k_N(i\mu)^N}
\sum_{j\in J} A_j \prod_{l\in J\setminus\{j\}}
\left(1-\frac{i\mu}{\sqrt[N]{\l}}\,\bar{\t}_l\right)
e^{-(i\mu-\t_j\!\!\sqrt[N]{\l}\,)(x-a)}.
\label{LFT-XM-inter}
\eeqan
Consequently, by putting~\refp{LFT-XM-inter} into~\refp{LFT-tau-Xtau-inter},
we obtain~\refp{LFT-tau-Xtau}.
\fin

\brem{
Choosing $\mu=0$ in~\refp{LFT-tau-Xtau} supplies the Laplace transforms of
$\tau_a^+$ and $\tau_a^-$:
\beqa
\Ex\!\!\left[e^{-\l\tau_a^+}\right]\!\!
&=\;
\dis\sum_{j\in J} A_j \,e^{\t_j\!\!\sqrt[N]{\l}\,(x-a)}
&\mbox{for $x\le a$,}
\\
\Ex\!\!\left[e^{-\l\tau_a^-}\right]\!\!
&=\;
\dis\sum_{k\in K} B_k \,e^{\t_k\!\!\sqrt[N]{\l}\,(x-a)}
&\mbox{for $x\ge a$.}
\eeqa
}
%
\brem{
An alternative method for deriving the distribution of
$(\tau_a^+,X(\tau_a^+))$ consists of computing the joint distribution
of $\left(X(t),\ind_{(-\infty,a)}(M(t))\right)$ instead of that of
$(X(t),M(t))$ and next to invert a certain Fourier transform.
This way was employed by Nishioka~\cite{nish2} in the case $N=4$
and may be applied to the general case \textit{mutatis mutandis}.
}
\brem{
The following relationship issued from fluctuation theory holds for Levy
processes: if $x\le a$,
\beq\label{fluct}
\Ex\!\!\left[e^{-\l\tau_a^-+i\mu X(\tau_a^-)}\right]\!
=e^{i\mu a}\,\frac{\int_0^{+\infty} e^{-\l t}\,
\Ex\!\!\left[e^{i\mu (M(t)-a)},M(t)\ge a\right]dt}
{\int_0^{+\infty} e^{-\l t}\,\Eo\!\!\left[e^{i\mu M(t)}\right]dt}.
\eeq
Let us check that~\refp{fluct} also holds, at least formally,
for the pseudo-process $X$. We have, by~\refp{densities-Mm},
\beqan
\lqn{\hspace{2.4em}
\int_0^{+\infty} e^{-\l t}\,\Ex\!\!\left[e^{i\mu (M(t)-a)},M(t)\ge a\right]dt
}\nonumber\\
&=&
\int_a^{+\infty} e^{i\mu(z-a)}\,dz \int_0^{+\infty}
e^{-\l t}\,dt\,\P_x\{M(t)\in dz\}/dz
\nonumber\\
&=&
\int_a^{+\infty} \l^{1/N-1} \sum_{j\in J}\t_jA_j\,
e^{i\mu(z-a)-\t_j\!\!\sqrt[N]{\l}\,(z-x)}\,dz
\nonumber\\
&=&
\frac{1}{\l} \sum_{j\in J} \frac{\t_jA_j}{\t_j-\frac{i\mu}{\sqrt[N]{\l}}}
\,e^{\t_j\!\!\sqrt[N]{\l}\,(x-a)}.
\label{fluct1}
\eeqan
For $x=a$, this yields, by~\refp{expansion},
\beq\label{fluct2}
\int_0^{+\infty} e^{-\l t}\,\Eo\!\!\left[e^{i\mu M(t)}\right]dt
=\frac{1}{\l} \sum_{j\in J} \frac{\t_jA_j}{\t_j-\frac{i\mu}{\sqrt[N]{\l}}}
=\frac{1}{\l}\Bigg[\prod_{j\in J}
\bigg(1-\frac{i\mu}{\sqrt[N]{\l}}\,\bar{\t}_j\bigg)\Bigg]^{-1}.
\eeq
As a result, by plugging~\refp{fluct1} and~\refp{fluct2} into~\refp{fluct},
we retrieve~\refp{LFT-tau-Xtau}.
}

\bex{\textsl{Case $N=2$:} we simply have
\beqa
\Ex\!\!\left[e^{-\l\tau_a^++i\mu X(\tau_a^+)}\right]
&=\;
e^{i\mu a+\sqrt{\l}\,(x-a)}
&\mbox{for $x\le a$,}
\\
\Ex\!\!\left[e^{-\l\tau_a^-+i\mu X(\tau_a^-)}\right]
&=\;
e^{i\mu a-\sqrt{\l}\,(x-a)}
&\mbox{for $x\ge a$.}
\eeqa
}
\bex{\textsl{Case $N=3$:}
\bitem
\item
In the case $\k_3=+1$, we have, for $x\le a$,
$$
\sum_{j\in J} A_j\prod_{l\in J\setminus\{j\}}
\left(1-\frac{i\mu}{\sqrt[3]{\l}}\,\bar{\t}_l\right)
e^{\t_j\!\sqrt[3]{\l}\,(x-a)}
=A_0\,e^{\t_0\!\sqrt[3]{\l}\,(x-a)}
=e^{\!\sqrt[3]{\l}\,(x-a)},
$$
and, for $x\ge a$,
\beqa
\lqn{\hspace{2.4em}
\sum_{k\in K} B_k\prod_{l\in K\setminus\{k\}}
\left(1-\frac{i\mu}{\sqrt[3]{\l}}\,\bar{\t}_l\right)
e^{\t_k\!\sqrt[3]{\l}\,(x-a)}
}\\
&=&
B_1\left(1-\frac{i\mu}{\sqrt[3]{\l}}\,\bar{\t}_2\right)
e^{\t_1\!\sqrt[3]{\l}\,(x-a)}
+B_2\left(1-\frac{i\mu}{\sqrt[3]{\l}}\,\bar{\t}_1\right)
e^{\t_2\!\sqrt[3]{\l}\,(x-a)}
\\
&=&
\frac{1}{\sqrt3}\,e^{-\frac12\sqrt[3]{\l}\,(x-a)}\left[
\left(e^{-i\pi/6}+\frac{\mu}{\sqrt[3]{\l}}\right)
e^{i\frac{\sqrt3}{2}\sqrt[3]{\l}\,(x-a)}
+\left(e^{i\pi/6}-\frac{\mu}{\sqrt[3]{\l}}\right)
e^{-i\frac{\sqrt3}{2}\sqrt[3]{\l}\,(x-a)}\right]\!.
\eeqa
Therefore, \refp{LFT-tau-Xtau} writes
\beqa
\Ex\!\!\left[e^{-\l\tau_a^++i\mu X(\tau_a^+)}\right]
&=&
e^{i\mu a+\sqrt[3]{\l}\,(x-a)}
\quad\mbox{for $x\le a$,}
\\
\Ex\!\!\left[e^{-\l\tau_a^-+i\mu X(\tau_a^-)}\right]
&=&
\frac{2}{\sqrt3}\,e^{i\mu a-\frac12\sqrt[3]{\l}\,(x-a)}\left[
\cos\bigg(\frac{\sqrt3}{2}\sqrt[3]{\l}\,(x-a)-\frac{\pi}{6}\bigg)\right.
\\
&&
+\left.\frac{i\mu}{\sqrt[3]{\l}} \sin\bigg(\frac{\sqrt3}{2}\sqrt[3]{\l}\,(x-a)\bigg)
\right]\quad\mbox{for $x\ge a$.}
\eeqa

\item
In the case $\k_3=-1$, we similarly have that
\beqa
\Ex\!\!\left[e^{-\l\tau_a^++i\mu X(\tau_a^+)}\right]
&=&
\frac{2}{\sqrt3}\,e^{i\mu a+\frac12\sqrt[3]{\l}\,(x-a)}\left[
\cos\bigg(\frac{\sqrt3}{2}\sqrt[3]{\l}\,(x-a)+\frac{\pi}{6}\bigg)\right.
\\
&&
-\left.\frac{i\mu}{\sqrt[3]{\l}} \sin\bigg(\frac{\sqrt3}{2}\sqrt[3]{\l}\,(x-a)\bigg)
\right]\quad\mbox{for $x\le a$,}
\\
\Ex\!\!\left[e^{-\l\tau_a^-+i\mu X(\tau_a^-)}\right]
&=&
e^{i\mu a-\sqrt[3]{\l}\,(x-a)}\quad\mbox{for $x\ge a$.}
\eeqa
\eitem
}

\bex{\textsl{Case $N=4$:} we have, for $x\le a$,
\beqa
\lqn{
\sum_{j\in J} A_j\prod_{l\in J\setminus\{j\}}
\left(1-\frac{i\mu}{\sqrt[4]{\l}}\,\bar{\t}_l\right)
e^{\t_j\!\sqrt[4]{\l}\,(x-a)}
}\\
&=&
A_2\left(1-\frac{i\mu}{\sqrt[4]{\l}}\,\bar{\t}_3\right)
e^{\t_2\!\sqrt[4]{\l}\,(x-a)}
+A_3\left(1-\frac{i\mu}{\sqrt[4]{\l}}\,\bar{\t}_2\right)
e^{\t_3\!\sqrt[4]{\l}\,(x-a)}
\\
&=&
\frac{1}{\sqrt2}\,e^{\frac{1}{\sqrt2}\sqrt[4]{\l}\,(x-a)}\left[
\left(e^{-i\pi/4}-\frac{\mu}{\sqrt[4]{\l}}\right)
e^{-i\frac{1}{\sqrt2}\sqrt[4]{\l}\,(x-a)}
+\left(e^{i\pi/4}+\frac{\mu}{\sqrt[4]{\l}}\right)
e^{i\frac{1}{\sqrt2}\sqrt[4]{\l}\,(x-a)}\right]\!,
\eeqa
and, for $x\ge a$,
\beqa
\lqn{
\sum_{k\in K} B_k\prod_{l\in K\setminus\{k\}}
\left(1-\frac{i\mu}{\sqrt[4]{\l}}\,\bar{\t}_l\right)
e^{\t_k\!\sqrt[4]{\l}\,(x-a)}
}\\
&=&
B_0\left(1-\frac{i\mu}{\sqrt[4]{\l}}\,\bar{\t}_1\right)
e^{\t_0\!\sqrt[4]{\l}\,(x-a)}
+B_1\left(1-\frac{i\mu}{\sqrt[4]{\l}}\,\bar{\t}_0\right)
e^{\t_1\!\sqrt[4]{\l}\,(x-a)}
\\
&=&
\frac{1}{\sqrt2}\,e^{-\frac{1}{\sqrt2}\sqrt[4]{\l}\,(x-a)}\left[
\left(e^{-i\pi/4}+\frac{\mu}{\sqrt[4]{\l}}\right)
e^{i\frac{1}{\sqrt2}\sqrt[4]{\l}\,(x-a)}
+\left(e^{i\pi/4}-\frac{\mu}{\sqrt[4]{\l}}\right)
e^{-i\frac{1}{\sqrt2}\sqrt[4]{\l}\,(x-a)}\right]\!.
\eeqa
Therefore, \refp{LFT-tau-Xtau} becomes
\beqa
\Ex\!\!\left[e^{-\l\tau_a^++i\mu X(\tau_a^+)}\right]
&=&
\sqrt2\,e^{i\mu a+\frac{1}{\sqrt2}\sqrt[4]{\l}\,(x-a)}\left[
\cos\bigg(\frac{1}{\sqrt2}\sqrt[4]{\l}\,(x-a)+\frac{\pi}{4}\bigg)\right.
\\
&&
+\left.\frac{i\mu}{\sqrt[4]{\l}} \sin\bigg(\frac{1}{\sqrt2}\sqrt[4]{\l}\,(x-a)\bigg)
\right]\quad\mbox{for $x\le a$,}
\\[2ex]
\Ex\!\!\left[e^{-\l\tau_a^-+i\mu X(\tau_a^-)}\right]
&=&
\sqrt2\,e^{i\mu a-\frac{1}{\sqrt2}\sqrt[4]{\l}\,(x-a)}\left[
\cos\bigg(\frac{1}{\sqrt2}\sqrt[4]{\l}\,(x-a)-\frac{\pi}{4}\bigg)\right.
\\
&&
+\left.\frac{i\mu}{\sqrt[4]{\l}} \sin\bigg(\frac{1}{\sqrt2}\sqrt[4]{\l}\,(x-a)\bigg)
\right]\quad\mbox{for $x\ge a$.}
\eeqa
We retrieve formula~(8.3) of~\cite{nish2}.
}

\subsection{Density functions}

We invert the Laplace-Fourier transform~\refp{LFT-tau-Xtau}.
For this, we proceed in two stages:
we first invert the Fourier transform with respect to
$\mu$ and next invert the Laplace transform with respect to $\l$.

\subsubsection{Inverting with respect to $\mu$}\label{subsubsect-invert-mu-tau-Xtau}

Let us expand the product $\prod_{l\in J\setminus\{j\}}
\left(1-\bar{\t}_l x\right)$ as
\beq\label{expansion-pol}
\prod_{l\in J\setminus\{j\}} \left(1-\bar{\t}_l x\right)
=\sum_{q=0}^{\card J-1} \bar{c}_{jq}(-x)^q
\eeq
where the coefficients $c_{jq}$, $0\le q\le\card J-1$, are the elementary
symmetric functions of the $\t_l$'s, $l\in J\setminus\{j\}$, that is,
more explicitly, $c_{j0}=1$ and for $1\le q\le \card J-1$,
$$
c_{jq}=\s_q\left(\t_l,l\in J\setminus\{j\}\right)
=\sum_{l_1,\ldots,l_q\in J\setminus\{j\}\atop l_1<\cdots<l_q}
\t_{l_1}\cdots\t_{l_q}.
$$
In a similar way, we also introduce $d_{k0}=1$ and for $1\le q\le \card K-1$,
$$
d_{kq}=\s_q\left(\t_l,l\in K\setminus\{k\}\right)
=\sum_{l_1,\ldots,l_q\in K\setminus\{k\}\atop l_1<\cdots<l_q}
\t_{l_1}\cdots\t_{l_q}.
$$
By applying expansion~\refp{expansion-pol} to $x=i\mu/\!\sqrt[N]{\l}$,
we see that~\refp{LFT-tau-Xtau} can be rewritten as
\beqa
\Ex\!\!\left[e^{-\l\tau_a^++i\mu X(\tau_a^+)}\right]
&=&
\sum_{j\in J} A_j
\sum_{q=0}^{\card J-1} \bar{c}_{jq}\left(-\frac{i\mu}{\sqrt[N]{\l}}\right)^q
e^{\t_j\!\!\sqrt[N]{\l}\,(x-a)} \,e^{i\mu a}
\\
&=&
\sum_{q=0}^{\card J-1} \frac{1}{\l^{q/N}}\Bigg[\sum_{j\in J} \bar{c}_{jq}A_j
e^{\t_j\!\!\sqrt[N]{\l}\,(x-a)}\Bigg] (-i\mu)^q\,e^{i\mu a}.
\eeqa
Now, observe that $(-i\mu)^q\,e^{i\mu a}$ is nothing but the
Fourier transform of the $q^{\mbox{\scriptsize th}}$ derivative of
the Dirac distribution viewed as a tempered Schwartz distribution:
\beq\label{Dirac-deriv}
(-i\mu)^q\,e^{i\mu a}=\int_{-\infty}^{+\infty} e^{i\mu z}\,\d_a^{(q)}(z)\,dz.
\eeq
Hence, we have obtained the following intermediate result for the distribution of
$(\tau_a^+,X(\tau_a^+))$ and also for that of $(\tau_a^-,X(\tau_a^-))$.
%
\bpr{
We have, for $\Re(\l)>0$,
\beqan
\Ex\!\!\left[e^{-\l\tau_a^+},X(\tau_a^+)\in dz\right]/dz
&=\;
\dis\sum_{q=0}^{\card J-1} \l^{-q/N}\Bigg[\sum_{j\in J} \bar{c}_{jq}A_j\,
e^{\t_j\!\!\sqrt[N]{\l}\,(x-a)}\Bigg]\d_a^{(q)}(z)
&\mbox{for $x\le a$,}
\nonumber\\[-1ex]
\label{LT-tau-Xtau}\\[-1ex]
\Ex\!\!\left[e^{-\l\tau_a^-},X(\tau_a^-)\in dz\right]/dz
&=\;
\dis\sum_{q=0}^{\card K-1} \l^{-q/N}\Bigg[\sum_{k\in K} \bar{d}_{kq}B_k\,
e^{\t_k\!\!\sqrt[N]{\l}\,(x-a)}\Bigg]\d_a^{(q)}(z)
&\mbox{for $x\ge a$.}
\nonumber\eeqan
}
%
The appearance of the successive derivatives of $\d_a$ suggests to view the
distribution of $(\tau_a^+,X(\tau_a^+))$ as a tempered Schwartz
distribution (that is a Schwartz distribution acting on the space
$\cS$ of the $\cC^{\infty}$-functions exponentially decreasing together with
their derivatives characterized by
$$
\forall \f,\psi\in\cS,\;
\int\!\!\!\int \f(t)\psi(z)\,\P_x\{\tau_a^+\in dt,X(\tau_a^+)\in dz\}
=\Ex\!\!\left[\f(\tau_a^+)\psi(X(\tau_a^+))\right]\!.
$$

\subsubsection{Inversion with respect to $\l$}\label{subsubsect-invert-lambda-tau-Xtau}

In order to extract the densities of $(\tau_a^+,X(\tau_a^+))$ and
$(\tau_a^-,X(\tau_a^-))$ from~\refp{LT-tau-Xtau}, we search functions
$I_{lq}$, $0\le q\le \max(\card I-1,\card J-1)$, such that, for $\Re(\t_l\xi)<0$,
\beq\label{LTI}
\int_0^{+\infty} e^{-\l t} I_{lq}(t;\xi)\,dt =\l^{-q/N} e^{\t_l\!\!\sqrt[N]{\l}\,\xi}.
\eeq
The rhs of~\refp{LTI} seems closed to the Laplace transform of the probability
density function of a completely asymmetric stable random variable, at least for $q=0$.
Nevertheless, because of the presence of the complex term $\t_l$ within
the rhs of~\refp{LTI}, we did not find any precise relationship between
the function $I_{lq}$ and stable processes.
So, we derive below an integral representation for $I_{lq}$.

Invoking Bromwich formula, the function $I_{lq}$ writes
\beqa
I_{lq}(t;\xi)
&=&
\frac{1}{2i\pi} \int_{-i\infty}^{i\infty}
\l^{-q/N} e^{t\l+\t_l\xi\!\sqrt[N]{\l}} \,d\l
\;\;=\;\;
\frac{1}{2\pi} \int_{-\infty}^{+\infty}
(i\l)^{-\frac{q}{N}} e^{it\l+\t_l\xi\!\sqrt[N]{i\l}} \,d\l
\\
&=&
\frac{1}{2\pi}\left[e^{-\frac{i\pi q}{2N}}  \int_0^{+\infty}
\!\l^{-\frac qN} e^{it\l+\t_l e^{\frac{i\pi}{2N}}\xi\!\sqrt[N]{\l}} \,d\l
+e^{\frac{i\pi q}{2N}} \int_0^{+\infty} \!\l^{-\frac qN}
e^{-it\l+\t_l e^{-\frac{i\pi}{2N}}\xi\!\sqrt[N]{\l}} \,d\l\right]\!.
\eeqa
The substitution $\l\longmapsto \l^N$ yields
\beqa
I_{lq}(t;\xi)
&=&
\frac{N}{2\pi}\left[e^{-\frac{i\pi q}{2N}}  \int_0^{+\infty}
\l^{N-q-1} e^{it\l^N+\t_l e^{\frac{i\pi}{2N}}\xi\l} \,d\l \right.
\\
&&
+\left.e^{\frac{i\pi q}{2N}} \int_0^{+\infty} \l^{N-q-1}
e^{-it\l^N+\t_l e^{-\frac{i\pi}{2N}}\xi\l} \,d\l\right]
\eeqa
and the substitutions $\l\longmapsto e^{\pm \frac{i\pi}{2N}}\l$
together with the residues theorem provide
\beqa
I_{lq}(t;\xi)
&=&
\frac{Ni}{2\pi}\left[e^{-\frac{i\pi q}{N}} \int_0^{+\infty}
\l^{N-q-1} e^{-t\l^N+\t_l e^{\frac{i\pi}{N}}\xi\l} \,d\l \right.
\\
&&
-\left.e^{\frac{i\pi q}{N}} \int_0^{+\infty} \l^{N-q-1}
e^{-t\l^N+\t_l e^{-\frac{i\pi}{N}}\xi\l} \,d\l\right]\!.
\eeqa
In particular, for $q=0$ we have, by integration by parts,
\beqan
I_{l0}(t;\xi)
&=&
\frac{i\t_l\xi}{2\pi t}\left[e^{\frac{i\pi}{N}} \int_0^{+\infty}
e^{-t\l^N+\t_l e^{\frac{i\pi}{N}}\xi\l} \,d\l
-e^{-\frac{i\pi}{N}} \int_0^{+\infty}
e^{-t\l^N+\t_l e^{-\frac{i\pi}{N}}\xi\l} \,d\l\right]\!.
\label{Ik0}
\eeqan
%
\brem{
The following relation holds between all the functions $I_{lq}$'s:
$$
\frac{\partial^m I_{lq}}{\partial \xi^m}(t;\xi)= \t_l^m I_{l\,q-m}(t;\xi)
\quad\mbox{for $0\le m\le q$}.
$$
}
%
Hence, \refp{LT-tau-Xtau} can be rewritten as an explicit Laplace transform
with respect to $\l$:
$$
\Ex\!\!\left[e^{-\l\tau_a^+},X(\tau_a^+)\in dz\right]/dz
=\int_0^{+\infty} e^{-\l t}\,dt \left[\sum_{q=0}^{\card J-1}
\Bigg(\sum_{j\in J} \bar{c}_{jq}A_j\,I_{jq}(t;x-a)\Bigg)\d_a^{(q)}(z)\right]\!.
$$
We are able to state the main result of this part.
\bth{
The joint ``distributional densities'' of the vectors $(\tau_a^+,X(\tau_a^+))$ and
\linebreak
$(\tau_a^-,X(\tau_a^-))$ are given by
\beqan
\P_x\{\tau_a^+\in dt,X(\tau_a^+)\in dz\}/dt\,dz
&=\;
\dis\sum_{q=0}^{\card J-1} \cJ_q(t;x-a)\,\d_a^{(q)}(z)
&\mbox{for $x\le a$,}
\nonumber\\[-1.5ex]
\label{density-tau_Xtau}\\[-1.5ex]
\P_x\{\tau_a^-\in dt,X(\tau_a^-)\in dz\}/dt\,dz
&=\;
\dis\sum_{q=0}^{\card K-1} \cK_q(t;x-a)\,\d_a^{(q)}(z)
&\mbox{for $x\ge a$,}
\nonumber
\eeqan
where
$$
\cJ_q(t;\xi)=\sum_{j\in J} \bar{c}_{jq}A_j\,I_{jq}(t;\xi)
\quad\mbox{and}\quad\cK_q(t;\xi)=\sum_{k\in K} \bar{d}_{kq}B_k\,I_{kq}(t;\xi).
$$
}
\brem{
Another expression for $\cJ_q(t;\xi)$, for instance, may be written. Indeed,
for $\xi\le 0$ and $0\le q\le\card J-1$,
\beqan
\cJ_q(t;\xi)
&=&
\frac{Ni}{2\pi} \Bigg[e^{-\frac{i\pi q}{N}} \int_0^{+\infty}
\Bigg(\sum_{j\in J} \bar{c}_{jq}A_j \,e^{\t_j e^{\frac{i\pi}{N}}\xi\l}\Bigg)
\l^{N-q-1} e^{-t\l^N} \,d\l
\nonumber\\
&&
-\,e^{\frac{i\pi q}{N}} \int_0^{+\infty}
\Bigg(\sum_{j\in J} \bar{c}_{jq}A_j \,e^{\t_j e^{-\frac{i\pi}{N}}\xi\l}\Bigg)
\l^{N-q-1} e^{-t\l^N} \,d\l\Bigg].
\label{integral-cJ}
\eeqan
The second integral displayed in~\refp{integral-cJ} is the conjugate
of the first one. In effect,
by introducing the symmetry $\s:j\in J\longmapsto \s(j)\in J$
such that $\t_{\s(j)}=\bar{\t}_j$, we can see that
$$
A_{\s(j)}=\prod_{l\in J\setminus\{\s(j)\}} \frac{\t_l}{\t_l-\t_{\s(j)}}
=\prod_{l\in J\setminus\{j\}} \frac{\t_{\s(l)}}{\t_{\s(l)}-\t_{\s(j)}}
=\prod_{l\in J\setminus\{j\}}
\frac{\bar{\t}_l}{\vphantom{\bar{\bar{\t}}}\bar{\t}_l-\bar{\t}_j}
=\bar{A}_j
$$
and
$$
c_{\s(j)q}=\s_q\left(\t_l,l\in J\setminus\{\s(j)\}\right)
=\s_q\left(\t_{\s(l)},l\in J\setminus\{j\}\right)
=\s_q\left(\bar{\t}_l,l\in J\setminus\{j\}\right)=\bar{c}_{jq}.
$$
So, the sum lying within the second integral in~\refp{integral-cJ} writes
$$
\sum_{j\in J} \bar{c}_{\s(j)q}A_{\s(j)} \,e^{\t_{\s(j)} e^{-\frac{i\pi}{N}}\xi\l}
=\sum_{j\in J} c_{jq}\bar{A}_j \,e^{\bar{\t}_j e^{-\frac{i\pi}{N}}\xi\l}
=\Bigg(\sum_{j\in J} \bar{c}_{jq}A_j \,e^{\t_j e^{\frac{i\pi}{N}}\xi\l}\Bigg)^-.
$$
As a result,
$$
\cJ_q(t;\xi)=-\frac{N}{\pi} \,\Im\Bigg[e^{-\frac{i\pi q}{N}} \int_0^{+\infty}
\Bigg(\sum_{j\in J} \bar{c}_{jq}A_j \,e^{\t_j e^{\frac{i\pi}{N}}\xi\l}\Bigg)
\l^{N-q-1} e^{-t\l^N} \,d\l \Bigg].
$$
In particular, $\cJ_q(t;\xi)$ is real and for $q=0$ we have, since $c_{j0}=1$
and $\sum_{j\in J} A_j=1$,
\beqa
\cJ_0(t;\xi)
&=&
-\frac{N}{\pi} \,\Im\Bigg[\int_0^{+\infty}
\Bigg(\sum_{j\in J} A_j \,e^{\t_j e^{\frac{i\pi}{N}}\xi\l}\Bigg)
\l^{N-1} e^{-t\l^N} \,d\l \Bigg]
\\
&=&
-\frac{\xi}{\pi t} \,\Im\Bigg[e^{\frac{i\pi}{N}}\int_0^{+\infty}
\Bigg(\sum_{j\in J} \t_j A_j \,e^{\t_j e^{\frac{i\pi}{N}}\xi\l}\Bigg)
e^{-t\l^N} \,d\l \Bigg]
\eeqa
which is nothing but $\P_x\{\tau_a^+\in dt\}/dt$.
}

\subsection{Distribution of the hitting places}\label{subsect-Xtau}

We now derive the distribution of the hitting places $X(\tau_a^+)$
and $X(\tau_a^-)$. To do this for $X(\tau_a^+)$ for example, we
integrate~\refp{density-tau_Xtau} with respect to $t$:
\beqan
\lqn{\P_x\{X(\tau_a^+)\in dz\}/dz
=\int_0^{+\infty} \P_x\{\tau_a^+\in dt,X(\tau_a^+)\in dz\}/dz
}\nonumber\\
&=&
\sum_{q=0}^{\card J-1} \left[\int_0^{+\infty} \cJ_q(t;x-a)\,dt\right]\d_a^{(q)}(z)
\nonumber\\
&=&
-\frac{N}{\pi} \sum_{q=0}^{\card J-1} |x-a|^q \Bigg[\int_0^{+\infty}
\Im\Bigg(e^{-\frac{i\pi q}{N}}
\sum_{j\in J} \bar{c}_{jq}A_j \,e^{-\t_j e^{\frac{i\pi}{N}}\l}\Bigg)
\frac{d\l}{\l^{q+1}} \Bigg]\,\d_a^{(q)}(z).
\label{Xtau+-intermediate}
\eeqan
We need two lemmas for carrying out the integral lying in~\refp{Xtau+-intermediate}.
%
\blm{\label{lemma-integral}
For any integers $m,n$ such that $1\le n\le m-1$ and any complexes $a_1,\ldots,a_m$
and $b_1,\ldots,b_m$ such that $\Re(b_j)\ge 0$
and \mbox{$\Im\Big(\sum_{j=1}^m a_j b_j^l\Big)=0$} for $0\le l\le n-1$,
$$
\int_0^{+\infty} \Im\Bigg(\sum_{j=1}^m a_j e^{-b_j \l}\Bigg)\frac{d\l}{\l^n}
=\frac{(-1)^n}{(n-1)!}\, \Im\Bigg(\sum_{j=1}^m a_j b_j^{n-1}\log b_j\Bigg).
$$
}
%
\dem
We proceed by induction on $n$.

For $n=1$, because of the condition $\Im\Big(\sum_{j=1}^m a_j\Big)=0$,
we can replace $\Im(a_m)$ by $-\Im\Big(\sum_{j=1}^{m-1} a_j\Big)$. This gives
$$
\int_0^{+\infty} \Im\Bigg(\sum_{j=1}^m a_j e^{-b_j \l}\Bigg)\frac{d\l}{\l}
=\int_0^{+\infty} \Im\Bigg[\sum_{j=1}^{m-1} a_j
\left(e^{-b_j \l}-e^{-b_m \l}\right)\Bigg]\frac{d\l}{\l}.
$$
The foregoing integral involves the elementary integral below:
$$
\int_0^{+\infty} \Im\left(e^{-b_j \l}-e^{-b_m \l}\right)
\frac{d\l}{\l}=\Im\left(\log b_m-\log b_j\right).
$$
Therefore,
$$
\int_0^{+\infty} \Im\Bigg(\sum_{j=1}^m a_j e^{-b_j \l}\Bigg)\frac{d\l}{\l}
=\Im\Bigg[\sum_{j=1}^{m-1} a_j(\log b_m-\log b_j)\Bigg]
=-\Im\Bigg[\sum_{j=1}^m a_j\log b_j\Bigg]
$$
which proves Lemma~\ref{lemma-integral} in the case $n=1$.

Assume now the result of the lemma valid for an integer $n\ge 1$.
Let $m$ be an integer such that $m\ge n+2$ and
$a_1,\ldots,a_m$ and $b_1,\ldots,b_m$ be complex numbers such that $\Re(b_j)\ge 0$
and $\Im\Big(\sum_{j=1}^m a_j b_j^l\Big)=0$ for $0\le l\le n$.
By integration by parts, we have
$$
\int_0^{+\infty} \Im\Bigg(\sum_{j=1}^m a_j e^{-b_j \l}\Bigg)\frac{d\l}{\l^{n+1}}
=\left[-\frac{1}{n\l^n} \Im\Bigg(\sum_{j=1}^m a_j e^{-b_j \l}\Bigg)\!\right]_0^{\infty}
\!\!-\frac 1n \int_0^{+\infty} \Im\Bigg(\sum_{j=1}^m a_j b_j e^{-b_j \l}\Bigg)
\frac{d\l}{\l^n}.
$$
Applying L'H\^opital's rule $n$ times, we see, using the condition
$\Im\Big(\sum_{j=1}^m a_j b_j^l\Big)=0$ for $0\le l\le n$, that
$\left[-\frac{1}{n\l^n}
\Im\left(\sum_{j=1}^m a_j e^{-b_j \l}\right)\right]_0^{\infty}=0$.
Putting $\tilde{a}_j=a_jb_j$, we get
$$
\int_0^{+\infty} \Im\Bigg(\sum_{j=1}^m a_j e^{-b_j \l}\Bigg)\frac{d\l}{\l^{n+1}}
=-\frac 1n \int_0^{+\infty} \Im\Bigg(\sum_{j=1}^m \tilde{a}_j e^{-b_j \l}\Bigg)
\frac{d\l}{\l^n}.
$$
We have
$\Im\Big(\sum_{j=1}^m \tilde{a}_j b_j^l\Big)
=\Im\Big(\sum_{j=1}^m a_j b_j^{l+1}\Big)=0$
for $0\le l\le n-1$. Then, invoking the recurrence hypothesis,
the intermediate integral writes
$$
\int_0^{+\infty} \Im\Bigg(\sum_{j=1}^m \tilde{a}_j e^{-b_j \l}\Bigg)\frac{d\l}{\l^n}
=\frac{(-1)^n}{(n-1)!}\, \Im\Bigg(\sum_{j=1}^m \tilde{a}_j b_j^{n-1}\log b_j\Bigg)
$$
and thus
$$
\int_0^{+\infty} \Im\Bigg(\sum_{j=1}^m a_j e^{-b_j \l}\Bigg)\frac{d\l}{\l^{n+1}}
=\frac{(-1)^{n+1}}{n!}\, \Im\Bigg(\sum_{j=1}^m a_j b_j^n\log b_j\Bigg)
$$
which achieve the proof of Lemma~\ref{lemma-integral}.
\fin

%
\blm{\label{lemma-sum}
We have, for $0\le p\le q\le \card J-1$,
$$
\sum_{j\in J} \bar{c}_{jq}\t_j^pA_j=
\left\{\begin{array}{ll}
0&\mbox{if $p\le q-1$,}
\\
(-1)^q &\mbox{if $p=q$.}
\end{array}\right.
$$
}
%
\dem
Consider the following polynomial:
\beqa
\sum_{q=0}^{\card J-1} \Bigg(\sum_{j\in J} \bar{c}_{jq}\t_j^pA_j\Bigg)(-x)^q
&=&
\sum_{j\in J} \t_j^pA_j \sum_{q=0}^{\card J-1} \bar{c}_{jq}(-x)^q
\\
&=&
\sum_{j\in J} \t_j^pA_j \prod_{l\in J\setminus\{j\}} (1-\bar{\t}_lx)
\\
&=&
\prod_{l\in J} (1-\bar{\t}_lx) \sum_{j\in J} \frac{\t_j^pA_j}{1-\bar{\t}_jx}.
\eeqa
We then obtain, due to~\refp{expansion}, if $p\le\card J-1$,
$$
\sum_{q=0}^{\card J-1} \Bigg(\sum_{j\in J} \bar{c}_{jq}\t_j^pA_j\Bigg)(-x)^q=x^p
$$
which entails the result by identifying the coefficients of both polynomials above.
\fin

Now, we state the following remarkable result.
%
\bth{
The ``distributional densities'' of $X(\tau_a^+)$ and $X(\tau_a^-)$
are given by
\beqan
\P_x\{X(\tau_a^+)\in dz\}/dz
&=\;
\dis\sum_{q=0}^{\card J-1} (-1)^q\frac{(x-a)^q}{q!}\,\d_a^{(q)}(z)
&\mbox{for $x\le a$,}
\nonumber\\[-1ex]
\label{density-Xtau}\\[-1ex]
\P_x\{X(\tau_a^-)\in dz\}/dz
&=\;
\dis\sum_{q=0}^{\card K-1} (-1)^q\frac{(x-a)^q}{q!}\,\d_a^{(q)}(z)
&\mbox{for $x\ge a$.}
\nonumber
\eeqan
}
%
It is worth that the distributions of $X(\tau_a^+)$ and $X(\tau_a^-)$
are linear combinations of the successive derivatives of the Dirac distribution
$\d_a$. This noteworthy fact has already been observed by Nishioka~\cite{nish1,nish2}
in the case $N=4$ and the author spoke of ``monopoles'' and ``dipoles''
respectively related to $\d_a$ and $\d'_a$ (see also~\cite{nish3} for more
account about relationships between monopoles/dipoles and different kinds of
absorbed/killed pseudo-processes).
More generally,~\refp{density-Xtau} suggests to speak of ``multipoles''
related to the $\d_a^{(q)}$'s.

In the case of Brownian motion ($N=2$), the trajectories are continuous,
so $X(\tau_a^{\pm})=a$ and then we classically write
$\P_x\{X(\tau_a^{\pm})\in dz\}=\d_a(dz)$
where $\d_a$ is viewed as the Dirac probability measure.
For $N\ge 4$, it emerges from~\refp{density-Xtau} that the distributional
densities of $X(\tau_a^{\pm})$ are concentrated at the point $a$ through a
sequence of successive derivatives of $\d_a$ where $\d_a$ is now viewed as a
Schwartz distribution. Hence, we could guess in~\refp{density-Xtau} a
curious and unclear kind of continuity.
In Subsection~\ref{before}, we study the distribution of $X(\tau_a^{\pm}-)$
which will reveal itself to coincide with that of $X(\tau_a^{\pm})$.
This will confirm this idea of continuity.

\dem
Let us evaluate the integral lying in~\refp{Xtau+-intermediate}.
We have, thanks to Lemma~\ref{lemma-sum},
$$
e^{-\frac{i\pi q}{N}} \sum_{j\in J} \bar{c}_{jq}A_j
\left(\t_j e^{\frac{i\pi}{N}}\right)^l
=e^{\frac{i\pi}{N}\,(l-q)} \sum_{j\in J} \bar{c}_{jq}A_j \t_j^l
=0 \mbox{ if $l\le q-1$.}
$$
Therefore, the conditions of Lemma~\ref{lemma-integral} are fulfilled and we get
\beqa
\lqn{\P_x\{X(\tau_a^+)\in dz\}/dz}
\\
&=&
\sum_{q=0}^{\card J-1} \frac{(-1)^qN}{\pi q!}\, |x-a|^q
\Bigg[\sum_{j\in J} \Im\left(\bar{c}_{jq}A_j\t_j^q
\log\Big(\t_je^{\frac{i\pi}{N}}\Big)\right)\Bigg]\,\d_a^{(q)}(z)
\\
&=&
\sum_{q=0}^{\card J-1} \frac{(-1)^qN}{\pi q!}\, |x-a|^q
\Bigg[\Re\Bigg(\sum_{j\in J} \bar{c}_{jq}A_j\t_j^q\arg(\t_j)\Bigg)
+\frac{\pi}{N}\,\Re\Bigg(\sum_{j\in J} \bar{c}_{jq}A_j\t_j^q \Bigg)\Bigg]
\,\d_a^{(q)}(z).
\eeqa
The second sum lying within the brackets is equal, by Lemma~\ref{lemma-sum},
to $(-1)^q$. The first one vanishes: indeed, by using the symmetry
$\s:j\in J\longmapsto \s(j)\in J$ such that $\t_{\s(j)}=\bar{\t}_j$,
\beqa
\Re\Bigg(\sum_{j\in J} \left(\bar{c}_{jq}A_j\t_j^q\right) \arg(\t_j)\Bigg)
&=&
\frac12\Bigg(\sum_{j\in J}\bar{c}_{jq}A_j\t_j^q\arg(\t_j)
+\sum_{j\in J}c_{jq}\bar{A}_j\bar{\t}_j^q\arg(\t_j)\Bigg)
\\
&=&
\frac12\Bigg(\sum_{j\in J}\bar{c}_{jq}A_j\t_j^q\arg(\t_j)
+\sum_{j\in J}c_{\s(j)q}\bar{A}_{\s(j)}\bar{\t}_{\s(j)}^q\arg(\t_{\s(j)})\Bigg).
\eeqa
The terms of the second last sum are the opposite of those of the first sum
since
$$
c_{\s(j)q}\bar{A}_{\s(j)}\bar{\t}_{\s(j)}^q=\bar{c}_{jq}A_j\t_j^q
\qquad\mbox{and}\qquad \arg(\t_{\s(j)})=-\arg(\t_j)
$$
which proves the assertion. As a result, we get~\refp{density-Xtau}.
\fin

\subsection{Fourier transforms of the hitting places}\label{subsect-FT-Xtau}

By using~\refp{density-Xtau} and~\refp{Dirac-deriv}, it is easy to derive the
Fourier transforms of the hitting places $X(\tau_a^+)$ and $X(\tau_a^-)$.
%
\bpr{
The Fourier transforms of $X(\tau_a^+)$ and $X(\tau_a^-)$ are given by
\beqan
\Ex\!\!\left[e^{i\mu X(\tau_a^+)}\right]
&=\;
\dis e^{i\mu a}\sum_{q=0}^{\card J-1} \frac{(x-a)^q}{q!}\,(i\mu)^q
&\mbox{for $x\le a$,}
\nonumber\\[-1ex]
\label{FT-Xtau}\\[-1ex]
\Ex\!\!\left[e^{i\mu X(\tau_a^-)}\right]
&=\;
\dis e^{i\mu a}\sum_{q=0}^{\card K-1} \frac{(x-a)^q}{q!}\,(i\mu)^q
&\mbox{for $x\ge a$.}
\nonumber\eeqan
}
%
In this part, we suggest to retrieve~\refp{FT-Xtau} by letting $\l$ tend to $0^+$
in~\refp{LFT-tau-Xtau}. We rewrite~\refp{LFT-tau-Xtau}, for instance for $x\le a$, as
\beqan
\Ex\!\!\left[e^{-\l\tau_a^++i\mu X(\tau_a^+)}\right]
&=&
e^{i\mu a}\prod_{l\in J} \left(1-\frac{i\mu}{\sqrt[N]{\l}}\,\bar{\t}_l\right)
\sum_{j\in J} \frac{A_j}{1-\frac{i\mu}{\sqrt[N]{\l}}\,\bar{\t}_j}
\,e^{\t_j\!\!\sqrt[N]{\l}\,(x-a)}
\nonumber\\
&=&
\left(-\frac{i\mu}{\sqrt[N]{\l}}\right)^{\card J-1}
\Bigg(\prod_{j\in J} \bar{\t}_j\Bigg) e^{i\mu a}
\nonumber\\
&&
\times\prod_{j\in J}
\left(1-\frac{\t_j\!\!\sqrt[N]{\l}}{i\mu}\right)
\sum_{j\in J} \frac{\t_jA_j}{1-\frac{\t_j\!\!\sqrt[N]{\l}}{i\mu}}
\,e^{\t_j\!\!\sqrt[N]{\l}\,(x-a)}.
\label{asympt-prelim}
\eeqan
Using the elementary expansions, as $\l\to 0^+$,
\beqa
\frac{1}{1-\frac{\t_j\!\!\sqrt[N]{\l}}{i\mu}}
&=&
\sum_{p=0}^{\card J-1} \left(\frac{\t_j\!\!\sqrt[N]{\l}}{i\mu}\right)^p
+o\left(\l^{(\card J-1)/N}\right)\!,
\\
e^{\t_j\!\!\sqrt[N]{\l}\,(x-a)}
&=&
\sum_{q=0}^{\card J-1} \frac{1}{q!}\,\left(\t_j\!\!\sqrt[N]{\l}\,(x-a)\right)^q
+o\left(\l^{(\card J-1)/N}\right)\!,
\eeqa
yields
\beqa
\sum_{j\in J} \frac{\t_jA_j}{1-\frac{\t_j\!\!\sqrt[N]{\l}}{i\mu}}
\,e^{\t_j\!\!\sqrt[N]{\l}\,(x-a)}
&=&
\sum_{j\in J} \t_jA_j \left[\sum_{r=0}^{\card J-1} \left(\sum_{q=0}^r
\frac{(x-a)^q}{q!(i\mu)^{r-q}}\right) \left(\t_j\!\!\sqrt[N]{\l}\right)^r
\right]\!+o\left(\l^{(\card J-1)/N}\right)
\\
&=&
\sum_{r=0}^{\card J-1}  \left(\sum_{j\in J} \t_j^{r+1} A_j\right)
\left(\sum_{q=0}^r \frac{(x-a)^q}{q!(i\mu)^{r-q}}\right) \l^{r/N}
+o\left(\l^{(\card J-1)/N}\right)\!.
\eeqa
On the other hand, applying~\refp{expansion} to $x=0$ gives
$$
\sum_{j\in J} \t_j^{r+1} A_j
=\left\{\begin{array}{ll}
0&\mbox{if $r\le\card J-2$,}
\\
(-1)^{\card J-1}\prod_{j\in J}\t_j
&\mbox{if $r=\card J-1$.}
\end{array}\right.
$$
Therefore,
\beq
\sum_{j\in J} \frac{\t_jA_j}{1-\frac{\t_j\!\!\sqrt[N]{\l}}{i\mu}}
\,e^{\t_j\!\!\sqrt[N]{\l}\,(x-a)}
\underset{\l \to 0^+}{\overset{}{\sim}}
(-1)^{\card J-1}\Bigg(\prod_{j\in J}\t_j\Bigg)
\Bigg(\sum_{q=0}^{\card J-1} \frac{(x-a)^q}{q!(i\mu)^{\card J-1-q}}\Bigg)
\l^{(\card J-1)/N}.
\label{asympt-inter}
\eeq
Consequently, the limit of
$\Ex\!\!\left[e^{-\l\tau_a^++i\mu X(\tau_a^+)}\right]$ as $\l\to 0^+$ ensues.
The constant arising when combining~\refp{asympt-prelim} and~\refp{asympt-inter} is
$$
(-1)^{\card J-1}\Bigg(\prod_{j\in J}\t_j\Bigg)
\times\Bigg(\prod_{j\in J}\bar{\t}_j\Bigg)
(-i\mu)^{\card J-1} e^{i\mu a}=(i\mu)^{\card J-1} e^{i\mu a}.
$$
In view of~\refp{FT-Xtau}, we have proved the equality
$$
\lim_{\l\to 0^+}\Ex\!\!\left[e^{-\l\tau_a^++i\mu X(\tau_a^+)}\right]
=\Ex\!\!\left[e^{i\mu X(\tau_a^+)}\right].
$$
%
\brem{
The distribution of $X(\tau_a^+)$ may also be deduced from the joint distribution
of $\left(\tau_a^+,X(\tau_a^+)\right)$ through~\refp{LT-tau-Xtau}.
Indeed, by letting $\l$ tend to $0$ in~\refp{LT-tau-Xtau} and using
elementary expansions together with Lemma~\ref{lemma-sum},
\beqa
\sum_{j\in J} \bar{c}_{jq}A_j\, e^{\t_j\!\!\sqrt[N]{\l}\,(x-a)}
&=&
\sum_{j\in J} \bar{c}_{jq}A_j \sum_{p=0}^q \frac{(\t_j(x-a))^p}{p!}
\,\l^{p/N} +o\left(\l^{q/N}\right)
\\
&=&
\sum_{p=0}^q \Bigg(\sum_{j\in J} \bar{c}_{jq}A_j\t_j^p \Bigg)
\frac{(x-a)^p}{p!} \,\l^{p/N}+o\left(\l^{q/N}\right)
\\
&\underset{\l \to 0^+}{\overset{}{\sim}}&
(-1)^q \frac{(x-a)^q}{q!} \,\l^{q/N}\!,
\eeqa
which, with~\refp{LT-tau-Xtau}, confirms~\refp{density-Xtau}.
}


\subsection{Strong Markov property for $\tau_a^{\pm}$}

We roughly state a strong Markov property related to the hitting times
$\tau_a^{\pm}$.
%
\bth{
For suitable functionals $F$ and $G$, we have
\beqan
\lqn{
\E_x\!\!\left[F\!\left((X(t))_{0\le t<\tau_a^{\pm}}\right)
G\!\left((X(t+\tau_a^{\pm}))_{t\ge 0}\right)\right]
}\nonumber\\[-1ex]
&=&
\E_x\!\!\left[F\!\left((X(t))_{0\le t<\tau_a^{\pm}}\right)
\E_{X(\tau_a^{\pm})}\left[G\!\left((X(t))_{t\ge 0}\right)\right]\right]\!,
\label{strong-markov1}\\
\E_x\!\!\left[G\!\left((X(t+\tau_a^+))_{t\ge 0}\right)\right]
&=&
\sum_{q=0}^{\card J-1} \frac{(x-a)^q}{q!}\left.\frac{\partial^q}{\partial z^q}\,
\E_z\!\!\left[G\!\left((X(t))_{t\ge 0}\right)\right]\right|_{z=a} \mbox{ if $x\ge a$},
\nonumber\\[-1ex]
\label{strong-markov2}
\\[-1ex]
\E_x\!\!\left[G\!\left((X(t+\tau_a^-))_{t\ge 0}\right)\right]
&=&
\sum_{q=0}^{\card K-1} \frac{(x-a)^q}{q!}\left.\frac{\partial^q}{\partial z^q}\,
\E_z\!\!\left[G\!\left((X(t))_{t\ge 0}\right)\right]\right|_{z=a} \mbox{ if $x\le a$}.
\nonumber
\eeqan
}
%
\dem
We first consider the step-process $X_n$ and we use the notations of
Subsection~\ref{subsect-LFT-tau-Xtau}. On the set $\{\tauan^+=k/2^n\}$, the
quantities $F\!\left((X_n(t))_{0\le t<\tauan^+}\right)$ and
$G\!\left((X_n(t+\tauan^+))_{t\ge 0}\right)$ depend respectively on
$X_{n,0},X_{n,1},\ldots,X_{n,k-1}$ and $X_{n,k},X_{n,k+1},\ldots$ So we can set,
if $\tauan^+=k/2^n$,
\beqa
F\!\left((X_n(t))_{0\le t<\tauan^+}\right)
&=&
F_k(X_{n,0},X_{n,1},\ldots,X_{n,k-1})\;=\;F_{n,k-1},
\\
G\!\left((X_n(t+\tauan^+))_{t\ge 0}\right)
&=&
G_k(X_{n,k},X_{n,k+1},\ldots)\;=\;G_{n,k}.
\eeqa
Therefore,
$$
F\!\left((X_n(t))_{0\le t<\tauan^+}\right) G\!\left((X_n(t+\tauan^+))_{t\ge 0}\right)
=\sum_{k=1}^{\infty} F_{n,k-1} G_{n,k} \ind_{\{M_{n,k-1}<a\le M_{n,k}\}}.
$$
Taking the expectations, we get for $x\le a$:
\beqan
\lqn{
\E_x\!\!\left[F\!\left((X_n(t))_{0\le t<\tauan^+}\right)G\!\left((X_n(t+\tauan^+))_{t\ge 0}\right)\right]
}&=&
\sum_{k=1}^{\infty} \E_x\!\!\left[F_{n,k-1} \ind_{\{M_{n,k-1}<a\le M_{n,k}\}}
\E_{X_{n,k}}(G_{n,0})\right]
\nonumber\\
&=&
\E_x\!\!\left[F\!\left((X_n(t))_{0\le t<\tauan^+}\right)
\E_{X(\tauan^+)}\left[G\!\left((X_n(t))_{t\ge 0}\right)\right]\right]
\label{strong-markov-step}
\eeqan
and~\refp{strong-markov1} ensues by taking the limit
of~\refp{strong-markov-step} as $n$ tends to $+\infty$ in the sense of Definition~\ref{def3}.

In particular, choosing $F=1$, \refp{strong-markov1} writes for $x\le a$
$$
\E_x\!\!\left[G\!\left((X(t+\tau_a^+))_{t\ge 0}\right)\right]
=\int_{-\infty}^{+\infty} \P_x\{X(\tau_a^+)\in dz\}
\E_z\!\!\left[G\!\left((X(t))_{t\ge 0}\right)\right]
$$
which, by~\refp{density-Xtau}, immediately yields~\refp{strong-markov2}.
\fin

The argument in favor of discontinuity evoked in~\cite{arcsine} should fail
since, in view of~\refp{LTI}, a term is missing when applying the strong Markov property.


\subsection{Just before the hitting time}\label{before}

In order to lighten the notations, we simply write $\tau_a^{\pm}=\tau_a$ and
we introduce the jump $\D_aX=X(\tau_a)-X(\tau_a-)$.
%
\bpr{\label{prop-jump}
The Laplace-Fourier transform of the vector $(\tau_a,X(\tau_a-),\D_a X)$
is related to those of the vectors $(\tau_a,X(\tau_a-))$ and
$(\tau_a,X(\tau_a))$ according as, for $\Re(\l)>0$ and $\mu,\nu\in\R$,
\beq\label{LT-jump}
\Ex\!\!\left[e^{-\l\tau_a+i\mu X(\tau_a-)+i\nu \D_aX}\right]
=\Ex\!\!\left[e^{-\l\tau_a+i\mu X(\tau_a-)}\right]
=\Ex\!\!\left[e^{-\l\tau_a+i\mu X(\tau_a)}\right]\!.
\eeq
}
%
\dem
The proof of Proposition~\ref{prop-jump} is similar to that of
Lemma~\ref{lemma-relation}. So, we outline the main steps with less details.
We consider only the case where $\tau_a=\tau_a^+$ and $x\le a$, the other one
is quite similar.

\noindent\textbf{$\bullet$ Step 1}

Recall that for the step-process $(X_n(t))_{t\ge 0}$, the first hitting time
$\tauan^+$ is the instant $\tnk$ with $k$ such that $M_{n,k-1}\le a$ and
$M_{n,k}> a$, and then $X(\tauan-)=X_{n,k-1}$ and $X(\tauan)=X_{n,k}$.
Set $\D_{n,k}=X_{n,k}-X_{n,k-1}$. We have, for $x\le a$,
\beqan
\lqn{
e^{-\l\tauan+i\mu X_n(\tauan-)+i\nu\D_aX_n}
}\nonumber\\[-3ex]
&=&
\sum_{k=1}^{\infty} e^{-\l \tnk+i\mu X_{n,k-1}+i\nu\D_{n,k}}
\ind_{\{M_{n,k-1}\le a< M_{n,k}\}}
\nonumber\\
&=&
e^{-\l \tnun+i\mu X_{n,0}+i\nu(X_{n,1}-X_{n,0})}
\nonumber\\
&&
+\sum_{k=1}^{\infty} \left[e^{-\l \tnkun+i\mu X_{n,k}+i\nu \D_{n,k+1}}
-e^{-\l \tnk+i\mu X_{n,k-1}+i\nu \D_{n,k}}\right]\ind_{\{M_{n,k}\le a\}}.
\label{exp_tauan_delta}
\eeqan

\noindent\textbf{$\bullet$ Step 2}

We take the expectation of~\refp{exp_tauan_delta}:
\beqa
\lqn{
\Ex\!\!\left[e^{-\l\tauan+i\mu X_n(\tauan-)+i\nu\D_aX_n}\right]
=e^{-\l/2^n+i\mu x+\k_N(i\nu)^N/2^n}
}\\[-2ex]
&&
+\sum_{k=1}^{\infty} e^{-\l \tnk}
\Ex\!\!\left[e^{i\mu X_{n,k-1}}\ind_{\{M_{n,k-1}\le a\}}
\left(e^{-\l/2^n+i\mu \D_{n,k}+i\nu \D_{n,k+1}} -e^{i\nu \D_{n,k}}\right)
\ind_{\{X_{n,k}\le a\}}\right]\!.
\eeqa
The expectation lying in the rhs of the foregoing equality can be evaluated
as follows:
\beqa
\lqn{
\Ex\!\!\left[e^{i\mu X_{n,k-1}}\ind_{\{M_{n,k-1}\le a\}}
\left(e^{-\l/2^n+i\mu \D_{n,k}+i\nu \D_{n,k+1}} -e^{i\nu \D_{n,k}}\right)
\ind_{\{X_{n,k}\le a\}}\right]
}\\[-1ex]
&=&
\int_{-\infty}^a e^{i\mu y} \,\P_x\{X_{n,k-1}\in dy,M_{n,k-1}\le a\}
\\
&&
\times\,\E_0\!\!\left[\left(e^{-\l/2^n+i\mu \D_{n,1}+i\nu \D_{n,2}}
-e^{i\nu \D_{n,1}}\right) \ind_{\{\D_{n,1}\le a-y\}}\right]
\\
&=&
\int_{-\infty}^a e^{i\mu y} \P_x\{X_{n,k-1}\in dy,M_{n,k-1}\le a\}
\\
&&
\times\!\left[e^{-\l/2^n} \E_0\!\!\left(e^{i\mu X_{n,1}}
\ind_{\{X_{n,1}\le a-y\}}\right) \E_0\!\!\left(e^{i\nu X_{n,1}}\right)
-\E_0\!\!\left(e^{i\nu X_{n,1}}\ind_{\{X_{n,1}\le a-y\}}\right)\right]\!.
\eeqa
For computing the term within brackets, we need the following quantities:
$$
\E_0\!\!\left(e^{i\mu\mbox{\scriptsize (or $\nu$)}
X_{n,1}}\ind_{\{X_{n,1}\le a-y\}}\right)
=\int_{-\infty}^{a-y} e^{i\mu\mbox{\scriptsize (or $\nu$)}z}\,p(1/2^n;-z)\,dz,
\quad\E_0\!\!\left(e^{i\nu X_{n,1}}\right)=e^{\k_N(i\nu)^N/2^n}.
$$
With these relations at hand, we get
\beqan
\lqn{
\Ex\!\!\left[e^{-\l\tauan+i\mu X_n(\tauan-)+i\nu\D_aX_n}\right]
}\nonumber\\[-2ex]
&=&
e^{-(\l-\k_N(i\nu)^N)/2^n+i\mu x}+\frac{1}{2^n}\sum_{k=1}^{\infty} e^{-\l \tnk}
\int_{-\infty}^a e^{i\mu y} \,\P_x\{X_{n,k-1}\in dy,M_{n,k-1}\le a\}
\nonumber\\
&&
\times 2^n\!\!\left[e^{-(\l-\k_N(i\nu)^N)/2^n}
\!\!\int_{-\infty}^{a-y} \!\! e^{i\mu z}\,p(1/2^n;-z)\,dz
-\!\!\int_{-\infty}^{a-y} \!\! e^{i\nu z}\,p(1/2^n;-z)\,dz\right]\!.
\label{LT-tauan-delta}
\eeqan

\noindent\textbf{$\bullet$ Step 3}

We now take the limit of~\refp{LT-tauan-delta} as $n$ tends to infinity:
\beqa
\lqn{
\Ex\!\!\left[e^{-\l\tau_a+i\mu X(\tau_a-)+i\nu \D_aX}\right]
}\\[-1ex]
&=&
e^{i\mu x}+\int_0^{\infty} e^{-\l t}\,dt
\int_{-\infty}^a e^{i\mu y} \,\f(\l,\mu,\nu;y)\,\P_x\{X(t)\in dy,M(t)\le a\}
\eeqa
where we set, for $y<a$,
\beq\label{function-phi}
\f(\l,\mu,\nu;y)=\lim_{\e\to 0^+} \frac{1}{\e}\left[e^{-(\l-\k_N(i\nu)^N)\e}
\int_{-\infty}^{a-y} e^{i\mu z}\,p(\e;-z)\,dz
-\int_{-\infty}^{a-y} e^{i\nu z}\,p(\e;-z)\,dz\right]\!.
\eeq
\noindent\textbf{$\bullet$ Step 4}

For evaluating the above function $\f$, we need two lemmas.
%
\blm{\label{moments}
For $0\le p\le N$, we have
$$
\E_0[X(t)^p]=\left\{\begin{array}{ll}
1 & \mbox{for $p=0$,}\\
0 & \mbox{for $1\le p\le N-1$,}\\
\k_NN!\,t & \mbox{for $p=N$.}
\end{array}\right.
$$
}
%
\dem
By differentiating $k$ times the identity $\E_0\Big(e^{iuX(t)}\Big)=e^{\k_N(iu)^Nt}$
with respect to $u$ and next substituting $u=0$, we have that
$$
\E_0[X(t)^k]=(-i)^k \left.\frac{\partial^k}{\partial u^k}
\Big[e^{\k_N(iu)^Nt}\Big]\right|_{u=0}.
$$
Fix a complex number $\a\neq0$. It can be easily seen by induction that there
exists a family of polynomials $(P_k)_{k\in\N}$ such that, for all $k\in\N$,
\beq\label{derivatives-exp-auN}
\frac{\partial^k}{\partial u^k} \Big(e^{\a u^N}\Big)=P_k(u)\,e^{\a u^N}.
\eeq
In particular, we have $P_0(u)=1$ and $P_1(u)=N\a u^{N-1}$.
Using the Leibniz rule, we obtain
\beqa
P_k(u)
&=&
e^{-\a u^N}\frac{\partial^k}{\partial u^k} \Big(e^{\a u^N}\Big)
\;=\;
e^{-\a u^N}\frac{\partial^{k-1}}{\partial u^{k-1}}
\Big(N\a \,u^{N-1} e^{\a u^N}\Big)
\\
&=&
N!\,\a \sum_{j=\max(0,k-N)}^{k-1}
\left(\!\!\begin{array}{c} k-1\\ j\end{array}\!\!\right)
\frac{u^{N+j-k}}{(N+j-k)!}\,P_j(u).
\eeqa
This ascertains the aforementioned induction and gives, for $u=0$,
$$
P_k(0)=\left\{\begin{array}{ll}
0 & \mbox{if $1\le k\le N-1$,}\\
N!\,\a\, P_0(0)=N!\,\a & \mbox{if $k=N$.}
\end{array}\right.
$$
Choosing $\a=\k_N i^Nt$ and $u=0$ in~\refp{derivatives-exp-auN}, we
immediately complete the proof of Lemma~\ref{moments}.
\fin
%
\blm{
For $\a<0<\beta$, the following expansion holds as $\e\to 0^+$:
\beq\label{expansion-phi}
\int_{\a}^{\beta} e^{i\mu z}\,p(\e;-z)\,dz= 1+\k_N(i\mu)^N\e+o(\e).
\eeq
}
%
\dem
Performing a simple change of variables and using some asymptotics
of~\cite{arcsine}, we get
\beqa
\int_{\a}^{\beta} e^{i\mu z}\,p(\e;-z)\,dz
&=&
\int_{\a/\e^{1/N}}^{\beta/\e^{1/N}} e^{i\mu\e^{1/N}z}\,p(1;-z)\,dz
\;=\;\int_{-\infty}^{+\infty} e^{i\mu\e^{1/N}z}\,p(1;-z)\,dz+o(\e)
\\
&=&
\sum_{p=0}^{\infty} \frac{(i\mu)^p}{p!}\,\e^{p/N}
\int_{-\infty}^{+\infty} z^p\,p(1;-z)\,dz+o(\e).
\eeqa
Observing that $\int_{-\infty}^{+\infty} z^p\,p(1;-z)\,dz=\E_0(X(1)^p)$,
we immediately derive from Lemma~\ref{moments} the expansion~\refp{expansion-phi}.
\fin

\noindent\textbf{$\bullet$ Step 5}

Now, plugging~\refp{expansion-phi} into~\refp{function-phi}, it comes
\beqa
\f(\l,\mu,\nu;y)
&=&
\lim_{\e\to 0^+} \frac{1}{\e}\left[
\left(1-\left(\l-\k_N(i\nu)^N\right)\e+o(\e)\right)\left(1+\k_N(i\mu)^N\e+o(\e)\right)
\right.
\\
&&
\left.-\left(1+\k_N(i\nu)^N\e+o(\e)\right)\right]=-\l+\k_N(i\mu)^N.
\eeqa
Therefore,
$$
\Ex\!\!\left[e^{-\l\tau_a+i\mu X(\tau_a-)+i\nu \D_aX}\right]
=e^{i\mu x}-\left(\l-\k_N(i\mu)^N\right) \int_0^{\infty} e^{-\l t}
\,\E_x\!\!\left[e^{i\mu X(t)},M(t)\le a\right]dt.
$$
Writing finally
\beqa
\lqn{
\int_0^{\infty} e^{-\l t} \,\E_x\!\!\left[e^{i\mu X(t)},M(t)\le a\right]dt
}\\[-1ex]
&=&
\int_0^{\infty} e^{-\l t} \,\E_x\!\!\left[e^{i\mu X(t)}\right]dt
- \int_0^{\infty} e^{-\l t} \,\E_x\!\!\left[e^{i\mu X(t)},M(t)> a\right]dt
\\
&=&
\frac{1}{\l-\k_N(i\mu)^N}
-\int_0^{\infty} e^{-\l t} \,\E_x\!\!\left[e^{i\mu X(t)},M(t)> a\right]dt,
\eeqa
we obtain~\refp{LT-jump} by invoking the relationship~\refp{LFT-tau-Xtau-inter}
and by noting that the result does not depend on $\nu$ (and hence we can choose $\nu=0$).
\fin

Choosing $\mu=0$ or $\nu=0$, we obtain the corollary below.
%
\bco{\label{cor-cont}
We have, for $\Re(\l)>0$ and $\mu,\nu\in\R$,
$$
\Ex\!\!\left[e^{-\l\tau_a+i\mu X(\tau_a-)}\right]
=\Ex\!\!\left[e^{-\l\tau_a+i\mu X(\tau_a)}\right]
\quad\mbox{and}\quad
\Ex\!\!\left[e^{-\l\tau_a+i\nu \D_aX}\right]
=\Ex\!\!\left[e^{-\l\tau_a}\right]\!.
$$
}
%
From Corollary~\ref{cor-cont}, we expect that $\D_aX=X(\tau_a)-X(\tau_a-)=0$
in the following sense: $\forall \f\in\cS,$ $\E_0[\f(\D_aX)]=\f(0).$ This
provides a new argument in favor of continuity.


\subsection{Particular cases}

\bex{\textsl{Case $N=3$:}

\bitem
\item
In the case $\k_3=+1$, densities~\refp{density-tau_Xtau} write
$$
\P_x\{\tau_a^+\in dt,X(\tau_a^+)\in dz\}/dt\,dz
=\cJ_0(t;x-a)\,\d_a(z)\mbox{ for $x\le a$}
$$
and
$$
\P_x\{\tau_a^-\in dt,X(\tau_a^-)\in dz\}/dt\,dz
=\cK_0(t;x-a)\,\d_a(z)+\cK_1(t;x-a)\,\d'_a(z)
\mbox{ for $x\ge a$.}
$$
Here, we have $d_{11}=\t_2=\bar{\t}_1$, $d_{21}=\t_1=\bar{\t}_2$ and
\beqa
\cJ_0(t;\xi)
&=&
-\frac{\xi}{\pi t} \,\Im\bigg[e^{\frac{i\pi}{3}}\int_0^{+\infty}
\t_0A_0 \,e^{\t_0 e^{\frac{i\pi}{3}}\xi\l}\,e^{-t\l^3} \,d\l \bigg]
\\
&=&
-\frac{\xi}{\pi t} \,\Im\bigg[e^{\frac{i\pi}{3}}\int_0^{+\infty}
e^{e^{\frac{i\pi}{3}}\xi\l-t\l^3} \,d\l \bigg]
\\
&=&
-\frac{\xi}{\pi t} \int_0^{+\infty} e^{\frac12\xi\l-t\l^3}
\sin\bigg(\frac{\sqrt3}{2}\,\xi\l+\frac{\pi}{3}\bigg)\,d\l;
\eeqa
\beqa
\cK_0(t;\xi)
&=&
-\frac{\xi}{\pi t} \,\Im\bigg[e^{\frac{i\pi}{3}}\int_0^{+\infty}
\Big(\t_1B_1 \,e^{\t_1 e^{\frac{i\pi}{3}}\xi\l}
+\t_2B_2 \,e^{\t_2 e^{\frac{i\pi}{3}}\xi\l}\Big)\,e^{-t\l^3} \,d\l \bigg]
\\
&=&
-\frac{\xi}{\pi\sqrt3\,t} \,\Re\bigg[e^{\frac{i\pi}{3}}\int_0^{+\infty}
\Big(e^{-\xi\l}-e^{e^{-\frac{i\pi}{3}}\xi\l}\Big)\,e^{-t\l^3} \,d\l \bigg]
\\
&=&
-\frac{\xi}{\pi\sqrt3\,t} \int_0^{+\infty}
\bigg[\frac12\,e^{-\xi\l}-e^{\frac12\xi\l}
\cos\bigg(\frac{\sqrt3}{2}\,\xi\l-\frac{\pi}{3}\bigg)\bigg]\,e^{-t\l^3}\,d\l;
\eeqa
\beqa
\cK_1(t;\xi)
&=&
-\frac{3}{\pi} \,\Im\bigg[e^{-\frac{i\pi}{3}}\int_0^{+\infty}
\bigg(\bar{d}_{11}B_1 \,e^{\t_1 e^{\frac{i\pi}{3}}\xi\l}
+\bar{d}_{21}B_2 \,e^{\t_2 e^{\frac{i\pi}{3}}\xi\l}\bigg)
\,\l\,e^{-t\l^3} \,d\l \bigg]
\\
&=&
-\frac{\sqrt3}{\pi} \,\Re\bigg[e^{-\frac{i\pi}{3}}\int_0^{+\infty}
\bigg(e^{-\xi\l}-e^{e^{-\frac{i\pi}{3}}\xi\l}\bigg)\,\l\,e^{-t\l^3} \,d\l \bigg]
\\
&=&
-\frac{\sqrt3}{\pi} \int_0^{+\infty}
\bigg[\frac12\,e^{-\xi\l}-e^{\frac12\xi\l}
\cos\bigg(\frac{\sqrt3}{2}\,\xi\l+\frac{\pi}{3}\bigg)\bigg]\,\l\,e^{-t\l^3}\,d\l.
\eeqa

\item
In the case $\k_3=-1$, densities~\refp{density-tau_Xtau} write
$$
\P_x\{\tau_a^+\in dt,X(\tau_a^+)\in dz\}/dt\,dz
=\cJ_0(t;x-a)\,\d_a(z)+\cJ_1(t;x-a)\,\d'_a(z)
\mbox{ for $x\le a$}
$$
and
$$
\P_x\{\tau_a^-\in dt,X(\tau_a^-)\in dz\}/dt\,dz
=\cK_0(t;x-a)\,\d_a(z)\mbox{ for $x\ge a$.}
$$
In this case, we have $c_{01}=\bar{\t}_2=\t_0$, $c_{21}=\bar{\t}_0=\t_2$ and
\beqa
\cJ_0(t;\xi)
&=&
-\frac{\xi}{\pi\sqrt3\,t} \int_0^{+\infty}
\bigg[\frac12\,e^{\xi\l}+e^{-\frac12\xi\l}
\cos\bigg(\frac{\sqrt3}{2}\,\xi\l+\frac{\pi}{3}\bigg)\bigg]\,e^{-t\l^3}\,d\l
\\
\cJ_1(t;\xi)
&=&
-\frac{\sqrt3}{\pi} \int_0^{+\infty}
\bigg[\frac12\,e^{\xi\l}-e^{-\frac12\xi\l}
\cos\bigg(\frac{\sqrt3}{2}\,\xi\l-\frac{\pi}{3}\bigg)\bigg]\,\l\,e^{-t\l^3}\,d\l
\\
\cK_0(t;\xi)
&=&
-\frac{\xi}{\pi t} \int_0^{+\infty} e^{-\frac12\xi\l-t\l^3}
\sin\bigg(\frac{\sqrt3}{2}\,\xi\l-\frac{\pi}{3}\bigg)\,d\l.
\eeqa
Let us point out that the functions $\cJ_0$, $\cJ_1$, $\cK_0$ and $\cK_1$
may be expressed by means of Airy functions.
\eitem}

\bex{\textsl{Case $N=4$:} formulas~\refp{density-tau_Xtau} read here
$$
\P_x\{\tau_a^+\in dt,X(\tau_a^+)\in dz\}/dt\,dz
=\cJ_0(t;x-a)\,\d_a(z)+\cJ_1(t;x-a)\,\d'_a(z)
\mbox{ for $x\le a$}
$$
and
$$
\P_x\{\tau_a^-\in dt,X(\tau_a^-)\in dz\}/dt\,dz
=\cK_0(t;x-a)\,\d_a(z)+\cK_1(t;x-a)\,\d'_a(z)
\mbox{ for $x\ge a$.}
$$
We have
$c_{21}=\t_3=\bar{\t}_2$, $c_{31}=\t_2=\bar{\t}_3$,
$d_{01}=\t_1=\bar{\t}_0=-\bar{\t}_2$, $d_{11}=\t_0=\bar{\t}_1=-\bar{\t}_3$
and
\beqa
\cJ_0(t;\xi)
&=&
-\frac{\xi}{\pi t} \,\Im\bigg[e^{\frac{i\pi}{4}}\int_0^{+\infty}
\Big(\t_2A_2 \,e^{\t_2 e^{\frac{i\pi}{4}}\xi\l}+
\t_3A_3 \,e^{\t_3 e^{\frac{i\pi}{4}}\xi\l}\Big)\,e^{-t\l^4} \,d\l \bigg]
\\
&=&
\frac{\xi}{\pi\sqrt2\,t} \,\Re\bigg[e^{\frac{i\pi}{4}}\int_0^{+\infty}
(e^{\xi\l}-e^{i\xi\l})\,e^{-t\l^4}\,d\l \bigg]
\\
&=&
\frac{\xi}{2\pi t} \int_0^{+\infty} \Big[e^{\xi\l}-\sqrt2\,
\cos\Big(\xi\l+\frac{\pi}{4}\Big)\Big]\,e^{-t\l^4}\,d\l
\\
&=&
\frac{\xi}{2\pi t} \int_0^{+\infty} \Big[e^{\xi\l}
-\cos(\xi\l)+\sin(\xi\l)\Big]\,e^{-t\l^4}\,d\l;
\eeqa
\beqa
\cJ_1(t;\xi)
&=&
-\frac{4}{\pi} \,\Im\bigg[e^{-\frac{i\pi}{4}}\int_0^{+\infty}
\Big(\bar{c}_{21}A_2 \,e^{\t_2 e^{\frac{i\pi}{4}}\xi\l}+
\bar{c}_{31}A_3 \,e^{\t_3 e^{\frac{i\pi}{4}}\xi\l}\Big)
\,\l^2\,e^{-t\l^4} \,d\l \bigg]
\\
&=&
\frac{2\sqrt2}{\pi} \,\Re\bigg[e^{-\frac{i\pi}{4}}\int_0^{+\infty}
(e^{\xi\l}-e^{i\xi\l})\,\l^2\,e^{-t\l^4}\,d\l \bigg]
\\
&=&
-\frac{2}{\pi} \int_0^{+\infty} \Big[e^{\xi\l}-\sqrt2\,
\cos\Big(\xi\l-\frac{\pi}{4}\Big)\Big]\,\l^2\,e^{-t\l^4}\,d\l
\\
&=&
\frac{2}{\pi} \int_0^{+\infty} \Big[\cos(\xi\l)+\sin(\xi\l)
-e^{\xi\l}\Big]\,\l^2\,e^{-t\l^4}\,d\l
\eeqa
and similarly
\beqa
\cK_0(t;\xi)
&=&
\frac{\xi}{2\pi t} \int_0^{+\infty} \Big[\cos(\xi\l)+\sin(\xi\l)
-e^{-\xi\l}\Big]\,e^{-t\l^4}\,d\l;
\\
\cK_1(t;\xi)
&=&
\frac{2}{\pi} \int_0^{+\infty} \Big[e^{-\xi\l}-\cos(\xi\l)+
\sin(\xi\l)\Big]\,\l^2\,e^{-t\l^4}\,d\l.
\eeqa
We retrieve formulas~(8.17) and~(8.18) of~\cite{nish2}.
}

\subsection{Boundary value problem}

We end up this work by exhibiting a boundary value problem satisfied by the
Laplace-Fourier transform
$U(x)=\Ex\Big[e^{-\l\tau_a^++i\mu X(\tau_a^+)}\Big]$, $x\in(-\infty,a)$.
%
\bpr{
The function $U$ satisfies the differential equation
\beq\label{ODE-U}
\cD_xU(x)=\l U(x)\quad\mbox{for $x\in(-\infty,a)$}
\eeq
together with the conditions
\beq\label{boundary-cond-U}
U^{(l)}(a^-)=(i\mu)^le^{i\mu a}\quad\mbox{for $0\le l\le \card J-1$.}
\eeq
}
%
\dem
The differential equation~\refp{ODE-U} is readily obtained by
differentiating~\refp{LFT-tau-Xtau} with respect to $x$.
Let us derive the boundary conditions~\refp{boundary-cond-U}:
by~\refp{LFT-tau-Xtau},
\beqa
U^{(l)}(a^-)
&=&
\l^{l/N}\prod_{j\in J}\left(1-\frac{i\mu}{\sqrt[N]{\l}}\,\bar{\t}_j\right)
\Bigg(\sum_{j\in J} \frac{\t_j^l A_j}{1-\frac{i\mu}{\sqrt[N]{\l}}
\,\bar{\t}_j}\Bigg) e^{i\mu a}.
\eeqa
By~\refp{expansion} we see that
$$
\sum_{j\in J} \frac{\t_j^l A_j}{1-\frac{i\mu}{\sqrt[N]{\l}}\,\bar{\t}_j}=
\frac{(i\mu)^l}{\l\prod_{j\in J}
\left(1-\frac{i\mu}{\sqrt[N]{\l}}\,\bar{\t}_j\right)}
$$
which proves Condition~\refp{boundary-cond-U}.
\fin

We also refer the reader to~\cite{nish3} for a very detailed account on
PDE's with various boundary conditions and their connections with different
kinds of absorbed/killed pseudo-processes.

\textbf{Acknowledgement.} This work has been inspired to the author while
sojourning at the university ``La Sapienza'' (Roma, Italy) where many
discussions with Pr.~E.~Orsingher and Dr.~L.~Beghin were very fruitful.


\end{document}